\documentclass[11pt]{amsart}
\usepackage{amsmath}
\usepackage{amsthm}
\usepackage{epsfig}
\usepackage{amsfonts}

\swapnumbers \theoremstyle{plain}

\newtheorem{thm}{Theorem}[section]

\newtheorem{lem}[thm]{Lemma}
\newtheorem{cor}[thm]{Corollary}

\newtheorem{prop}[thm]{Proposition}

\newtheorem{prob}[thm]{Problem}

\theoremstyle{definition}

\newtheorem{ex}[thm]{Example}

\theoremstyle{remark}
\newtheorem*{rem}{Remark}

\newtheorem*{claim}{Claim}

\swapnumbers

\newcommand{\T}{\mathcal T}

\renewcommand{\P}{{\mathcal P}}
\renewcommand{\S}{{\mathcal S}}
\newcommand{\bdy}{\partial}

\title{Decision problems in the space of  Dehn fillings }
\author{William Jaco }
\author{Eric Sedgwick}
\date{\today}
\address{mathematics department, oklahoma state university, stillwater, ok, 74078}
\email{jaco@math.okstate.edu}
\email{sedgwic@math.okstate.edu}
\thanks{The first author was partially supported by NSF Grant DMS9704833,
and the Grayce B. Kerr Foundation, along with sabbatical support
from Oklahoma State University, University of California - Davis,
California Institute of Technology, and University of Texas -
Austin.}

\usepackage{amssymb}

\begin{document}

\begin{abstract}
In this paper, we use normal surface theory to study Dehn filling
on a knot-manifold.  First, it is shown that there is a finite
computable set of slopes on the boundary of a knot-manifold that
bound normal and almost normal surfaces in a one-vertex
triangulation of that knot-manifold. This is combined with
existence theorems for normal and almost normal surfaces to
construct algorithms to determine precisely which manifolds
obtained by Dehn filling: 1) are reducible, 2) contain two--sided
incompressible surfaces, 3) are Haken, 4) fiber over $S^1$, 5) are
the 3--sphere, and 6) are a lens space.  Each of these algorithms
is a finite computation.

Moreover, in the case of essential surfaces, we show that the
topology of the filled manifolds is strongly reflected in the
triangulation of the knot-manifold.  If a filled manifold contains
an essential surface then the knot-manifold contains an essential
{\it vertex solution} that caps off to an essential surface of the
same type in the filled manifold.  (Vertex solutions are the
premier class of normal surface and are computable.)
\end{abstract}

\maketitle

\section{Introduction}

A compact, connected, orientable 3--manifold with connected
boundary a torus is called a {\it knot-manifold}. Dehn filling is
a method of obtaining closed 3--manifolds from a knot-manifold. It
is a special case of a general construction from which one can
obtain all closed, orientable, 3--manifolds \cite{lick, wall}.
Specifically, if $X$ is a knot-manifold, we call an isotopy class
of a simple closed curve in $\bdy X$ a {\it slope}. If $\alpha$ is
a slope, the {\it Dehn filling of $X$ along $\alpha$}, denoted
$X(\alpha)$, is the closed, orientable 3--manifold obtained from
$X$ by attaching a solid torus $V_{\alpha}$ to $X$ via a
homeomorphism from $\bdy X$ to $\bdy V_{\alpha}$ which takes a
simple closed curve of slope $\alpha$ to the meridian of
$V_{\alpha}$; i.~e., to an essential curve in $\bdy V_{\alpha}$
that bounds a disk in $V_{\alpha}$. The homeomorphism type of
$X(\alpha)$ is completely determined by the identification of the
slope $\alpha$ to a meridian of $V_{\alpha}$.

If $X$ is a knot-manifold and we select a homology basis, say
$\mu$, $\lambda$, for $H_{1}(\bdy X)$, then each slope $\alpha$
can be written $\alpha = p\mu + q\lambda$ where p and q are
integers. Hence, if we include $\infty$ and forget orientation
(sign) of a homology class, the slope $\alpha$ is uniquely
associated with a rational number p/q with $1/0$ associated with
$\infty$. Hence, for a given knot-manifold $X$ we obtain a family
of closed, orientable manifolds $X(\alpha), \alpha \in \mathbb{Q}
\cup \{\infty\}$. The collection of such manifolds is called the
{\it space of Dehn fillings on $X$}.

A great deal of work has been done to understand the manifolds in
the space of Dehn fillings on a knot-manifold. In particular, for
a hyperbolic knot-manifold $X$, a knot-manifold whose interior
admits a complete Riemannian metric of constant sectional
curvature $-1$, it has been shown \cite{thurston:1982} that
$X(\alpha)$ is hyperbolic for all but finitely many slopes
$\alpha$. For the past decade some of the most interesting work in
low-dimensional topology has been toward understanding exceptions
to $X(\alpha)$ being hyperbolic. In this sense, the exceptions
include the possibilities that $X(\alpha)$ is reducible, toroidal,
or a lens space.

If $\alpha$ and $\beta$ are slopes, we let $\Delta (\alpha,\beta)$
denote the absolute value of the homology intersection between
$\alpha$ and $\beta$ and call $\Delta (\alpha,\beta)$ the {\it
distance between $\alpha$ and $\beta$}.  If for some homology
basis of $H_{1}(\bdy X)$ we have $\alpha = p/q$ and $\beta = r/s$,
then $\Delta (\alpha,\beta) = |ps - qr|$. Now, if $X$ is a
hyperbolic knot-manifold and $\alpha$ and $\beta$ are exceptional
slopes, then in many situations bounds can be placed on $\Delta
(\alpha,\beta)$ and thereby one obtains bounds on the numbers of
exceptional Dehn fillings on $X$ \cite{go1}. The remarkable and
very satisfactory consequence of these methods is that the bounds
obtained are global; they do not depend on $X$. For example, it is
conjectured that for $X$ a hyperbolic knot-manifold and $X$ not
one of a finite number of exceptions formed by Dehn filling on a
component of the Whitehead link in $S^3$, then  $\Delta
(\alpha,\beta) \leq 5$ if $\alpha$ and $\beta$ are exceptional
slopes \cite{go1}. It is known that  $\Delta (\alpha,\beta) \leq
1$ if $X(\alpha)$ and $X(\beta)$ are reducible \cite{g-l2,bz1};
$\Delta (\alpha,\beta) \leq 5$ if $X(\alpha)$ and $X(\beta)$ have
finite fundamental group \cite{bz1}; and for all but the
aforementioned exceptions on $X$, $\Delta (\alpha,\beta) \leq 5$
if $X(\alpha)$ and $X(\beta)$ are toroidal \cite{go}. Results for
mixed outcomes of exceptional Dehn fillings are given in
\cite{go1}. For a knot-manifold embedded in $S^3$, a preferred
basis for $H_{1}(\bdy X)$ is the unique meridian and longitude
pair, $\mu$ and $\lambda$, respectively.

Our work addresses most of these same issues about Dehn filling
but from a different point of view. Namely, given a knot-manifold
$X$, we are interested in determining precisely those slopes
$\alpha$ on $\bdy X$ for which Dehn filling leads to
``interesting" phenomena for $X(\alpha)$. In particular, we
consider for precisely what slopes $\alpha$ is $X(\alpha)$
reducible; is $X(\alpha)$ toroidal; does $X(\alpha)$ contain an
embedded, incompressible, two-sided surface; is $X(\alpha)$ a
Haken-manifold; does $X(\alpha)$ fiber over $S^{1}$; is
$X(\alpha)$ homeomorphic with $S^3$; and is $X(\alpha)$ a lens
space. Recall that given a 3--manifold $M$ there are algorithms to
answer each of these questions regarding $M$. Namely, given a
compact 3--manifold $M$, it can be decided if $M$ is reducible
\cite{rubinstein,thompson,j-t,j-re}; it can be decided if $M$ is
toroidal \cite{haken,j-o}; it can be decided if $M$ contains an
embedded, incompressible two-sided surface \cite{haken,j-o,j-t};
it can be decided if $M$ fibers over $S^1$ \cite{jac}; it can be
decided if $M$ is homeomorphic with $S^3$
\cite{rubinstein,thompson}; and it can be decided if $M$ is a lens
space \cite{rubinstein}. Now, one might think that the existence
of these algorithms will solve our problem; however, for a given
knot-manifold $X$, the family of manifolds in the Dehn filling
space of $X$ is infinite. Hence, we have the situation that
knowing that there is a manifold in the space of Dehn fillings of
$X$ that is of interest, then we can find one (our problem is
recursively enumerable); however, without {\it a priori}
information, these algorithms, alone, will not necessarily
determine if there is an interesting manifold, let alone determine
all slopes for which such interesting phenomena occur. In this
paper we provide the additional ingredients and algorithms to
determine precisely the slopes, or manifolds in the space of Dehn
fillings of $X$, that exhibit the various ``interesting" phenomena
mentioned above.

We will assume 3-manifolds are given via triangulations or cell
subdivisions.  In most settings we use either one-vertex
triangulations of the manifolds under considerations or at least a
triangulation that restricts to a one-vertex triangulation on each
torus component of the boundary. The existence of such
triangulations is straight forward and discussion of these and
other useful triangulation environments are given in \cite{j-r4}.
We use normal and almost normal surface theory for these
triangulations.

The study of Dehn fillings has exhibited strong relationships
between the topology of $X$ and those manifolds in the space of
fillings of $X$. Our methods re-enforce this relationship in a
remarkable way. Given a knot-manifold $X$  via a triangulation
$\T$ that restricts to a one-vertex triangulation on $\bdy X$, we
use the methods of \cite{j-r3, j-r4} to extend $\T$ to a
triangulation of $X(\alpha)$; that is, for each slope $\alpha$, we
construct a triangulation $\T (\alpha)$ of $X(\alpha)$ that
restricts to $\T$ on $X$. Furthermore, the triangulation $\T
(\alpha)$ restricts to a well understood one-vertex triangulation
of $V_{\alpha}$, the attached solid torus. Each of the problems we
consider is to determine precisely the slopes $\alpha$ for which a
certain type of surface exists in the manifold $X(\alpha)$. For
example, reducibility is the existence of an embedded 2--sphere
that does not bound  a 3--cell; and to determine if $X(\alpha)$ is
$S^3$ or a lens space is to find a genus zero or genus one
Heegaard surface, respectively. Normal and almost normal surface
theory provide a parameterization of ``interesting" surfaces by
rational points in a computable, compact, convex, linear cell in
$\mathbb R^{n}$, the {\it projective solution space}. If $X$ is a
knot-manifold and $\T$ is a triangulation of $X$, we denote the
projective solution space of $X$ with respect to $\T$ by $\P
(X,\T)$. In this situation, if $S$ is a properly embedded, normal
surface in $(X,\T)$, then either $\bdy S =\emptyset$ or $\bdy S
\neq \emptyset$ and $\bdy S$ is a collection of pairwise disjoint,
normal curves in $\bdy X$. If $\bdy S = \emptyset$ or is a
collection of trivial and, hence, vertex-linking curves, then for
any slope $\alpha$, $S$ determines a unique normal surface
$S(\alpha)$ in $X(\alpha)$ ($S(\alpha)$ is obtained from $S$ by
capping off $\bdy S$ with copies of the vertex-linking normal
disks in the special triangulation of $V_{\alpha}$ determined by
$\T (\alpha)$). If $\bdy S$ contains a nontrivial component and
determines a unique boundary slope $\alpha$, then $S$ determines a
unique normal surface $S(\alpha)$ in $X(\alpha)$ just for the
slope $\alpha$ ($S(\alpha)$ is obtained from $S$ by capping off
$\bdy S$ with copies of the vertex-linking normal disks and
meridional normal disks in the special triangulation of
$V_{\alpha}$ determined by $\T (\alpha)$). We show a Dehn filling
$X(\alpha)$ will contain one of our ``interesting" surfaces,
listed above, if and only if there is a normal or almost normal
surface $S$ in $(X,\T)$ whose projective class is a vertex
solution of $\P (X,\T)$ and $S(\alpha)$ has the same interesting
property in $X(\alpha)$. Hence, for any triangulation $\T$ of $X$
that restricts to a one-vertex triangulation on $\bdy X$, the
normal and almost normal surfaces in $(X,\T)$ whose projective
classes in $\P (X,\T)$ are vertex solutions completely determine
for all slopes $\alpha$ whether the manifold $X(\alpha)$ is
reducible, toroidal, contains an embedded, incompressible,
two-sided surface, fibers over $S^1$, is $S^3$, or is a lens
space. The vertex solutions of $\P (X,\T)$ form a finite,
computable set. It is this set which plays the fundamental role in
most of our algorithms.

In Section 2, we recall material from normal and almost normal
surface theory. We limit this to material that is directly
relevant to this paper and assume the reader has some familiarity
with this theory. More sweeping introductions from our point of
view may be found in \cite{j-r1, j-t, j-re}.

In Section 3, we introduce one of the fundamental features of
using one-vertex triangulations: the relationship between
interesting slopes and normal and almost normal surfaces. Also, we
provide techniques to compute interesting slopes. We compute the
projective solution space of normal curves in a one-vertex
triangulation of a torus. It is represented as the standard
2--simplex, $\Delta = \{(x_{1},x_{2},x_{3}) \in \mathbb R^{3} :
\sum_{i=1}^{3} x_{i} = 1, x_{i} \geq 0\}$. The rational points of
$\Delta$ represent projective classes of families of embedded,
normal curves in the torus; the barycenter of $\Delta$ represents
the trivial, vertex linking family and the edges (some $x_{i} =
0$) represent the various slopes of the families of embedded
curves having nontrivial components. Now, given a triangulation
$\T$ of a knot-manifold $X$ that restricts to a one-vertex
triangulation on $\bdy X$, we call a slope $\alpha$ a {\it
boundary slope} if there is a normal or almost normal surface $S$
properly embedded in $X$ and a component of $\bdy S$ represents a
curve of slope $\alpha$. We prove that for such a triangulation
there are only finitely many slopes $\alpha$ which are boundary
slopes; furthermore, we do this by showing that the boundary
slopes are completely determined by the boundary slopes of normal
or almost normal surfaces in $(X,\T)$ whose projective class is a
vertex solution of $\P (X,\T)$, a finite set. This result
generalizes the results of \cite{hatcher} (using similar
techniques) and gives a new proof of the main theorem in
\cite{hatcher} that there are only finitely many boundary slopes
for embedded, incompressible and $\partial$-incompressible
surfaces in $X$. However, a distinguishing feature of our work is
a means to actually compute precisely the relevant boundary slopes
from the triangulation $\T$ of $X$.

In Section 4 we generalize the one-vertex triangulations of solid
tori introduced in \cite{j-r3} to one-vertex triangulations called
{\it layered triangulations}. We analyze the embedded, planar,
normal surfaces in these layered triangulations, classifying such
surfaces and obtaining lower bounds for their weights (the {\it
weight} of a normal surface is the cardinality of its intersection
with the one-skeleton of the triangulation). These results lead to
the special triangulations we use when studying Dehn fillings and
enable us to relate the existence of interesting normal surfaces
in the manifolds in the space of Dehn fillings to interesting
surfaces in the original knot-manifold.

In Sections 5 and 6 we consider the central problems of this
paper: given a  knot-manifold $X$, to determine precisely those
slopes $\alpha$ on $\bdy X$ for which Dehn filling leads to
``interesting" phenomena for $X(\alpha)$. We divide this work into
two parts. In Section 5, we look at phenomena associated with
embedded {\it essential} surfaces and in Section 6, we look at
phenomena associated with {\it Heegaard} surfaces. We have
organized the presentation so that following proofs of the
existence of certain algorithms, we give step by step outlines of
the algorithms.

In Section 5, Theorem \ref{T-essinx} provides one of the major
ingredients for our algorithms. It ties the topology of a
knot-manifold $X$ quite tightly with that of the manifolds
obtained by Dehn filling on $X$. For a triangulation $\T$ that
restricts to a one-vertex triangulation of $\bdy X$, Theorem
\ref{T-essinx} gives that if $X(\alpha)$ is reducible, then a
vertex-solution $S$ of $\mathcal P$($X$,$\mathcal T$) must be
either an embedded, essential 2--sphere or planar surface and
$S(\alpha)$ is an embedded, essential 2--sphere in $X(\alpha)$; if
$X(\alpha)$ contains an embedded, incompressible, two-sided
surface, then a vertex-solution $F$ of $\mathcal P$($X$,$\mathcal
T$) must be an embedded, essential, non-planar surface and
$F(\alpha)$ is an embedded, incompressible, two-sided surface in
$X(\alpha)$.  In the latter case, if $X(\alpha)$ contains an
embedded, incompressible torus, then a vertex-solution $T$ of
$\mathcal{P}(X,\mathcal{T})$ must be an embedded, essential torus
or punctured-torus and $T(\alpha)$ is an embedded, incompressible
torus in $X(\alpha)$; and if $X(\alpha)$ fibers over $S^{1}$, then
a vertex-solution $F$ of $\mathcal{P}(X,\mathcal{T})$ has the
property that $F(\alpha)$ is a fiber in a fibration of $X(\alpha)$
over $S^{1}$.

Theorem \ref{T-reducibleslopes} and Theorem
\ref{T-irreducibleslopes} show that there exists an algorithm to
determine precisely those manifolds in the space of Dehn fillings
of a knot-manifold $X$ that are reducible. In particular we have

\vspace{.2 in}

\noindent {\bf Algorithm R.} {\it Given a knot-manifold $X$,
determine precisely those slopes $\alpha$ for which the Dehn
filling $X(\alpha)$ is reducible.}

\vspace{.2 in} Theorem \ref{T-surfacesincompress} and Theorem
\ref{T-surfacescompress} show that there exists an algorithm to
determine precisely those manifolds which contain an embedded,
incompressible, two-sided surface in the space of Dehn fillings of
a knot-manifold $X$. Of particular importance in the proof of this
theorem, and of independent interest, is Lemma \ref{closedess}
which provides an algorithm to determine for a {\it given} closed,
two-sided, normal surface $S$ in $(X,\T)$ precisely those slopes
$\alpha$ for which $S$ is incompressible in $X(\alpha)$. The proof
of this results investigates when one can determine if a
3--manifold contains an embedded, essential punctured disk. It
uses a new estimate for curve length of the boundary of a normal
surface discovered by Jaco and Rubinstein; we call this the ALE,
average length estimate. It is used in \cite{j-ru-sedg} to give
algorithms for the existence of planar surfaces and their
relationship to the Word Problem for 3--manifold groups. Our work
here provides two algorithms: one that determines if a given
closed surface is incompressible in a Dehn filling, Algorithm I,
and one that determines precisely those slopes $\alpha$ for which
the associated Dehn filling contains an embedded, incompressible,
two-sided surface, Algorithm S.

\vspace{.2 in}

\noindent {\bf Algorithm I.} {\it Suppose $X$ is a knot-manifold
with a triangulation $\T$ which restricts to a one-vertex
triangulation on $\partial X$. Given an embedded, two-sided,
closed, normal surface in $(X,\T)$, determine precisely those
slopes $\alpha$ for which the surface compresses in the Dehn
filling $X(\alpha)$.}

\vspace{.2 in}

\noindent {\bf Algorithm S.} {\it Given a knot-manifold $X$,
determine precisely those slopes $\alpha$ for which the Dehn
filling $X(\alpha)$ contains an embedded, incompressible,
two-sided surface.}

\vspace{.2 in} We use Algorithm R and Algorithm S to give an
algorithm to determine precisely those slopes for which the
associated Dehn filling is a Haken-manifold, Algorithm H.

\vspace{.2 in} \noindent{\bf Algorithm H.} {\it Given a
knot-manifold $X$, determine precisely those slopes $\alpha$ for
which the Dehn filling $X(\alpha)$ is a Haken-manifold.}

\vspace{.2 in}

At particular points in the application of Algorithm S, one may
consider the alternative questions as to those slopes $\alpha$ for
which the Dehn filling $X(\alpha)$ is either toroidal, the
existence of an embedded, incompressible torus, or fibers over
$S^1$, the existence of an embedded, incompressible surface that
is a fiber in such a fibration.

Finally, in Section 6, we apply our techniques to similar
considerations for Heegaard surfaces. We use almost normal
surfaces introduced by H.~Rubinstein \cite{rubinstein} and thin
position introduced by D.~Gabai \cite{gabai} as presented in the
papers of Rubinstein \cite{rubinstein} and A.~Thompson
\cite{thompson}. The two main results of this section are given in
Theorem 6.4 and Theorem 6.7 which provide algorithms to determine
for a given knot-manifold $X$ precisely those slopes $\alpha$ for
which the Dehn filling $X(\alpha)$ is either $S^3$ or a lens
space, respectively.

\vspace{.2in}

\noindent {\bf Algorithm $\S$.} {\it Given a knot-manifold $X$,
determine precisely those slopes $\alpha$ for which $X(\alpha)$ is
the 3--sphere.}

\vspace{.2in}

\noindent {\bf Algorithm L.} {\it Given a knot-manifold $X$,
determine precisely those slopes $\alpha$ for which $X(\alpha)$ is
a lens space.}

\vspace{.2in}

The authors wish to thank J.~Hyam Rubinstein, whose collaborations
with the first author have lead to many useful ideas and tools
used in this work, as well as the members of the topology group at
Oklahoma State University, who contributed to this paper through
many useful formal and informal discussions. The first author also
acknowledges, with thanks, the support of the Departments of
Mathematics, and especially the topology groups, at University of
California--Davis, California Institute of Technology, and
University of Texas--Austin during his sabbatical leave from
Oklahoma State University.


\section{Normal curves and normal surfaces}
\label{s-normal}

\begin{figure}[h]
{\epsfxsize = 1.5 in \centerline{\epsfbox{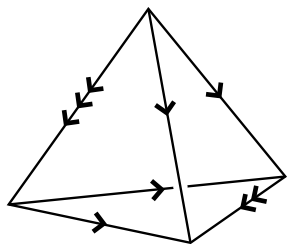}} } \caption{A
one tetrahedron triangulation of $S^3$.} \label{f-s3}
\end{figure}

\begin{figure}[h]
{\epsfxsize = 1.5 in \centerline{\epsfbox{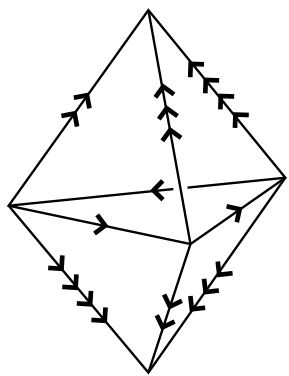}} } \caption{A
two tetrahedra triangulation of  $\mathbb RP^3$.} \label{f-rp3}
\end{figure}

Throughout this paper a 3--manifold will be given via a
triangulation, where a {\it triangulation } $\T$ of a 3--manifold
$M$ is a pairwise disjoint collection of tetrahedra,
$\boldsymbol{\Delta} = \{ \Delta_1, \dots, \Delta_t \}$, along
with a family $\Phi$ of face identifications having $M$ the
underlying point set of the identification space
$\boldsymbol{\Delta}/\Phi$. Under this definition the tetrahedra
may not be embedded in $M$ and two distinct tetrahedra may meet in
more than a face of each.  Figure \ref{f-s3} shows a one
tetrahedron triangulation of the 3--sphere, $S^3$, and Figure
\ref{f-rp3} shows the two tetrahedra triangulation of the familiar
lens space presentation of real projective 3--space, $\mathbb
RP^3$.

Triangulations of surfaces are considered in the same generality;
that is, a triangulation $\T$ of a surface $S$ is a pairwise
disjoint collection of triangles $\Lambda = \{ \lambda_1,
\dots,\lambda_s \}$, along with  a family $\Psi$ of edge
identifications having $S$ the underlying point set of the
quotient space $\Lambda/\Psi$. Figure \ref{f-torus} shows a
one-vertex triangulation of the torus $S^1 \times S^1$.

\begin{figure}[h]
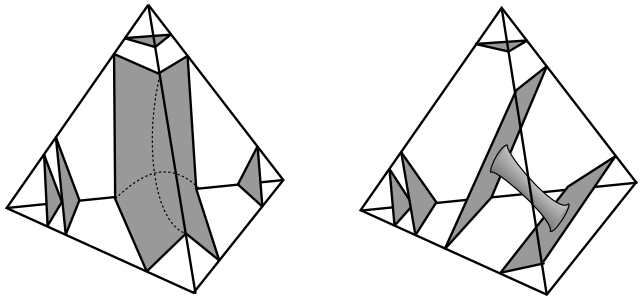

{\epsfxsize = 1.5 in \centerline{
  \epsfxsize = 1.5 in \epsfbox{F-OCT.AI} \qquad
  \epsfxsize = 1.5 in \epsfbox{F-TUBES.AI} \qquad
}} \caption{Exceptional pieces - an octagon and a tube.}
\label{f-almostnormal}
\end{figure}

We shall assume the reader has a basic understanding of normal
surface theory as well as the application of this theory to curves
in 2--manifolds.  The references \cite{j-r1} and \cite{j-t} are
sources to review normal surface theory.  We also use the concept
of an almost normal surface introduced by H. Rubinstein in
\cite{rubinstein}. If $\T$ is a triangulation of the 3--manifold
$M$, a surface $F$ is {\it almost normal} (with respect to $\T$)
if $F$ meets each tetrahedron of $\T$ in a collection of normal
triangles and quadrilaterals and intersects one of the tetrahedra
in one component that is either a normal octagon or a normal tube
and possibly some normal triangles.  See Figure
\ref{f-almostnormal}.  In Section \ref{s-heegaard} we prove the
existence of almost normal surfaces using octagons only.  Note
however, that our restrictions on the slopes bounding almost
normal surfaces developed in Section \ref{s-triang} also apply to
almost normal surfaces possessing tubes.

If $t$ is the number of tetrahedra in $\T$, then a normal isotopy
class of a normal surface has a parameterization as an $n$--tuple
of non-negative integers $(x_1,\dots,x_n)$ in $\mathbb R^n
~(n=7t)$, where $x_i$ is the number of elementary triangles and
quadrilaterals of type $i$.  Similarly,  there is a
parameterization of the normal isotopy classes of almost normal
surfaces, but in this case $n$ is larger as there are 3 normal
octagon types and 25 normal tube types in each tetrahedron.

Associated with the triangulation $\T$ is a system of linear
equations.  Non-negative integer solutions to this system give the
parameterization of the normal isotopy classes of normal and
almost normal surfaces.  We add the equation $\sum_{i=1}^n x_i =
1$ along with the condition $x_i \geq 0, \forall i$ and obtain a
compact, convex linear cell.  The rational points in this cell
correspond to projective classes of normal isotopy classes of
normal and almost normal surfaces in $(M,\T)$. We denote this
compact, convex linear cell by $\P(M,\T)$ and call it the {\it
projective solution space} (of $(M,\T)$).

If $S$ is a normal or almost normal surface in $M$ we do not
distinguish and let $S$ denote the surface $S$, its normal isotopy
class, and its representation as an $n$--tuple in $\mathbb R^n$.
We denote the projective class of $S$ by $\bar{S} \in \P(M,\T)$.
The {\it carrier} of a normal surface $S$ is the unique minimal
face of $\P(M,\T)$ that contains $\bar{S}$.  Two normal or almost
normal surfaces $S$ and $S'$ are {\it compatible} if and only if
each component of $S \cap S'$ is an embedded regular curve.
Compatibility is equivalent to the normal sum $S+S'$ being
defined. If $S$ and $S'$ are embedded normal or almost normal
surfaces then $S$ and $S'$ are compatible if and only if they do
meet in a tetrahedron in distinct normal quadrilateral types, a
quadrilateral and an exceptional piece, or two exceptional pieces,
where an exceptional piece is either an octagon or a tube.    If
$S$ is an embedded normal surface, then every normal surface with
projective class in the carrier of $S$ is embedded and any two
such normal surfaces are compatible.

W. Haken has observed \cite{haken} that there is a finite set of
embedded normal surfaces $F_1,\dots,F_N$ so that any normal
surface $S$ can be written as a non-negative integer linear
combination of the $F_i$'s; i.e. $$S = \sum_1^N n_i F_i,
\thickspace \text{ each $n_i$ is a non-negative integer.}$$ There
is a unique minimal such set, called the set of {\it fundamental
surfaces}.  A surface is fundamental if it cannot be written as a
non-trivial sum of surfaces. Among these fundamental surfaces is
an important set that have projective classes at the vertices of
$\P(M,\T)$. These latter surfaces are called {\it vertex
solutions}. A surface is a vertex solution if no multiple of the
surface can be written as a non-trivial sum of distinct surfaces.
Note that the sum notation is used for both normal (or geometric)
sum as well as coordinate-wise addition of $n$--tuples in $\mathbb
R^n$.

We remind the reader that when normal surface theory is applied to
curves in 2--manifolds; then every solution is realizable as an
embedded family of properly embedded arcs and simple closed
curves; i.e. there is always a unique embedded representative for
a solution. This is not the situation for normal surfaces in
3--manifolds; and solutions that do not have embedded
representatives (no realizable solutions) are not understood.  In
this paper we work only with embedded families of curves in
2--manifolds and with embedded surfaces in 3--manifolds.

For normal curves and normal surfaces there is a notion of
complexity analogous to geodesic curves and least are surfaces. If
$\T$ is a triangulation of the surface $S$ and $C$ a family of
normal curves in $S$ (with respect to $\T$) then we define the
{\it
  length of} $C$, written $L(C)$ to be the number of times $C$ meets
the 1--skeleton of $\T$; $$L(C) = | C \cap \T^{(1)}|.$$

Similarly, if $\T$ is a triangulation of the 3--manifold $M$ and
$S$ is a normal surface or almost normal surface in $M$, then we
define the {\it weight of} $S$, written $wt(S)$, to be the number
of times $S$ meets the 1--skeleton of $\T$; $$wt(S) = | S \cap
\T^{(1)}|.$$

If $S$ and $S'$ are embedded compatible normal surfaces, then the
normal sum $S+S'$ is defined and is a normal surface and we have:

\begin{enumerate}
\item If $S$ corresponds to the $n$--tuple $(x_1,\dots,x_n)$ and $S'$
  corresponds to the $n$--tuple $(x_1',\dots.x_n')$; then $S+S'$
  corresponds to the $n$--tuple $(x_1 + x_1',\dots,x_n+x_n')$.
\item $\chi(S + S') = \chi (S) + \chi(S')$, where $\chi$ is the Euler
  characteristic.
\item $wt(S+S')=wt(S) + wt(S')$.
\item $L(\bdy (S+ S')) = L(\bdy S) + L(\bdy S')$.
\end{enumerate}

The properties outlined in this section demonstrate that there is
a nice theory of computation using normal (almost normal)
surfaces. However, these computations are useful only if there
exist interesting surfaces with normal or almost normal
representatives.  In most situations, this is the case.  If $M$ is
an irreducible 3--manifold then every essential surface has a
normal representative in any triangulation of $M$, and every
strongly irreducible Heegaard surface has an almost normal
representative in any triangulation of $M$.

Unfortunately, when $M$ is a reducible manifold it may be
necessary to alter an essential surface before finding a normal
representative. Suppose $S$ is a surface properly embedded in the
3--manifold $M$ and $D'$ is a disk embedded in $M$ with $D' \cap S
= \partial{D'}$. Furthermore, suppose $\partial{D'}$ bounds a disk
$D \subset S$. Then $S' = (S \setminus D) \cup D'$ is a surface
topologically equivalent to $S$.  We say $S'$ is obtained from $S$
by a {\it disk-swap}. The two surfaces $S$ and $S'$ are said to be
{\it equivalent} (in $M$) if and only is there is a sequence $S =
S_{1}, \ldots, S_{n} = S'$ with $S = S_{1}$ and $S' = S_{n}$ where
$S_{i+1}$ is obtained from $S_{i}$ by a disk swap and/or isotopy.
Hence, if two surfaces $S$ and $S'$ are isotopic, then they are
equivalent. Equivalent and isotopic are the same when the ambient
manifold, $M$, is irreducible. The concept of ``disk-swapping"
applies to ``$\bdy$-compressing disks" as well and is a necessary
extension of this concept in the case that the manifold $M$ has
boundary and the surfaces in question are
$\partial$--incompressible. Note that any surface that is
equivalent to an incompressible and $\bdy$--incompressible surface
is also incompressible and $\bdy$--incompressible.

Let $S$ be a normal (or almost normal) surface in $(M,\T)$.  Then
$S$ is {\it least weight} if every normal (almost normal) surface
$S'$  that is equivalent to $S$ in $M$, we have that $wt(S) \leq
wt(S')$.

We now list the existence results mentioned above.  The first is
known from the work of H. Kneser \cite{kneser}.

\begin{thm}
Let $M$ be a 3--manifold.  If there is a 2--sphere embedded in $M$
that does not bound a 3--cell in $M$, then for any triangulation
of $M$ there is a  normal 2--sphere embedded in $M$ that does not
bound a 3--cell in $M$.
\end{thm}

The next theorem is from the work of W. Haken \cite{haken}.

\begin{thm}
Let $M$ be a 3--manifold.  If there is an incompressible and
$\bdy$--incompressible surface $S$ embedded in $M$, then for any
triangulation of $M$ there is a normal surface $S'$ embedded in
$M$ that is equivalent to $S$.
\end{thm}

Another reference where details can be found for the proofs of
these results is \cite{j-r1}.  Finally, we note M. Stocking's
result \cite{stocking} for Heegaard surfaces.

\begin{thm}
Let $M$ be an irreducible 3--manifold.  If there is a non-trivial
genus $g$ strongly irreducible Heegaard splitting of $M$, then for
any triangulation of  $M$ there is an almost normal genus $g$
surface isotopic to the Heegaard surface.
\end{thm}


\section{One-vertex triangulations and boundary slopes}
\label{s-triang}

\subsection{One-vertex triangulations.}

Computations in normal curve and normal surface theory can often
be simplified by selecting a special triangulation, in particular,
by choosing a triangulation with a minimal number of top
dimensional simplices. For surfaces with non-positive Euler
Characteristic, a minimum triangulation (a triangulation with the
minimal number of faces) requires a {\it one-vertex
triangulation}, a triangulation of the surface having just one
vertex; and while not so obvious, a minimum triangulation of a
closed, orientable 3--manifold (triangulation with the minimal
number of tetrahedra) requires a one-vertex triangulation, except
for $S^{3}$, and the lens spaces $\mathbb{R}P^{3}$ and $L(3,1)$.
It turns out, however, that by using one-vertex triangulations we
not only have the computational benefits but also can draw many
topological conclusions from their nice combinatorial properties.

\begin{thm} Every closed surface with $\chi \leq 0$ admits a
one-vertex triangulation.
\end{thm}

For example, any closed, orientable surface with genus $g \geq 1$
is the quotient of a 4g-gon in the plane, formed by identifying
edges in a way to give only one vertex. We can triangulate the
4g-gon by adding no additional vertices and $4g-3$ edges. This
induces a triangulation  of the genus g surface with one vertex,
$6g-3$ edges and $4g-2$ faces. The same construction also works
for closed, non-orientable surfaces with $\chi \leq 0$.

For 3-manifolds, it is not as easy to show that they admit
one-vertex triangulations and not as obvious (Euler characteristic
arguments do not work) to show that with the exceptions noted
above, a minimum triangulation must be a one-vertex triangulation.
We have the following result from \cite{j-r4}.

\begin{thm}
Every closed, orientable 3-manifold admits a one-vertex
triangulation. Furthermore, a compact, orientable 3-manifold with
non-empty boundary, no component of which is a 2-sphere, admits a
triangulation having all its vertices in the boundary and
precisely one vertex in each boundary component.
\end{thm}

There is a simpler version of this result that is satisfactory for
most of our work.

\begin{thm}
Given a triangulation $\T$ of a compact, orientable 3--manifold
with non-empty boundary, no component of which is a 2-sphere, then
$\T$ can be modified to a triangulation $\T'$ where $\T'$ has
precisely one vertex in each boundary component.
\end{thm}

From the previous theorem, we see that each compact, orientable
3-manifold (no boundary component a 2-sphere) admits a
triangulation that restricts to a one-vertex triangulation on each
boundary component. We shall exploit this and especially use such
triangulations for our study of knot-manifolds and Dehn fillings.

\subsection{Normal curves in a one-vertex triangulation of a torus}

In particular, we rely on some particularly nice properties of the
space of normal curves in a one-vertex triangulation of a torus,
$S^1 \times S^1$.  We will assume that our knot-manifold has such
a triangulation of its boundary and this will simplify
computations involving properly embedded surfaces.

\begin{figure}
{\epsfxsize = 5 in \centerline{\epsfbox{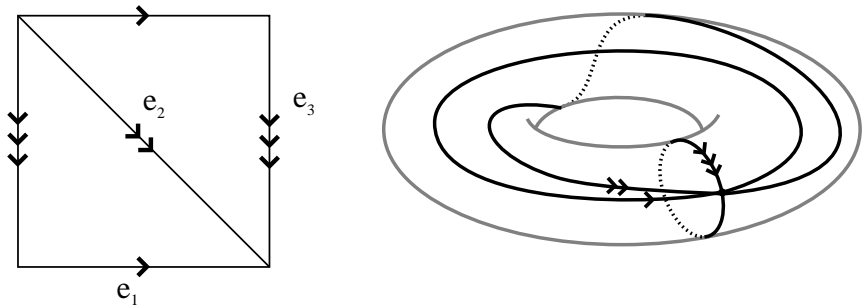}} }
\caption{The one-vertex triangulation of the torus $S^1 \times
S^1.$} \label{f-torus}
\end{figure}

Pictured in Figure \ref{f-torus} is the one-vertex triangulation
of a torus; the Euler characteristic of the torus determines that
it has 2 triangles and three edges. Note that we refer to `the'
one-vertex triangulation; any other one-vertex triangulation of
the torus is combinatorially equivalent to this one. The three
edges are essential curves which meet in a single point.   Any
other triangulation also has three edges which meet in a single
point and we can choose a homeomorphism of the torus mapping the
edges of the new triangulation, hence the triangulation itself, to
the triangulation of Figure \ref{f-torus}.

Among the nice properties is that there is a 1-1 identification
between normal curves and isotopy classes of curves on the torus.

\begin{lem}
\label{l:triv}
 In the one-vertex triangulation of a torus every
trivial normal curve is vertex-linking.
\end{lem}

\begin{proof}
A trivial normal curve $C$ bounds a disk $D$ on the torus.
Consider the intersection of this disk with the one-skeleton of
the triangulation. If there is an arc of intersection which is not
incident in $D$ to the vertex then it splits $D$ into two pieces,
at least one of which does not contain the vertex.  An outermost
arc of intersection with this subdisk demonstrates that the
trivial curve $C$ is not normal.  Therefore the intersection is a
collection of arcs each of which is incident in $D$ to the vertex.
This describes a vertex linking trivial curve.
\end{proof}

\begin{lem}
\label{l:isotopy} In a one-vertex triangulation of the torus two
normal curves are normally isotopic if and only if they are
isotopic.
\end{lem}

\begin{proof}
It suffices to consider the case where the two normal curves $C_1$
and $C_2$ are connected, essential, in general position with
respect to each other, and have been normally isotoped to
intersect minimally.  If $C_1$ and $C_2$ are disjoint then they
cobound 2 annuli on the torus.  One annulus , call it $A$, does
not contain the vertex.  As the boundary of $A$ consists of the
normal curves $C_1$ and $C_2$, each edge of the triangulation must
intersect $A$ in arcs running from $C_1$ to $C_2$. Thus, the edges
meet $A$ in a parallel collection of such arcs and we may use $A$
to perform a normal isotopy of $C_1$ to $C_2$.

When the curves do intersect, there must be at least two bigons on
the torus which are bounded by subarcs of $C_1$ and $C_2$. One of
these bigons does not contain a vertex and an innermost such
bounds a disk in which all edges of the triangulation intersect in
arcs joining $C_1$ to $C_2$.  We can use the bigon to construct a
normal isotopy reducing the number of intersections between $C_1$
and $C_2$.
\end{proof}

\begin{rem}
While Lemma \ref{l:triv} remains true for one-vertex
triangulations of any surface, Lemma \ref{l:isotopy} is never true
in a one-vertex triangulation of a surface of genus $\geq 2$. For
example, each separating curve on such a surface possess at least
two distinct normal representatives, determined by the side of the
curve to which the vertex lies.
\end{rem}

Recall that the isotopy class of an  essential simple closed curve
on the torus is called a {\it slope} on the torus.  If $C \subset
S^1 \times S^1 $ is a collection of pairwise disjoint curves with
at least one non-trivial component, then the {\it slope} of $C$,
denoted $slope(C)$, is the slope of one of the non-trivial
components. By the preceding lemma, when using a one-vertex
triangulation, it is equivalent to define {\it slope} as the
normal isotopy class of an essential simple closed curve.

In the one-vertex triangulation $\T$ of a torus there are six
normal arc types yielding variables, so the solution space and
projective solution space are embedded in six-dimensional space,
$\mathbb R^6$.  However, in computing these spaces the system
reduces to one with only three degrees of freedom and it becomes
more natural to think of the solution space and the projective
solution space as being embedded in $\mathbb R^3$. We will denote
these {\it representations of the solution and projective solution
spaces} by $\S_\T \subset \mathbb R^3$, and $\P_\T \subset \mathbb
R^3$, respectively.

\begin{thm}
Normal curves in a one-vertex triangulation $\T$ of a torus are
projectively parameterized by $\P_\T$, the set of rational points
in the 2-simplex $$ \{(x_1,x_2,x_3) | x_1+x_2+x_3=1, x_i \geq 0 \}
\subset \mathbb R^3.$$ The vertices of this simplex represent the
projective classes of the 3 edges of $\T$.
\end{thm}

\begin{proof}
Any normal curve in the one-vertex triangulation of the torus will
meet the two simplices of $\T$ in a collection of normal arcs from
the six types labeled $a_i$ in  Figure \ref{f-bound}.

\begin{figure}
{\epsfxsize = 2 in \centerline{\epsfbox{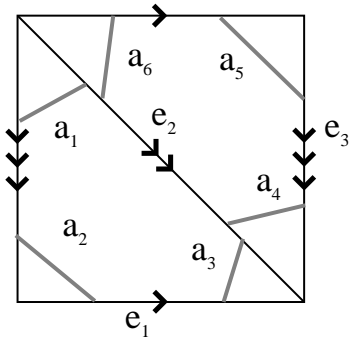}} }
\caption{Normal arcs in the one-vertex triangulation of the torus
$S^1 \times S^1$.} \label{f-bound}
\end{figure}

Therefore, a normal curve can be identified by a point
$(x_1,x_2,x_3,$ $x_4,$ $ x_5,x_6) \in \mathbb Z^6, x_i \geq 0$,
where $x_i$ denotes the number of arcs of the given type $a_i$.
Furthermore, for each of the three edges,
 each of the two triangles must have the same number arcs which
intersect that edge.  This yields three matching equations (see
Figure \ref{f-bound})
\begin{gather*}
x_1 + x_2 = x_4 + x_5 \text{, along edge } e_3, \\ x_1 + x_3 = x_4
+ x_6 \text{, along edge } e_2,
\\ x_2 + x_3 = x_5 + x_6 \text{, along edge } e_1,
\end{gather*}
which reduce to
\begin{gather*}
x_1 = x_4 \\ x_2 = x_5 \\ x_3 = x_6.
\end{gather*}

The solution space $$\S (S^1 \times S^1,\T)$$ is the set of points
with non-negative integer coordinates in the cone $$\{
(x_1,x_2,x_3,x_4,x_5,x_6) | x_i \geq 0, x_1 = x_4, x_2 = x_5, x_3
= x_6\} \subset \mathbb R^6.$$ However, it is more natural to
forget about the coordinates $x_4,x_5$ and $x_6$ and represent the
solution space by $\S_\T$ the set of points with non-negative
integer coordinates  $$ \{ (x_1,x_2,x_3) | x_i \geq 0 \} \subset
\mathbb R^3.$$ It can be seen from Figure \ref{f-bound} that the
solutions $(1,0,0), (0,1,0)$ and $(0,0,1)$ are the normal
coordinates of the isotopy classes of the three edges,
$e_1,e_2,e_3$, respectively, of the triangulation.  Moreover,
every solution to the normal equations can be written as a linear
combination of these three solutions using non-negative integer
coefficients, so they are the set of {\it fundamental solutions}
of $\S_\T$.

We projectivize the solution space $\S_\T$ by adding the
normalizing equation $$x_1 + x_2 + x_3 = 1.$$  Any solution to the
normal equations will have a unique projective representation as a
triple of non-negative rational numbers.   The resulting
projective space $\P_\T$  is the set of points in $$ \{
(x_1,x_2,x_3) | x_i \geq 0, x_1 + x_2 + x_3 = 1 \} \subset \mathbb
R^3.$$ See Figure \ref{f-psoln}.

\begin{figure}
{\epsfxsize = 2.5 in \centerline{\epsfbox{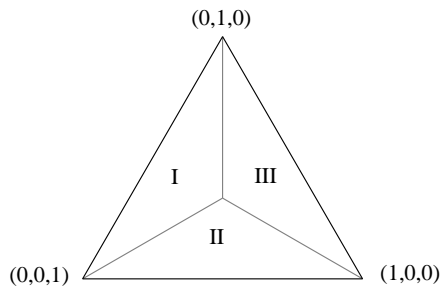}} }
\caption{The projective solution space $\P_\T$.} \label{f-psoln}
\end{figure}

Thus, the normal curves in $(S^1 \times S^1,\T)$ are projectively
parameterized by the set of rational points in the 2-simplex in
$\mathbb R^3$ spanned by $(1,0,0), (0,1,0)$ and $(0,0,1)$.  The
vertices are the projective classes of each of the fundamental
solutions, the three edges of the triangulation.

\end{proof}

In the remaining discussion we will refer to a normal curve by its
representation in $\S_\T$ or projective representation in $\P_\T.$
A normal curve will be called {\it Type I} if its $x_1$ coordinate
is less than or equal to each of its $x_2$ and $x_3$ coordinates.
{\it Type II}, and {\it Type III} are defined analogously.  See
Figure \ref{f-psoln}. Note that a curve may be of more than one
type. For example, a collection of trivial curves is
simultaneously all three types, and normal representatives of the
edges $e_1,e_2,e_3$ are two types. If $C$ is a family of normal
curves then we will let $\tau(C)$ denote the number of trivial
curves in $C$. Two slopes, $\alpha$ and $\beta$, will be said to
be {\it complementary} if $\alpha + \beta$ is a collection of
trivial curves.

We now state without proof some useful, elementary facts about
normal curves in a one-vertex triangulation of a torus.
\begin{enumerate}

\item The set of slopes on the torus is projectively
represented by the points in the boundary  of the projective space
$\P_\T \subset \mathbb R^3$.

\item If the normal curve $C$ has representation the triple
$(x_1,x_2,x_3) \in \S_\T$ then $\tau(C) = \min \{ x_1,x_2,x_3\}$.

\item The projective class of a collection of trivial curves is the barycenter
$(1/3,1/3,1/3)$ $\in \P_\T$.

\item If $C_1$ and $C_2$ are normal curves which are not the same type
then $\tau(C_1 + C_2) > \tau (C_1) + \tau(C_2) $.  If $C_1$ and
$C_2$ are normal curves which are the same type then $\tau(C_1 +
C_2) = \tau (C_1) + \tau(C_2).$

\item If $C$ is a normal curve with projective class $\bar{C}$ then
the slope of $C$ is determined by projecting the point $\bar{C}$
from the barycenter to the boundary of $\P_\T$. (Figure
\ref{f-slope}.)

\begin{figure}
{\epsfxsize = 2.5 in \centerline{\epsfbox{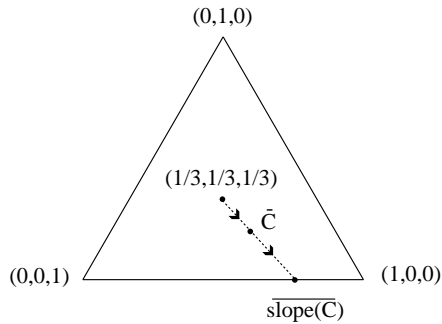}} }
\caption{Computing the slope of a normal curve.} \label{f-slope}
\end{figure}

\item The slopes $\alpha$ and $\beta$ are complementary if and only if
for any curve $C_\alpha$ with slope $\alpha$ and any curve
$C_\beta$ with slope $\beta$, the line segment in $\P_\T$
connecting $\bar{C_\alpha}$ and $\bar{C_\beta}$ passes through the
barycenter $(1/3,1/3,1/3)$.  Thus each slope has a unique
complement.

\item Suppose $C$ is a normal curve with parameterization
$(x_1,x_2,x_3) \in S_\T$.  If $\mu(C) = \max \{x_1,x_2,x_3\}$ and
$\tau(C) = \min \{x_1,x_2,x_3\},$ then $slope(C)$ has projective
class that of $(x_1 - \tau(C),x_2 - \tau(C),x_3 -\tau(C))$ and the
slope complementary to $slope(C)$ has projective class that of
$(\mu(C) - x_1, \mu(C) - x_2, \mu(C) - x_3).$
\end{enumerate}

\subsection{Boundary slopes}

In \cite{hatcher} Hatcher used the theory of incompressible
branched surfaces developed by Floyd and Oertel
\cite{floyd-oertel} to show that the slopes bounding
incompressible and $\bdy$-incompressible surfaces in a
knot-manifold are finite in number. Here we adapt Hatcher's
argument to normal surfaces in a one-vertex triangulation and show
that the result holds more generally for the slopes bounding
normal and almost normal surfaces; hence, our results imply
Hatcher's result for incompressible and $\bdy$-incompressible
surfaces as well.

\begin{prop}
\label{same-or-comp} Let $M$ be an orientable 3-manifold having a
boundary component a torus, $T$, and let $\T$ be a triangulation
of $M$ that restricts to a one-vertex triangulation of $T$.
Suppose that $S_1$ and $S_2$ are embedded normal or almost normal
surfaces and $\bdy S_1 \subset T$.  If $S_1$ and $S_2$ are
compatible and both meet $T$ in non-trivial slopes, then these
slopes are either equal or complementary.
\end{prop}

\begin{proof}
Let $\T'$ denote the induced one-vertex triangulation of the
boundary torus $T$.  We proceed in two steps. First, we show that
if the slopes of the surfaces $S_1$ and $S_2$ are the same type in
the one-vertex triangulation $\T'$ of the boundary torus $T$ then
they have the same slope in $T$.  Next, we show that if they have
different types in the triangulated torus $\T'$ then their slopes
in $T$ are complementary in $\T'$.

So we first assume that $S_1$ and $S_2$ are compatible surfaces
which intersect non-trivially and $\partial S_1 \cap T$ and
$\partial S_2 \cap T$ are of the same type in $\T'$, say type II.
First perform a normal isotopy of $S_1$ and $S_2$ so that all of
the trivial curves of $\partial S_1 \cap T$ and $\partial S_2 \cap
T$ are disjoint from all other curves.   All remaining
intersections on $T$ between $S_1$ and $S_2$ lie on the
non-trivial components of $\bdy S_1 \cap T$ and $\bdy S_2 \cap T$;
denote these by $C_1$ and $C_2$, respectively, and let
$(x_1,0,x_3)$ and $(y_1,0,y_3)$ denote their respective normal
coordinates in $\T'$.

The normal sum $S_1 + S_2$ restricts to a normal sum of the
boundary curves $\bdy S_1 + \bdy S_2$, hence to a sum of the
essential boundary curves in $T$, $C_1 + C_2$.  The normal curves
$C_1$,$C_2$ and $C_1 + C_2$ can be given an orientation by
orienting the normal arcs $a_1,a_3,a_4$ and $a_6$ as indicated in
Figure \ref{f-orient}. (Here we are using that $x_2 = y_2 = 0$.)

\begin{figure}
 {\epsfxsize = 2 in
\centerline{\epsfbox{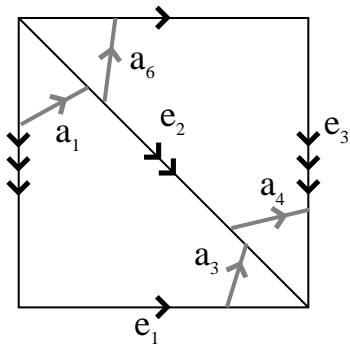}} } \caption{Orienting curves of
type II.} \label{f-orient}
\end{figure}

Consider an intersection between the normal curves $C_1$ and $C_2$
that lies along an intersection of normal arcs of type $a_1 \cap
a_1$, $a_1 \cap a_3$, or $a_3 \cap a_3$. The regular switch
performed at such an intersection is the switch that follows the
given orientation of the normal arcs, see Figure \ref{f-bsum}.  It
is easily verified that this is also true for intersections of
type, $a_4 \cap a_4$, $a_4 \cap a_6$ and $a_6 \cap a_6$, i.e. for
all intersections between $S_1$ and $S_2$ on the boundary
component $T$.

\begin{figure}
{\epsfxsize = 2 in \centerline{\epsfbox{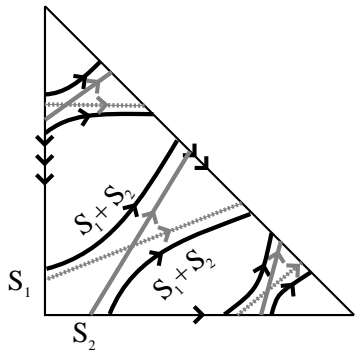}} } \caption{The
normal sum of curves of the same type.} \label{f-bsum}
\end{figure}

\begin{figure}
{\epsfxsize = 2.5 in \centerline{\epsfbox{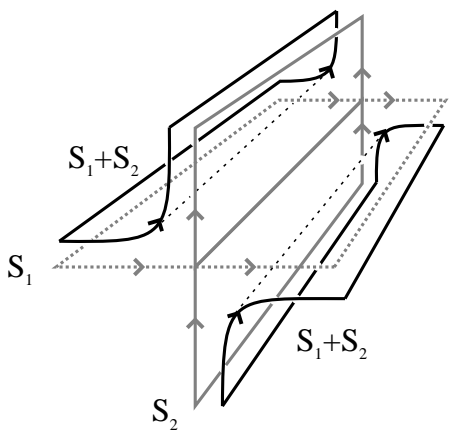}} }
\caption{An arc of  intersection between  $S_1$ and $S_2$.}
\label{f-intsum}
\end{figure}

If $a$ is an arc of intersection between the surfaces $S_1$ and
$S_2$, both of its endpoints are in $C_1 \cap C_2$ ($\bdy S_1
\subset T)$. At each endpoint we know that the regular switch
along $a$ follows the given orientation on the boundary curves. In
Figure \ref{f-intsum}, we follow the regular switch along $a$
through the interior of the orientable manifold $M$ and see that
the two endpoints of $a$ are intersections between $C_1$ and $C_2$
with opposite algebraic sign. If we consider all of the arcs of
intersection, hence all points in $C_1 \cap C_2$, we see that the
algebraic intersection between $C_1$ and $C_2$ sums to 0.  Since
$T$ is a torus the surfaces $S_1$ and $S_2$ have the same slope on
$T$.

We now assume that $S_1$ and $S_2$ are compatible, intersect
non-trivially on $T$ and do not have the same type in the
triangulated torus $\T'$.  As $S_1$ and $S_2$ are compatible each
member of the collection $ \{ n_1 S_1 + n_2 S_2 : n_1,n_2 \geq 0
\}$ is an embedded normal surface contained in the same
compatibility class and with non-empty boundary in $T$.
Representing this collection of surfaces by their normal
boundaries in $\T'$ and projectivizing, we obtain the set of
rational points on the segment joining $\partial S_1$ and
$\partial S_2$ in $\P_{\T'}$, where the endpoints have different
types. We can choose some surface $S_1'=n_1' S_1 + n_2' S_2$ so
that $\bdy S_1' \cap T$ has the same type as $\partial S_1 \cap
T$. See Figure \ref{f-incomp}.

\begin{figure}
{\epsfxsize = 2.25 in \centerline{\epsfbox{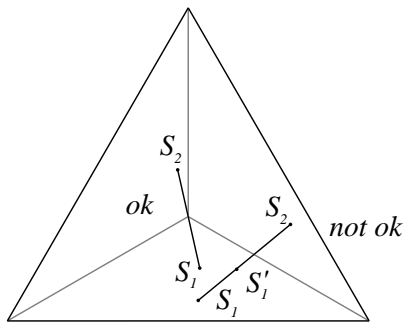}} }
\caption{Compatible surfaces of different types are
complementary.} \label{f-incomp}
\end{figure}

By the first step of the proof, the slopes of $S_1$ and $S_1'$ on
$T$ are identical.  This implies that the segment joining $S_1$
and $S_2$ passes through the point $(1/3,1/3,1/3) \subset
\P_{\T'}$ and in particular means that the slopes of $S_1$ and
$S_2$ on $T$ are complementary.  (Recall the elementary facts 5
and 6).
\end{proof}

\begin{rem}
In Section \ref{s-filling} we will give an example of compatible
normal surfaces with complementary slopes; hence it is necessary
to include the two possibilities unlike the situation of
\cite{hatcher}.  Our methods include more general surfaces, normal
and almost normal as opposed to incompressible and
$\bdy$-incompressible surfaces, and hence the corresponding
branched surfaces may have monogons and bigons in their
boundaries.
\end{rem}

Since a normal or almost normal surface $S$ is compatible with all
surfaces in its carrier we obtain the following.

\begin{cor}
\label{carrier-slopes} Let $X$  be a knot-manifold with a
triangulation $\T$ that restricts to a one-vertex triangulation on
$\bdy X$.   Suppose $S$ is an embedded normal or almost normal
surface and $\bdy S \neq \emptyset$. There are at most two slopes
(complementary ones) for all surfaces in the carrier of $S$,
$\mathcal C(S) \subset \P(X,\T)$.
\end{cor}

We note that if $S$ has no non-trivial boundary components, then
every surface in the carrier of $S$ has the same slope as $S$,
this is the case for $S$ incompressible. Now, if $S$ is a normal
or almost normal surface then some multiple of $S$ can be written
as a sum of embedded surfaces represented at the vertices in the
carrier of $S$, $$k S = \sum k_i V_i.$$ Hence, in a knot-manifold,
there can be at most two distinct boundary slopes for these
summands, from which the slope of $S$ is inherited.

\begin{cor}
\label{slopes-at-vertices} Let $X$ be a knot-manifold with a
triangulation $\T$ that restricts to a one-vertex triangulation on
$\bdy X$.  All possible slopes for the boundaries of embedded
normal or almost normal surfaces in $X$ are realized by the slopes
of embedded surfaces represented at the vertices of $\P(X,\T)$.
\end{cor}

\begin{cor} \label{slopes-are-finite}
Let $X$ be a knot-manifold with a triangulation $\T$ that
restricts to a one-vertex triangulation on $\bdy X$. Then there
are only a finite number of slopes realized as the slopes of
normal and almost normal surfaces in $(X,\T)$.
\end{cor}

The number of tetrahedra, $t$, in the triangulation $\T$ yields a
very  rough upper bound on the number of slopes bounding normal
and almost normal surfaces.  Let $\{ S_1,$ $\dots,S_n\}$ be a
maximal collection of normal surfaces with distinct slopes.  There
is a sub-collection with at least $n/2$ surfaces, no two of which
are compatible. For each pair of these $n/2$ surfaces there is a
tetrahedron in which they possess distinct quadrilateral types. In
the worst case possible, each surface in the sub-collection
possess a quadrilateral in each tetrahedron, implying that $n/2
\leq 3^t$. Thus, $2 (3^t)$ is an upper bound on the number of
slopes bounding normal surfaces.

A similar computation works for almost normal surfaces. Let
$\{S_1, \dots,$ $S_n\}$ be a maximal collection of almost normal
surfaces with distinct slopes.   If $S_i$  possesses a tube in
some tetrahedron, then we may compress the tube to obtain a normal
surface with the same slope.  There is a sub-collection of at
least $n/2$ surfaces, no pair of which are compatible.  In the
worst case  each surface possesses an octagon in one tetrahedron
and quadrilaterals in all others.  There are $3t$ choices for the
octagon and in each of the remaining $t-1$ tetrahedra we choose a
quadrilateral. Therefore, $n/2 \leq 3 t 3^{t-1}$ and there are at
most  $2 t 3^t$ slopes of almost normal surfaces.

We have already noted that every knot-manifold possesses a
one-vertex triangulation.  The slope of every incompressible and
$\bdy$-incompres\-sible surface is realized by the slope of a
normal surface, surface (see \cite{j-t} and Section
\ref{s-essential}). This implies Hatcher's theorem \cite{hatcher}.

\begin{cor}
Let $X$ be a knot-manifold.  Then there are a finite number of
slopes bounding incompressible and $\bdy$-incompressible surfaces
in $X$.
\end{cor}


\section{Layered triangulations of the  solid torus}
\label{s-filling} \label{s-layered}

In this section we give a method, also used in \cite{j-r3} and
\cite{j-r4}, for extending a one-vertex triangulation of a
knot-manifold $X$ to that of a manifold $X(\alpha)$ obtained by
Dehn filling.  This is accomplished by showing that a one-vertex
triangulation on the boundary of a solid torus can be extended to
a special one-vertex triangulation of the solid torus (Theorem
\ref{exist-layered}).

The special one-vertex triangulations referred to above are {\it
layered triangulations} of solid tori. We are able to give a
classification of the embedded, planar, normal surfaces in a
layered triangulation of a solid torus (Proposition
\ref{layered-planars}), which will be of use in Section
\ref{s-essential}.  The classification is given in terms of the
three types of normal surfaces, defined as follows:

\begin{enumerate}
\item $D_\mu$ will designate any
normal disk with essential boundary, i.e., a meridional disk for
the solid torus.
\item $D_\tau$ will designate any normal disk
with trivial boundary, i.e., a disk which is parallel into the
boundary of the solid torus.
\item $A_\alpha$ will designate any annulus with essential boundary
which is parallel into an annulus in the boundary,  so that the
parallel annulus in the boundary contains the vertex of the
triangulation.
\end{enumerate}

\begin{figure}[h]
 {\epsfxsize = 2 in
\centerline{\epsfbox{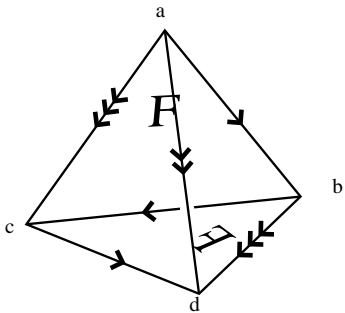}} } \caption{The one-tetrahedron
solid torus. } \label{f-storus}
\end{figure}

\begin{figure}[h]
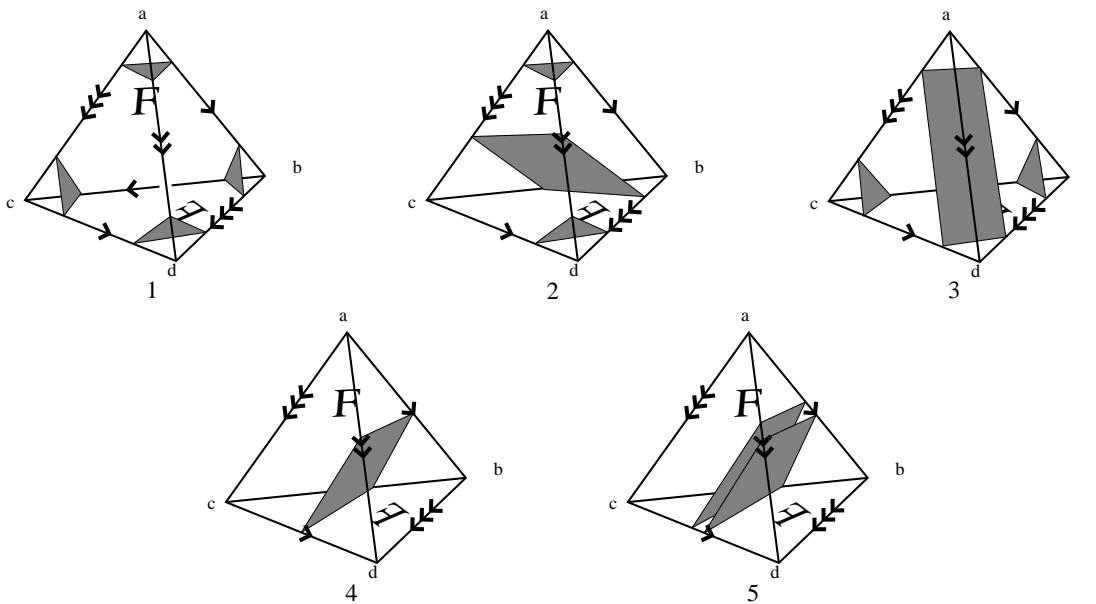

{\epsfxsize = 1.75 in \centerline{
  \epsfxsize = 1.75 in \epsfbox{F-STDT.AI} \qquad
  \epsfxsize = 1.75 in \epsfbox{F-STDM.AI} \qquad
  \epsfxsize = 1.75 in \epsfbox{F-STA1.AI} \qquad
}} {\epsfxsize = 1.75 in \centerline{
  \epsfxsize = 1.75 in \epsfbox{F-STMB.AI} \qquad
  \epsfxsize = 1.75 in \epsfbox{F-STA2.AI} \qquad
}}
\vspace{.2in}
\caption{The normal surfaces in a one-tetrahedron solid torus.
} \label{f-stsurf}
\end{figure}

Consider the tetrahedron $\Delta$ pictured in Figure
\ref{f-storus}. Glue the back two faces of the tetrahedron
together by making the ordered identification $\langle a,b,c
\rangle \sim \langle b,c,d \rangle$. All four vertices are
identified to a single vertex and the induced edge identifications
are indicated in the figure. It is easy to check that $M =
\Delta/\sim$ is a manifold with a single torus boundary component,
where the boundary torus has a one-vertex triangulation consisting
of the two front faces of the tetrahedron.  The normal surface
consisting of two triangles cutting of the vertices $\langle a
\rangle$ and $\langle d \rangle$ along with the quadrilateral
which separates the edges $\langle a,b \rangle$ and $\langle c,d
\rangle$ is a properly embedded disk, $D$. See Figure
\ref{f-stsurf}(2). Moreover, after cutting along the disk $D$ the
resulting manifold $M-N(D)$ is a ball.  We have therefore
constructed a triangulation of a solid torus with one tetrahedron,
three faces, three edges (each contained in the boundary), and a
single vertex. This triangulation of the solid torus will be
referred to as {\it the one-tetrahedron solid torus}.  Any other
triangulation of a solid torus with a one-tetrahedron is
combinatorially equivalent to this one: such a triangulation must
be obtained by gluing two (adjacent) faces of a tetrahedron
together with an orientation-reversing identification; the ordered
identification $\langle a,b,c \rangle \sim \langle c,d,b \rangle$
is equivalent by relabelling $\Delta$, $\langle b \rangle
\leftrightarrow \langle d \rangle$; and the ordered identification
$\langle a,b,c \rangle \sim \langle d,b,c \rangle$ forms a
triangulation of the ball.

The connected normal surfaces contained within the one-tetrahedron
solid torus are determined by their  quadrilateral type (or lack
thereof).   Pictured in Figure \ref{f-stsurf} are all of the
connected normal surfaces which can be properly embedded in the
one-tetrahedron solid torus.  They are:
\begin{enumerate}
\item A disk of type $D_\tau$ with boundary the trivial curve
$(1,1,1)$.
\item A disk of type $D_\mu$ with boundary
  $(2,0,1)$.
\item An annulus of type $A_\alpha$ with boundary $(0,2,0)$.
\item A M\"obius band with boundary $(0,0,1)$.
\item An annulus of type $A_\alpha$ which is the double of
  the M\"obius band. It has
  boundary $(0,0,2)$.
\end{enumerate}

\begin{figure}[h]
 {\epsfxsize = 2.15 in
\centerline{\epsfbox{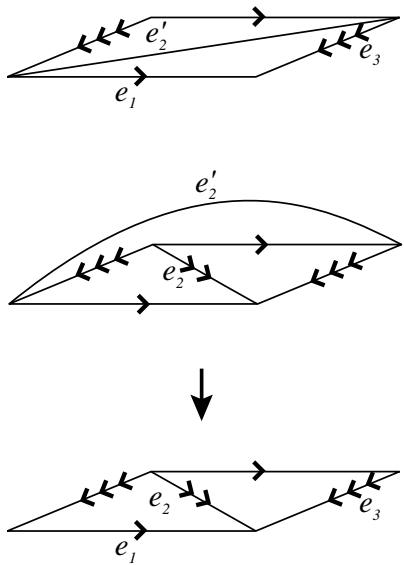}} } \caption{A Pachner move via
layering.} \label{f-layer}
\end{figure}

If $\T$ is a one-vertex triangulation of a solid torus then the
boundary torus has a one-vertex triangulation.  Any one of the
three edges $e$ in this triangulation may be thought of as the
diagonal of the rectangle bounded by the other two edges.  We can
change the boundary triangulation to a new one  by exchanging $e$
for $e'$, the other diagonal of the rectangle, this is known as a
{\it Type $I$ Pachner move}.  Fortunately, we can realize the Type
$I$ Pachner move by gluing an additional tetrahedron $\Delta$ to
the boundary of $\T$. See Figure \ref{f-layer} with $e=e_2$.

Glue the edge $e$ in the boundary torus to an edge $e$ of a
disjoint tetrahedron $\Delta$.  In addition glue the two faces on
the boundary torus that are adjacent to $e$ to the faces adjacent
to $e$ on $\Delta$. The result is a one-vertex triangulation of a
solid torus with the boundary changed by a Type $I$ Pachner move.

This move on a triangulation of a solid torus will be called {\it
layering} at the edge $e$ and we denote the new triangulation by
$\T' = \T \cup_e \Delta$.  Inductively, define a {\it layered
triangulation of a solid torus with $t$ layers}, $\T_t$, to be any
triangulation of a solid torus so that,
\begin{enumerate}
\item $\T_1 = \T$, a one-tetrahedron solid torus,
\item $\T_t = \T_{t-1} \cup_e \Delta_t, t \geq 2$, a layering at
$e$ of a layered triangulation with $t-1$ layers.
\end{enumerate}
Note that layering a solid torus has the effect of covering the
boundary edge $e$ and adding a new boundary  edge $e'$.  Thus
$\T_t$ will possess one vertex in the boundary torus and $t+2$
edges, 3 contained in the boundary torus.

We now give a theorem of \cite{j-r4}, that layered triangulations
are general enough to perform an arbitrary Dehn filling.

\begin{thm}
\label{exist-layered} Suppose $\T$ is a one-vertex triangulation
of the torus $S^1 \times S^1$.  For any slope $\mu$ on $\T$ there
is an algorithm to extend $\T$ to a layered triangulation of a
solid torus in which $\mu$ bounds a meridional disk. Furthermore,
for any positive integer $N$ there is such a layered triangulation
with greater than $N$ tetrahedra.
\end{thm}

\begin{proof}
We construct the layered triangulation in reverse order, layering
tetrahedra on the prescribed boundary, altering the normal
representative of $\mu$ until it is a $(2,0,1)$ curve.  This curve
(perhaps after relabelling) bounds a meridional disk in the one
tetrahedron solid torus, so we may glue our layers to the one
tetrahedron solid torus in reverse order and obtain the desired
layered triangulation.

We will keep track of $\mu$  by its {\it intersection numbers}, a
triple which indicates the number of intersections between $\mu$
and each of the three edges of the triangulation of the boundary
torus, $$[y_1,y_2,y_3] = [\#(\mu \cap e_1),\#(\mu \cap e_2),\#(\mu
\cap e_3)].$$  Note that these are {\it not} the normal
coordinates of $\mu$ in the solution space $\S_\T$ (see Section
\ref{s-triang}). We may convert from normal coordinates to
intersection numbers as follows $$[y_1,y_2,y_3] =
[x_2+x_3,x_1+x_3,x_1+x_2].$$ The length of a normal curve,
$L(\mu)$, is the sum of its intersection numbers $L(\mu) = y_1 +
y_2 + y_3 = 2 x_1 + 2 x_2 + 2x_3,$ an even number.

The slope $\mu$ is uniquely represented by a normal curve in the
triangulation of the torus.  Attaching a layer at the edge $e_2$,
is equivalent to a Type $I$ Pachner move, replacing the set of
edges $\{e_1,e_2,e_3\}$ by the edges $\{e_1,e_2',e_3\}$. Choose an
orientation on the torus and $\mu$ and orient $e_1$ and $e_2$ so
that the oriented intersection numbers $<\mu,e_1> =y_1$ and
$<\mu,e_3> = y_3$.  The edges $e_1$ and $e_3$ are a basis for the
homology of the boundary torus and the edge $e_2$ intersects each
once, so we may orient $e_2$ so that either $e_2 = e_1 + e_3$ or
$e_2 = e_1 - e_3$ with respect to homology. Thus, $y_2 = y_1 +
y_3$ or $y_2 =| y_1 - y_3|$. Then $e_2'$ can be oriented so that
$e_2'=e_1 - e_3$ or $e_2' = e_1 + e_3$, respectively, and then
$y_2' =| y_1 - y_3|$ or $y_2' = y_1 + y_3$, respectively.  Thus,
layering a tetrahedron on the boundary at $e_2$ changes the
intersection numbers of $\mu$, $$[y_1,y_1+y_3,y_3] \leftrightarrow
[y_1,|y_1-y_3|,y_3],$$ with the direction of the map determined by
whether $y_2 = y_1 + y_3$ or $y_2 = |y_1 - y_3|$.

Layering tetrahedra at the other edges alters the corresponding
intersection coordinate in precisely the same manner. By attaching
to the edge with highest intersection coordinate, $L(\mu)$ will
strictly decrease, unless with respect to some ordering of the
edges, $y_1 + y_3 =|y_1 - y_3|$, i.e. one intersection coordinate
is zero. This means that $\mu$ is in fact disjoint from some edge
and is therefore the normal representative of that edge. With
respect to some ordering of the edges, $\mu$ has intersection
numbers $[0,1,1]$ and $L(\mu) =2$.

If we ever have that $L(\mu) =6$ then the intersection triple of
$\mu$ is $[1,3,2]$, up to ordering. In normal coordinates this is
$(2,0,1)$, the curve that bounds a meridional disk in the one
tetrahedron solid torus. We may choose an ordering of the edges so
that we are able to glue the layered tetrahedra, in reverse order,
to the one-tetrahedron solid torus and obtain a layered
triangulation of the solid torus in which $\mu$ bounds a
meridional disk.

If the original length is $L(\mu)=4$ then the intersection triple
for $\mu$ is $[2,1,1]$ after a choice of edges.  By layering a
tetrahedron at either $e_2$ (or $e_3$, assuming this ordering) we
obtain the triple $[2,3,1]$.  As noted in the previous paragraph
we may then attach the one-tetrahedron solid torus. (Assuming this
ordering, layering at $e_1$ lowers $L(\mu)$.)  If the original
length is $L(\mu)=2$ then its intersection triple is $[0,1,1]$, up
to ordering.  In this ordering, by layering at the first edge we
change the intersection triple $[0,1,1] \mapsto [2,1,1]$.
(Layering at the edges with intersection values 1 only changes the
ordering). We then add one more layer as in the previous case.

If the original length $L(\mu) > 6$ then we may layer a sequence
of tetrahedra from the boundary, always attached to the edge with
highest intersection coordinate.  Continue this process, strictly
decreasing $L(\mu)$ until $L(\mu) \leq 6$.    In fact, this
process must terminate with $L(\mu) = 6$.  Because $L(\mu)$ is
strictly decreasing, and by the remarks of the previous paragraph,
if the process terminates at $L(\mu) =2$ or $L(\mu) =4$ then
$L(\mu) = 6$ was a previous step.

It is easy to obtain such a triangulation with an arbitrary number
of tetrahedra.  First layer $N$ tetrahedra on the boundary torus
in any fashion (keeping track of $\mu$).  Then, as specified
above, layer tetrahedra which reduce the length of $\mu$ until the
boundary can be capped off with the one-tetrahedron torus. A
one-vertex triangulation of a solid torus with at least $N+1$
tetrahedra is obtained.
\end{proof}

In Section \ref{s-essential} it will be necessary for us to
understand the normal  surfaces that can be embedded in a layered
triangulation of a solid torus $\T_t = \T_{t-1} \cup_e \Delta_t$.
Layering identifies two faces of $\Delta_t$ with the boundary of
$\T_{t-1}$ and leaves the other two faces of $\Delta_t$ as the new
boundary torus.  The one-tetrahedron solid torus has no closed
normal surfaces and each elementary disk type in $\Delta_t$ meets
both the old boundary and the new boundary.  It follows that in a
layered triangulation of a solid torus, there can be no closed
normal surfaces and each normal surface intersects each
tetrahedron.

So if $P_t \subset \T_t$ is a normal surface, then it was obtained
by attaching a non-empty collection of elementary disks in
$\Delta_t$ to a normal surface $P_{t-1} \subset \T_{t-1}$.  Call
the elementary quadrilateral type in $\Delta_t$ which separates
the attaching edge $e$ and the new edge $e'$ the {\it banding
quad} in $\Delta_t$.  See Figure \ref{f-aa1}.  The boundary of
$P_{t-1}$ and the number of banding quads attached in $\Delta_t$
completely determine the number of each of the other disk types
that are used in $\Delta_t$.  In particular, if no banding quads
are attached then the surface $P_t$ is homeomorphic to the surface
$P_{t-1}$, and we say that $P_t$ was obtained by {\it pushing
$P_{t-1}$ through $\Delta_t$}.   See Figure \ref{f-push}.

We are particularly interested in the planar normal surfaces
embedded in $\T_t$.  The one-tetrahedron solid torus $\T_1$
contains a unique surface of type $D_\mu$, a unique surface of
type $D_\tau$ and two distinct surfaces of type $A_\alpha$.  How
can new normal planar surfaces be obtained by the process of
layering ?  We list some (all, by Proposition
\ref{layered-planars}) ways to construct new planar surfaces $P_t
\subset \T_t$ from planar surfaces $P_{t-1} \subset \T_{t-1}$:
\begin{enumerate}
\item In any layer, each surface of type $D_\mu,D_\tau$ or $A_\alpha$
in $\T_{t-1}$ may be pushed through $\Delta_t$ to obtain a surface
of the same type in $\T_t$. See Figure \ref{f-push}.

\begin{figure}[h]
 {\epsfxsize = 2 in \centerline{\epsfbox{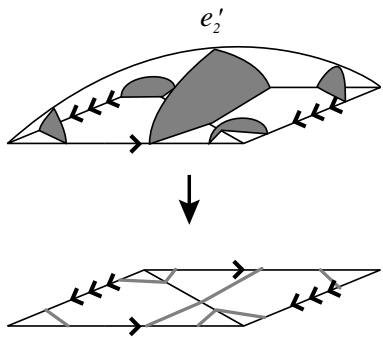}}
} \caption{Pushing through.} \label{f-push}
\end{figure}

\item In any layer $\Delta_t$ we may attach a banding quad and two
elementary triangles to a surface of type $D_\tau \subset
\T_{t-1}$ to produce a surface of type $A_\alpha \subset \T_t$,
see Figure \ref{f-aa1}. The band was attached along the edge $e'$
so this annulus is parallel to a neighborhood of the edge $e'$, an
annulus in the boundary containing the vertex.

\begin{figure}[h]
 {\epsfxsize = 2 in \centerline{\epsfbox{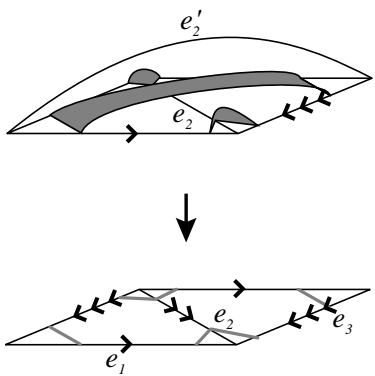}} }
\caption{Banding a trivial disk $D_\tau$ to create an annulus
$A_\alpha$. } \label{f-aa1}
\end{figure}

\item If the attaching edge $e$ happens to be the slope of the
meridional disk $D_\mu \subset \T_{t-1}$, then attach a single
banding quad and two triangles in $\Delta_t$ to 2 copies of
$D_\mu$ to obtain a surface of type $D_\tau \subset \T_t$.  See
Figure \ref{f-dt}. (Banding two meridional disks in a solid torus
produces a trivial disk.)

\begin{figure}[h]
 {\epsfxsize = 2 in \centerline{\epsfbox{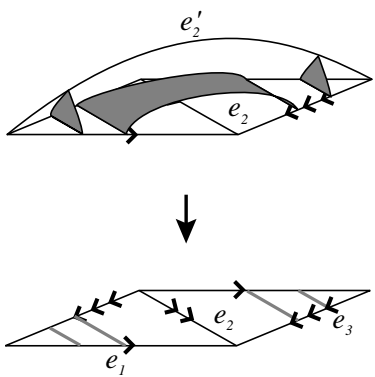}} }
\caption{Banding 2 meridional disks $D_\mu$ to create a trivial
disk $D_\tau$.} \label{f-dt}
\end{figure}

\item If the attaching edge $e$ happens to be the slope of the
meridional disk $D_\mu \subset \T_{t-1}$, then attach two banding
quads to two copies of $D_\mu$ to obtain a surface of type
$A_\alpha \subset \T_t$.  See Figure \ref{f-aa2}.  This annulus is
also parallel to a neighborhood of the edge $e'$, an annulus in
the boundary containing the vertex.

\begin{figure}[h]
{\epsfxsize = 2 in \centerline{\epsfbox{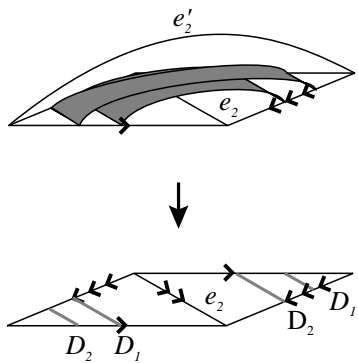}} }
\caption{Attach two bands to two meridional disks $D_\mu$ to
create an annulus $A_\alpha$.} \label{f-aa2}
\end{figure}
\end{enumerate}

In the proof of the following proposition we show that the above
moves are sufficient to generate all normal planar surfaces in a
layered solid torus.

\begin{prop}
\label{layered-planars} Let $\T_t$ be a layered triangulation of a
solid torus $T$ with $t$ layers.  Then the only  connected planar
normal surfaces which can be properly embedded in $(T,\T_t)$ are
of type $D_\mu, D_\tau$, and $A_\alpha$.  Furthermore, there is
unique (up to normal isotopy) surface of type $D_\mu$. The weight
of any surface of these types is bounded below by:
\begin{enumerate}
\item $wt(D_\mu) \geq t+4$,
\item $wt(D_\tau) \geq 2(t+2)$, and
\item $wt(A_\alpha) \geq 2(t+1)$.
\end{enumerate}
\end{prop}

\begin{proof}
We first give an inductive proof that all normal planar surfaces
in $\T_t$ are of type $D_\mu,D_\tau$ or $A_\alpha$.  When $t=1$
the triangulation $\T_1$ is the one-tetrahedron triangulation of
the solid torus. The normal planar surfaces were shown in Figure
\ref{f-stsurf}, each is of type $D_\mu, D_\tau$ or $A_\alpha$.
Assume that the result holds for any layered triangulation
$\T_{t-1}$ with $t-1$ tetrahedra, $t \geq 2$. We now show that the
result also holds for any layered triangulation with $t$
tetrahedra, $\T_t=\T_{t-1} \cup_e \Delta_t$.  With no loss of
generality we may assume that $e = e_2$.

Let $P_t$ be a connected normal planar surface in $\T_t$. Then
$P_{t-1} = P_t \cap \T_{t-1}$ is a (possibly disconnected) planar
normal surface in $\T_{t-1}$.

\begin{claim} We may assume that every component of  $P_{t-1}$ has
a banding quad in $\Delta_t$ attached to it.
\end{claim}

If any component of $P_{t-1}$ does not have a banding quad
attached to it, then it is merely pushed through the tetrahedron
$\Delta_t$ to a surface which is of precisely the same type in
$\T_t$. This component satisfies the conclusion of the theorem and
in particular is homeomorphic to $P_t$ ($P_t$ is connected). We
therefore assume that every component of $P_{t-1}$ has a banding
quad in $T$ attached to it.

\vspace{.2in}

A banding quad is a band along the edge $e_2'$ and joins a normal
arc of $\bdy P_{t-1}$ of type $a_2$ to a normal arc of type $a_5$.
The number of arcs of type $a_2$ is equal to the number of type
$a_5$ for any normal curve (recall that $x_2 = x_5$).  Number the
arcs of type $a_2$ with the numbers $1,\dots,x_2$ counting from
the vertex to the edge $e_2$, and number the arcs labeled $a_5$
from $1,\dots,x_2$ also counting from the vertex to the edge
$e_2$, see Figure \ref{f-arcs}.   If a banding quad is attached to
an arc of type $a_2$ labeled $i$ then all arcs of type $a_2$ with
greater labels must also have banding quads attached, it is
impossible to attach a normal triangle, see Figure \ref{f-dt}. The
same holds true for arcs of type $a_5$, and it follows that each
banding quad joins an arc of type $a_2$ to an arc of type $a_5$
with the same label.

\begin{figure}[h]
{\epsfxsize = 1.5 in \centerline{\epsfbox{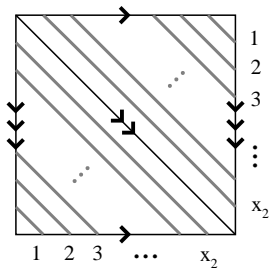}} }
\caption{Numbering the arcs of type $a_2$ and $a_5$.}
\label{f-arcs}
\end{figure}

\begin{claim}
If some component of $P_{t-1}$ is of type $D_\tau$  then $P_t$ is
an annulus $A_\alpha$.
\end{claim}

We are assuming that there is a banding quad attached to the
component of $P_{t-1}$ of type $D_\tau$.   The boundary of this
component consists of a single trivial curve, so it possesses one
$a_2$ arc and one $a_5$ arc. Moreover they must have the same
label as there can only be trivial curves between $\bdy D_\tau$
and the vertex, and each trivial curve possesses an equal number
of arcs of type $a_2$ and $a_5$. So the banding quad attached to
$D_\tau$ joins $D_\tau$ to itself.    This is move (2) described
before the proposition, which creates an annulus $A_\alpha \subset
\T_t$, see Figure \ref{f-aa1}.  Since $P_t$ is connected, $D_\tau$
was the only component of $P_{t-1}$.

\begin{claim}
If no component of $P_{t-1}$ is of type $D_\tau$ then $\bdy
P_{t-1}$ consists of essential curves parallel to the edge $e_2$.
\end{claim}

We are assuming that $\bdy P_{t-1}$ has a banding quad attached so
it must possess a normal arc of type $a_2$ hence the coordinate
$x_2 > 0$. Then $\bdy P_{t-1}$ has no trivial curves: by the
inductive hypothesis a trivial curve implies that $P_{t-1}$
contains a trivial disk $D_\tau$ and by the conclusions of the
previous claim $P_{t-1}$ must itself be a trivial disk $D_\tau$.
Therefore one of the two coordinates $x_1,x_3$ must be 0. No
normal arc of $\bdy P_{t-1}$ of type $a_1$ can be connected across
$e_2$ to an arc of type $a_4$ (nor $a_3$ to $a_6$), for this would
force a quad which is not the banding quad to be attached, which
prohibits any banding quads from being attached and means that the
surface was pushed through $\Delta_t$, see Figure \ref{f-push}.
Therefore all arcs of type $a_1$ are connected across $e_2$ to
those of type $a_6$ and $x_1 = x_6 (=x_3)$. Then both $x_1=x_3=0$
and only $x_2
>0$. The normal curve $\bdy P_{t-1}$ consists entirely of curves
parallel to the edge $e_2$.

\begin{claim}
No component of $P_{t-1}$ is an annulus $A_\alpha$ (or $P_t$ has
positive genus).
\end{claim}

By our previous claim,  any annulus $A_\alpha \subset P_{t-1}$ has
boundary disjoint from $e_2$ and is parallel into the boundary
annulus containing the vertex, i.e. parallel to a neighborhood of
the edge $e_2$. Thus any collection of such annuli is nested, and
we may choose an innermost annulus $A_\alpha$ with respect to the
edge $e_2$.  No component of $P_{t-1}$ is of type $D_\tau$ and any
component of type $D_\mu$ cannot have boundary contained in the
annulus in the boundary to which the $A_\alpha$'s are parallel. So
the boundary of the innermost $A_\alpha$ is adjacent to the edge
$e_2$, i.e. it has arcs of type $a_2$ and $a_5$ with label $x_2$.
A banding quad is attached from $\bdy A_\alpha$ to itself along
these arcs yielding a once punctured torus. We may either attach
two elementary triangles or a banding quad to the remaining
boundary arcs of type $a_2$ and $a_5$ (Figure \ref{f-dt} or Figure
\ref{f-aa2}). However, in either case the surface $P_t$ has
positive genus, a contradiction.

We are left with the case that $P_{t-1}$ is a collection of
meridional disks $D_\mu$.

\begin{claim}
$P_{t-1}$ is not a single copy of $D_\mu$ (for then $P_t$ would be
a M\"obius band).
\end{claim}

If $P_{t-1}$ is a single copy of $D_\mu$ then the banding quad is
glued from the single normal arc of $\bdy P_{t-1}$ of type $a_2$
to that of type $a_5$. The surface produced has a single boundary
component and $\chi = 0$, hence it is a M\"obius band (with
boundary $e_2'$).

\begin{claim}
If $P_{t-1}$ is 2 copies of $D_\mu$ and a single banding quad is
attached then $P_t$ has type $D_\tau$. If $P_{t-1}$ is 2 copies of
$D_\mu$ and  2 banding quads are attached then $P_t$ has type
$A_\alpha$.
\end{claim}

These cases are the moves (3) and (4) listed preceding the
theorem.

\begin{claim}
$P_{t-1}$ does not contain more than 2 meridional disks $D_\mu$
(for then $P_t$ would be disconnected).
\end{claim}

Let $D_1,\dots,D_{x_2}, x_2 \geq 2$ be a collection of meridional
disks $D_\mu$ numbered to induce our previous labeling of the arcs
of type $a_2$.  See Figure \ref{f-arcs} and Figure \ref{f-disks}.
The boundary of $D_i$ is parallel to $e_2$ and consists of an arc
of type $a_2$ labeled $i$ along with an arc of type $a_5$ labeled
$x_2 - i +1$. Since at least one banding quad is attached, there
is necessarily a banding quad attached to the two arcs labeled
$x_2$, this quad bands $D_{x_2}$ to $D_1$. There is either a
banding quad attached to the two arcs labeled 1 or elementary
triangles are added to each. In either event, $D_1$ and $D_{x_2}$
are attached to each other and to no other disk. Then  there must
be no other disks, for then $P_t$ would be disconnected.

\begin{figure}[h]
{\epsfxsize = 1.5 in \centerline{\epsfbox{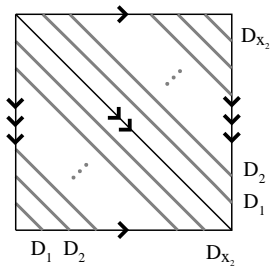}} }
\caption{A collection of meridional disks.} \label{f-disks}
\end{figure}

Note that a surface of type $D_\mu$ was created only by pushing
through each layer $\Delta_i$, a unique process.  For a given
layered solid torus $\T_t$, there is a unique surface of type
$D_\mu$.

So we have that a normal planar surface $P_t \subset \T_t$ is one
of the three types, $D_\mu, D_\tau$ and $A_\alpha$.  We now obtain
lower bounds on their weights.  Let $P_i = P_t \cap \T_i$.
Typically $P_t$ meets each of the $t+2$ edges of the
triangulation. However, there are three ways that a planar surface
$P_t$ can miss an edge of the triangulation $\T_t$:
\begin{enumerate}
\item $P_t$ does not intersect some edge $e$ in the core
triangulation $\T_1$. This happens only when the surface $P_1$ was
one of the two annuli of type $A_\alpha \subset \T_1$.  It follows
that $P_t$ is of type $A_\alpha$ and was obtained by pushing
through every subsequent layer.

\item In some layer the surface $P_i$ is obtained by attaching a
banding quad to the surface $P_{i-1}$ and $P_i$ misses the new
edge. Then $P_i$ has the same slope as the new edge $e'$, and is
therefore an annulus $A_\alpha$. (A trivial disk $D_\tau$ has
trivial boundary which intersects the new edge, and the type
$D_\mu$ cannot be created through banding.)

\item The surface $P_{i-1}$ is pushed through some layer $\Delta_i$
and the new surface $P_i$ misses the new edge $e'$.  Then both
$P_{i-1}$ and $P_i$ have slope $e'$ and are either copies of
meridional disks $D_\mu$ or an annulus $A_\alpha$. Note that the
edge $e$ which was covered by $\Delta_i$ intersects the slope of
the new edge $e'$ twice, hence every boundary component of the
surface $P_i$ intersects $e$ twice.  Moreover, each edge $e'$
missed in this fashion determines a distinct edge $e$ that is
covered.  So although, the edge $e'$ is missed, an earlier edge
$e$ makes up for the deficit and we may count the edge $e'$ as if
it was intersected by each boundary component of $P_i$ .

\end{enumerate}

If $P_t$ is a meridional disk $D_\mu$ then it was obtained by
pushing through each layer $\Delta_i$; every surface $P_i$ is also
of type $D_\mu$. Every edge in $\T_1$ is intersected by the
original disk $P_1$.  See Figure \ref{f-stsurf}.  If any layered
edge $e'$ is missed then by (3) above, some earlier edge is
intersected twice. We can therefore count 1 intersection for each
edge.  We may also count an extra two intersections because $P_1 =
D_\mu$ hit edge $e_2$ three times, and we have only counted 1
(being hit three times means that it can not correspond to an edge
buried by reason (3) above). We have, $wt(D_\mu) \geq t+4$.

If $P_t$ is a trivial disk $D_\tau$, then each intermediate
surface $P_i$ is either of type $D_\tau$ or 2 copies of $D_\mu$.
In any event, any edge that is met is met twice.   Both $D_\tau$
and $2 D_\mu$ meet each edge of $\T_1$.  If any subsequent edge is
missed, it is due to reason (3) listed above, and the surface
$P_i$ is 2 copies of $D_\mu$.   Each boundary curve of $2 D_\mu$
intersects some earlier edge twice and we count 2 intersections
for each edge of the triangulation, $wt(D_\tau) \geq 2(t+2)$.

Suppose that $P_t$ is type $A_\alpha$.  If $P_1$ was an annulus
then $P_t$ was obtained by pushing through each layer and each
$P_i$ is of type $A_\alpha$. Then $P_1$ misses one of the initial
edges of $\T_1$, and by reason (3) above if any subsequent edge is
missed then an earlier edge was met twice.  For each edge except
one we count 2 intersections, one for each boundary component of
$A_\alpha$, $wt(P_t) \geq 2(t+1)$. If some $P_i$ is a surface of
type $D_\tau$ then its weight was computed in the previous
paragraph.   A subsequent edge can be missed only when the banding
quad is attached, or, after the banding quad is attached and due
to reason (3) above. Thus, we can count 2 for all but one of the
subsequent edges, $wt(P_\tau) \geq 2(t+1)$. The final case is that
$A_\alpha$ was obtained by attaching two bands to $2 D_\mu$ in a
single layer. In this case, the three initial edges are met twice
each.  Using (3) above, we count all subsequent edges for two
intersections  except for the edge corresponding to the bands
attached,  $wt(P_t) \geq 2 (t+1)$. Regardless, of the construction
we have the bound $wt(A_\alpha) \geq 2(t+1)$.
\end{proof}

Our understanding of layered triangulations allows us to construct
an example of compatible surface with complementary slopes.

\begin{ex} Consider the annulus of type
$A_\alpha$ contained in $\T_1$ pictured in Figure
\ref{f-stsurf}(3); it is disjoint from the edge $e_2$ of the
triangulation. Attach a new layer $\Delta_2$ at the edge $e_2$ and
use triangles to push the annulus through $\Delta_2$ to obtain an
annulus $A_1 \subset \T_2$, see Figure \ref{f-a1a2}.

\begin{figure}[h]
{\epsfxsize = 4.5 in \centerline{\epsfbox{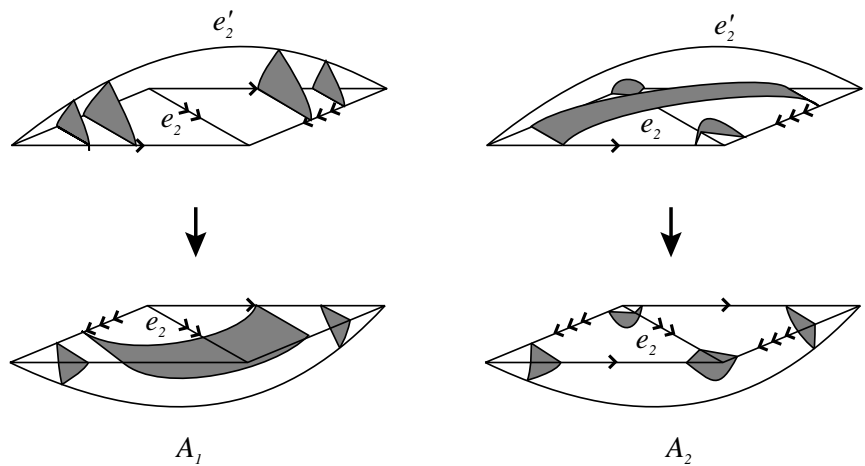}} }
\caption{Annuli $A_1$ and $A_2$ in $\T_2$.} \label{f-a1a2}
\end{figure}

Construct another annulus, $A_2 \subset \T_2$, by taking $D_\tau
\subset \P_1$ and attaching a banding quad, and two triangles in
$\Delta_2$.  These surfaces have distinct slopes, $\bdy A_1 =
(2,0,2)$ and $\bdy A_2 = (0,2,0)$ in normal coordinates with
respect to $e_1,e_2',e_3$.   Yet, their quads are in different
tetrahedra and the surfaces are thus compatible. Indeed, their
slopes are complementary, $\bdy A_1 + \bdy A_2 = \bdy (A_1 + A_2)
=(2,2,2) $ is two trivial curves. We can also see this by
constructing $A_1 + A_2$ by recombining the same pieces, see
Figure \ref{f-tt}.

\begin{figure}[h]
{\epsfxsize = 4.5 in \centerline{\epsfbox{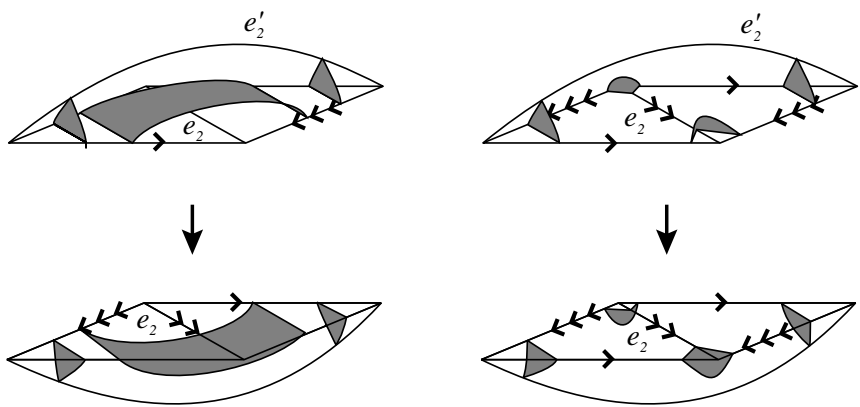}} } \caption{A
once-punctured torus and a trivial disk in $\T_2$.} \label{f-tt}
\end{figure}

Band the annulus to itself by attaching the quad and two triangles
and push $D_\tau$ through $\Delta_2$ by using all 4 triangles
types.  The former surface is a once punctured torus and the
latter a  vertex linking disk.  These normal surfaces have trivial
boundary and are disjoint.  The normal surface $A_1 + A_2$ is the
disjoint union of a vertex linking disk and a once punctured
torus.
\end{ex}


\section{Decision Problems in the Space of Dehn Fillings: Essential Surfaces.}
\label{s-essential}

In this section we consider the existence of certain interesting
surfaces in Dehn fillings of a knot-manifold $X$. Recall that a
surface $S$ properly embedded in a 3--manifold $M$ is {\it
compressible} if there is an embedded disk $D \subset M$ so that
$\bdy D \subset S$ is a non-trivial curve in $S$. If $S \neq
S^{2}$ is not compressible, we say $S$ is {\it incompressible}. A
properly embedded surface in $M$ is $\partial${\it -compressible}
if there is an embedded disk $D \subset M$ so that $\partial D = a
\cup b$, where $a$ and $b$ are arcs in $\partial D$, $a \cap b =
\partial a = \partial b$, $a \subset \partial M$, and $b \subset
S$ is not parallel into $\partial S$. If $S$ is not a disk and $S$
is not $\partial$-compressible, we say $S$ is {\it
$\partial$-incompressible}. A properly embedded surface is {\it
essential} if it is either a 2--sphere not bounding a 3--cell in
$M$, a disk not equivalent to a disk parallel into $\partial M$,
or it is two-sided, incompressible, $\partial$-incompressible and
not equivalent to a surface which is parallel into $\partial M$.
If $M$ contains an essential 2--sphere, then $M$ is said to be
{\it reducible}; otherwise $M$ is {\it irreducible}. The
3-manifold is {\it toroidal} if it contains an essential, embedded
torus; otherwise, it is {\it atoroidal}. Finally, a 3-manifold is
said to be a {\it Haken-manifold} if it is irreducible and
contains an embedded, incompressible surface. An irreducible
3-manifold with nonempty boundary is a Haken-manifold.

If a knot-manifold $X$ is given, we provide an algorithm to
determine precisely those slopes for which a Dehn filling is
reducible or those slopes for which a Dehn filling contains an
embedded, incompressible, two-sided surface.  Putting these
results together, we determine precisely those slopes for which a
Dehn filling is a Haken-manifold.  In the case of the
incompressible two-sided surface, our algorithm also may be used
to distinguish those slopes for which a Dehn filling is toroidal,
and those slopes for which a Dehn filling is fibered over $S^{1}$.

In \cite{g-l}, it is shown for $X$ an irreducible knot-manifold,
there are at most 3 reducible Dehn fillings of $X$; also, bounds
are given in \cite{go} for toroidal Dehn fillings when $X$ is
atoroidal. In \cite{wu}, it is shown that when $X$ contains an
embedded, essential surface, and when there is no embedded annulus
having one boundary a non-trivial curve in this surface and the
other a curve in $\partial X$, then there are at most 3 Dehn
fillings in which this surface compresses. Again, we comment that
we do not get such {\it a priori} global bounds; however, our
methods do give new proofs that for a given manifold bounds do
exist and for a given knot-manifold we give a method to compute
precisely the slopes for which these interesting phenomena happen.
The output of these algorithms will be a set of slopes described
by a finite set of points and/or by a {\it line} in the Dehn
filling space. If $\alpha$ is a slope on $\bdy X$ then the {\it
line of slopes determined by $\alpha$}, $L_\alpha$, is the
infinite set of slopes which intersect $\alpha$ precisely once,
i.e., $L_\alpha = \{\beta | \Delta(\alpha,\beta)\}$.

We begin this section by recalling results from normal surface
theory on deciding if a given manifold contains an essential
2--sphere or if it contains an embedded, incompressible, two-sided
surface.

\begin{thm}
\label{jo-jt} \cite{j-o,j-t} Let $\mathcal{T}$ be a triangulation
of the irreducible 3--manifold $X$.  Suppose $S$ is a normal
surface in ($X$,$\mathcal{T}$) that is least weight in its isotopy
class. If $S$ is two-sided, incompressible and
$\partial$-incompressible, then every rational point in the
carrier of $S$ in $\mathcal{P}$($X$,$\mathcal{T}$), is the
projective class of an embedded, incompressible and
$\partial$-incompres\-sible, two-sided, normal surface in
($X$,$\mathcal{T}$).
\end{thm}

The preceding theorem, in the case for embedded closed,
incompressible, two-sided surfaces, is one of the main results of
\cite{j-o}.  The theorem was extended to include embedded
incompressible and $\partial$-incompressible surfaces (extended to
the bounded case) in \cite{j-t}. We need analogous results for
embedded, essential, normal 2--spheres and for embedded,
incompressible, two-sided, closed, normal surfaces when the
3-manifold may not be irreducible. The desired result for
2-spheres follows from recent work of W.~Jaco and L.~Reeves
\cite{j-re} where the assumption on the 2-sphere is that it is an
absolute least weight, embedded, essential, normal 2--sphere,
Theorem \ref{T-spherecarrier}. Similar results appear in
\cite{j-t}. The latter case, involving incompressible surfaces,
requires modification of the proof in \cite{j-o} and consideration
of a possibly larger equivalence class of embedded,
incompressible, two-sided, least weight, normal surfaces. The
result we need is given in Theorem \ref{T-esscarrier}, below. The
proof of Theorem \ref{T-esscarrier}, including the case with
nonempty boundary and embedded, incompressible and
$\partial$-incompressible surfaces, can be obtained from straight
forward modification of the proof in \cite{j-o}.

\begin{thm}
\label{T-spherecarrier} \cite{j-re} Let $\mathcal{T}$ be a
triangulation of the 3--manifold $M$. If $\Sigma$ is a least
weight, embedded, essential, normal 2--sphere in
($M$,$\mathcal{T}$), then every rational point in the carrier of
$\Sigma$ in $\mathcal{P}$($M$,$\mathcal{T}$), is the projective
class of a normal surface each component of which is an embedded,
essential, normal 2--sphere in ($M$,$\mathcal{T}$).
\end{thm}

We have given the conclusion of the preceding theorem to allow for
the possibility that some projective class in the carrier of
$\Sigma$ in $\mathcal{P}$($M$,$\mathcal{T}$) may have no
representative that is connected. By using projective classes we
also have the possibility that some representative may be an
embedded projective plane; however, its double, also a
representative of the same projective class, will then be a
2--sphere.

In what follows, we use disk swapping, which was defined in
Section \ref{s-normal}, for equivalence between surfaces. Hence,
if two surfaces $S$ and $S'$ are isotopic, they are equivalent.
Being equivalent and isotopic are the same when the ambient
manifold is irreducible. The concept of ``disk swapping" applies
to ``$\bdy$-compressing disks" as well and is a necessary
extension of this concept in the case that the manifold $X$ has
boundary and the surfaces in question are
$\partial$-incompressible. Furthermore, any two embedded 2-spheres
are equivalent via disk-swapping and so an embedded, essential,
normal 2-sphere that is least weight in its equivalence class is a
least weight, embedded, essential, normal 2-sphere. Note that the
word essential is crucial, as a least weight normal 2--sphere may
not be essential and a least weight 2--sphere is not normal and
has zero weight.

\begin{thm}
\label{T-esscarrier} Let $\mathcal{T}$ be a triangulation of the
3--manifold $M$. Suppose $S$ is an embedded normal surface in
($M$,$\mathcal{T}$) that is least weight in its equivalence class.
If $S$ is two-sided, incompressible, and
$\partial$-incompressible, then every rational point in the
carrier of $S$ in $\mathcal{P}$($M$,$\mathcal{T}$), is the
projective class of an embedded, incompressible,
$\partial$-incompressible, two-sided, normal surface in $M$.
\end{thm}

The following theorem is the primary tool for many of the results
of this section. We obtained the results of this section prior to
discovering this theorem.  While it simplifies our earlier proofs,
its major appeal, however, is that of greatly simplifying the
algorithms and exhibiting the fundamental roll of the topology of
$X$ to that of $X(\alpha)$. Specifically, using special one-vertex
triangulations for Dehn fillings, as in \cite{j-r3}, which fix a
triangulation $\mathcal{T}$ of $X$ for all the Dehn fillings of
$X$, we show that $X(\alpha)$ contains an essential surface if and
only if one of the vertex-solutions of
$\mathcal{P}$($X$,$\mathcal{T}$) is an embedded, essential surface
in $X$ and is either closed or ``caps off" to a surface which is
essential in $X(\alpha)$. It follows that there are a finite
number of surfaces in $X$ (all computable) which determine the
existence (or nonexistence) of an essential surface in all Dehn
fillings of $X$.

\begin{thm}
\label{T-essinx} Suppose $X$ is a knot-manifold and $\mathcal{T}$
is a triangulation of $X$ that restricts to a one-vertex
triangulation of $\bdy X$. If $X(\alpha)$ contains an embedded,
essential surface, then there is an embedded, essential, normal
surface $G$ in ($X$,$\mathcal{T}$) such that the projective class
of $G$ is a vertex-solution of $\mathcal P$($X$,$\mathcal T$), the
boundary slope of $G$ is $\alpha$ (if $\partial{G} \neq
\emptyset$), and $G(\alpha)$ is an embedded, essential, normal
surface in ($X(\alpha)$,$\mathcal{T}(\alpha)$).

In fact, if $X(\alpha)$ is reducible, then a vertex-solution $S$
of $\mathcal P$($X$,$\mathcal T$) must be either an embedded,
essential 2--sphere or planar surface and $S(\alpha)$ is an
embedded, essential 2--sphere in $X(\alpha)$; if $X(\alpha)$
contains an embedded, incompressible, two-sided surface, then a
vertex-solution $F$ of $\mathcal P$($X$,$\mathcal T$) must be an
embedded, essential, non-planar surface and $F(\alpha)$ is an
embedded, incompressible, two-sided surface in $X(\alpha)$; and,
in the latter case, if $X(\alpha)$ contains an embedded,
incompressible torus, then a vertex-solution $T$ of
$\mathcal{P}(X,\mathcal{T})$ must be an embedded, essential torus
or punctured-torus and $T(\alpha)$ is an embedded, incompressible
torus in $X(\alpha)$, and if $X(\alpha)$ fibers over $S^{1}$, then
a vertex-solution $F$ of $\mathcal{P}(X,\mathcal{T})$ must be an
embedded, essential, two-sided surface and $F(\alpha)$ is a fiber
in a fibration of $X(\alpha)$ over $S^{1}$.
\end{thm}
\begin{proof}

We are given that $X(\alpha)$ contains an embedded, essential
surface. Hence, $X(\alpha)$ contains an embedded, essential
2--sphere or an embedded, incompressible, two-sided surface or
both. We have organized the proof to handle the general situation;
however, we indicate the specific considerations and give the
details needed to arrive at the special conclusions given in the
second part of the statement of the theorem.

Suppose $\Gamma$ is an embedded, essential surface in $X(\alpha)$.
Among all essential surfaces in $X(\alpha)$ that are equivalent
with $\Gamma$ (recall that equivalence means equivalent via
disk-swapping and isotopy) consider those that meet $V_\alpha$ in
the smallest number of components. We can find such a surface that
meets $V_\alpha$ in a collection of pairwise disjoint copies of
the meridional disk or not at all. Furthermore, assuming notation
has been chosen so that $\Gamma$ is itself such a surface, then
for $G$ = $X \cap \Gamma$, $G$ is an embedded, essential surface
with boundary slope $\alpha$ ( $\partial G \neq \emptyset$) or a
closed essential surface in $X$. There is no loss in generality to
assume that $G$ is also normal in ($X$,$\mathcal T$).

\renewcommand{\theenumi}{\roman{enumi}}

Having made these observations, it follows that there is an
embedded, essential, normal, surface $G$ with boundary slope
$\alpha$ (if $\partial G \neq \emptyset$) in ($X$,$\mathcal T$)
such that:
\begin{enumerate}
\item $G(\alpha)$ is defined and is equivalent to $\Gamma$ in
$X(\alpha)$,
\item $G(\alpha)$ meets $V_\alpha$ in the minimal number of
components among all embedded, essential surfaces in
$X(\alpha)$ that are equivalent to
$\Gamma$, and
\item if $G'$ is an embedded, essential, normal surface
that is either closed or has boundary slope $\alpha$ in
($X$,$\mathcal T)$ and $G'(\alpha)$ satisfies i and ii, then
$wt(G) \leq wt(G')$; i~.e~., $G$ is least weight in $(X,\T)$ with
respect to conditions i and ii.
\end{enumerate}

It follows from Theorem \ref{T-esscarrier} above that every
surface with projective class in the carrier of $G$ in $\mathcal
P$($X$,$\mathcal T$),  is an embedded, essential, normal surface
in ($X$, $\mathcal T$); furthermore, such a surface, if it has
boundary, has essential boundary and, therefore, by Corollary
\ref{carrier-slopes}, each boundary component has slope $\alpha$.
Hence, all the surfaces with projective class in the carrier of
$G$ are either closed in $X$ or cap off with meridional disks in
$V_\alpha$ to give closed surfaces in $X(\alpha)$. In particular,
the surfaces with projective classes at the vertices of the
carrier of $G$ cap off to closed surfaces in $X(\alpha)$. What we
need to show is that surfaces in the carrier of $G$ in $\mathcal
P(X, \mathcal T)$ cap off to essential surfaces in $X(\alpha)$. In
addition, to achieve the specific conclusions of the second part
of the theorem, we need to show that if $G$ is a 2--sphere or is
planar, then the surfaces with projective classes at the vertices
of the carrier of $G$ are 2--spheres or are planar and cap off to
essential 2--spheres; if $G$ is non-planar, the surfaces with
projective classes at the vertices of the carrier of $G$ are
non-planar and cap off to incompressible surfaces; and if $G$ is a
torus or punctured torus, then the surfaces with projective
classes at the vertices of the carrier of $G$ are either tori or
punctured tori and cap off to incompressible surfaces; and, finally, if
$G(\alpha)$ is a fiber in a fibration over $S^{1}$, then the
surfaces with projective classes at the vertices of the carrier of
$G$ cap off to fibers in fibrations over $S^{1}$.

The triangulation $\T$ induces a one-vertex triangulation on
$\partial X$ and so, a one-vertex triangulation on $\partial
V_\alpha$. By Theorem \ref{exist-layered} and Proposition
\ref{layered-planars} there is a layered, one-vertex triangulation
of $V_\alpha$, extending this triangulation on $\partial V_\alpha$
so that any planar, normal surface in $V_\alpha$ has weight $ \geq
wt(G)$.  We can extend the triangulation $\mathcal{T}$ to a
triangulation, say $\mathcal{T(\alpha)}$, of $X(\alpha)$ using
such a layered, one-vertex triangulation of $V_\alpha$.  If $F$ is
a normal surface whose projective class is in the carrier of $G$
then we may write $k G = F + F'$, where $F'$ is some other normal
surface whose projective class is in the carrier of $G$.  Then $F$
and $F'$ both have slope $\alpha$ and cap off to normal surfaces
$F(\alpha)$ and $F'(\alpha)$ in $(X(\alpha),\T(\alpha))$.
Furthermore, we may write $k G(\alpha) = F(\alpha) + F'(\alpha)$
so it follows that $F(\alpha)$ and $F'(\alpha)$ are surfaces whose
projective classes are in the carrier of $G(\alpha)$ in
$\P(X(\alpha),\T(\alpha))$.  We want to show that the surfaces in
the carrier of $G$ in $\mathcal P$($X$,$\mathcal{T}$) cap off to
essential surfaces in $X(\alpha)$. Furthermore, we will arrive at
such a conclusion using Theorem \ref{T-esscarrier} above in
$X(\alpha)$ and  results from \cite{j-re} (generalizing \cite{j-o,
j-t}) which will give us the special conclusions of the second
part of the theorem.

We claim $G(\alpha)$ is least weight in its equivalence class in
($X(\alpha), \mathcal{T}(\alpha)$).  For suppose $\Gamma'$ is a
normal surface equivalent to $G(\alpha)$ in ($X(\alpha),
\mathcal{T}(\alpha)$) and $wt(\Gamma') < G(\alpha)$. It was
observed in \cite{j-r3} that each component of $\Gamma' \cap
V_\alpha$ must be a (normal) planar surface in $V_\alpha$. Now, by
Proposition \ref{layered-planars}, each component of $\Gamma' \cap
V_\alpha$ is either a normal annulus or a normal disk ($V_\alpha$
has a layered triangulation) and therefore by the choice of the
layered triangulation of $V_\alpha$ each component of $\Gamma'
\cap V_\alpha$ has weight $\geq$ wt($G$).

By our choice of $G$, it follows that $\Gamma'$ meets $V_\alpha$
in at least as many meridional disks as $G(\alpha)$. If the number
of components in $\Gamma' \cap V_\alpha$ were more, then by the
choice of $\mathcal{T}(\alpha)$, $wt(\Gamma') \geq wt(G(\alpha))$.
Hence, we must have that $\Gamma' \cap V_\alpha$ has precisely the
same number of components as $G(\alpha) \cap V_\alpha$ and each
component of intersection is a meridional disk; for otherwise the
number of components of $\Gamma' \cap V_\alpha$ could be reduced,
contradicting our choice of $G$. Let $G' = \Gamma' \cap X$. Thus
$G'$ is an embedded, essential, normal surface in ($X$,$\mathcal
T)$ with boundary slope $\alpha$ and $G'(\alpha)$ satisfies i and
ii above. So, $wt($G$) \leq wt(G')$ and, therefore, $wt(G(\alpha))
\leq wt(G'(\alpha)) = wt(\Gamma')$. Hence, $G(\alpha)$ is a least
weight, embedded, essential, normal surface in ($X(\alpha)$,
$\mathcal{T}(\alpha)$).

If $G(\alpha)$ is an essential 2--sphere, then by \cite{j-re},
Theorem \ref{T-spherecarrier} above, every rational point in the
carrier of $G(\alpha)$ in
$\mathcal{P}$($X(\alpha)$,$\mathcal{T}(\alpha$)), is the
projective class of a normal surface each component of which is an
embedded, essential, normal 2--sphere in
($X(\alpha)$,$\mathcal{T}(\alpha)$). However, any normal surface
in ($X$,$\mathcal{T}$) with projective class in the carrier of $G$
in $\mathcal{P}$($X$,$\mathcal{T}$) can be capped off to a normal
surface in ($X(\alpha)$,$\mathcal{T}(\alpha)$) whose projective
class is in the carrier of $G(\alpha)$ in
$\mathcal{P}$($X(\alpha)$,$\mathcal{T}(\alpha)$). It follows that
any normal surface in ($X$,$\mathcal{T}$) with projective class in
the carrier of $G$ in $\mathcal{P}$($X$,$\mathcal{T}$) is the
projective class of a normal surface each component of which is an
embedded, essential, 2--sphere or planar surface in
($X$,$\mathcal{T}$) and caps off to an embedded, essential,
2--sphere in $X(\alpha)$; in particular, this is true for any
normal surface whose projective class is a vertex-solution of the
carrier of $G$ in $\mathcal{P}$($X$,$\mathcal{T}$).

If $G(\alpha)$ is an embedded, incompressible, two-sided surface,
then by Theorem \ref{T-esscarrier} above every rational point in
the carrier of $G(\alpha)$ in
$\mathcal{P}$($X(\alpha)$,$\mathcal{T}(\alpha$)), is the
projective class of an embedded, incompressible, two-sided, normal
surface in $X(\alpha)$. It follows that any normal surface in
($X$,$\mathcal{T}$) with projective class in the carrier of $G$ in
$\mathcal{P}$($X$,$\mathcal{T}$) is the projective class of an
embedded, essential, non-planar surface in ($X$,$\mathcal{T}$) and
caps off to an embedded, incompressible, two-sided, normal surface
in $X(\alpha)$; in particular, this is true for any normal surface
whose projective class is a vertex-solution of the carrier of $G$
in $\mathcal{P}$($X$, $\mathcal{T}$). If $G(\alpha)$ is a torus,
then every surface in the carrier of $G(\alpha)$ in
$\mathcal{P}$($X(\alpha)$,$\mathcal{T}(\alpha$)), is the
projective class of an embedded, incompressible, two-sided, normal
torus. Hence, any normal surface in ($X$,$\mathcal{T}$) with
projective class in the carrier of $G$ in
$\mathcal{P}$($X$,$\mathcal{T}$) is the projective class of an
embedded, essential, punctured torus or torus and caps off to an
embedded, incompressible, two-sided, normal torus in $X(\alpha)$.
Finally, if $G(\alpha)$ is a fiber in a fibration of $X(\alpha)$
over $S^{1}$, then every surface in the carrier of $G(\alpha)$ in
$\P(X(\alpha),\T(\alpha))$ is the projective class of a fiber in a
fibration of $X(\alpha)$ over $S^{1}$ \cite{jac}.  So, any normal
surface in $(X,\T)$ with projective class in the carrier of $G$ in
$\P(X,\T)$ caps off to a fiber in a fibration of $X(\alpha)$ over
$S^{1}$.  This completes the proof.
\end{proof}

\subsection{Reducible Manifolds in Dehn Surgery Space.}
Given a knot-manifold $X$, we consider the problem of determining
precisely those slopes $\alpha$ for which the Dehn filling
$X(\alpha)$ is reducible. We consider two distinct situations. The
first is when the knot-manifold $X$ is irreducible. In this
situation most (all but a finite number) Dehn fillings are
irreducible. If the knot-manifold $X$ is reducible, then we show
that only in very special situations does one get an irreducible
Dehn filling.

From above we have that a Dehn filling $X(\alpha)$ is reducible if
and only if at least one of a finite number of constructable
planar surfaces in $X$ leads to an essential 2--sphere in
$X(\alpha)$ or $X$ itself contains a constructable essential
2-sphere that remains essential in $X(\alpha)$. However, for an
algorithm to decide these issues, we need a result of H.
Rubinstein, which provides a method to recognize if a given normal
2--sphere is essential (\cite{rubinstein} and \cite{thompson}, a
solution to the 3--sphere recognition problem).

\begin{thm}
\label{T-recognize-1}\cite{rubinstein,thompson} Suppose $\T$ is a
triangulation of the 3--manifold $M$. Given a normal 2--sphere
$\Sigma$ in $(M,\T)$ it can be decided if $\Sigma$ bounds a
3--cell in $M$.
\end{thm}

\begin{thm}
\label{T-recognize-2} \cite{rubinstein,thompson} Given a compact
3--manifold $M$, it can be decided if $M$ is irreducible;
furthermore, \cite{j-t, j-re} if $M$ is not irreducible, there is
an algorithm to construct an irreducible (a minimal irreducible or
prime) decomposition of $M$.
\end{thm}

It is known for $X$ an irreducible knot-manifold there are only
finitely many slopes $\alpha$ for which $X(\alpha)$ is reducible;
\cite{wu} showed that if $\alpha$ and $\beta$ are both slopes for
which $X(\alpha)$ and $X(\beta)$ are reducible then $\Delta
(\alpha,\beta) \leq 2$. Later in \cite{g-l}, it was shown that
$\Delta (\alpha,\beta) \leq 1$ holds. Hence, there is a global
finite bound; namely, $X(\alpha)$ is reducible for at most 3
slopes. We do not get a global bound but do get a new proof that
the number is finite for any  knot-manifold $X$ and show that
there is an algorithm to determine precisely those slopes $\alpha$
for which $X(\alpha)$ is reducible.

\begin{thm}
\label{T-reducibleslopes} Given an irreducible knot-manifold $X$,
there is an algorithm to determine precisely those slopes $\alpha$
for which the Dehn filling $X(\alpha)$ is reducible; in
particular, it follows that there are only finitely many slopes
$\alpha$ for which $X(\alpha)$ is reducible.
\end{thm}

\begin{proof}
We assume $X$ is given via a triangulation $\mathcal{T}$ that
restricts to a one-vertex triangulation on $\bdy X$. By Theorem
\ref{T-essinx}, $X(\alpha)$ is reducible if and only if there is a
vertex-solution $S$ of $\mathcal P$($X$, $\mathcal T$) that is
planar ($X$ is assumed to be irreducible) and $S(\alpha)$ is an
embedded, essential 2--sphere in $X(\alpha)$. Let $\mathcal{A} =
\{\alpha_{1},\ldots,\alpha_{n}\}$ be the set of boundary slopes of
embedded, planar, normal surfaces with projective classes at a
vertex of $\mathcal P$($X$,$\mathcal T$).  If $\mathcal{A} =
\emptyset$, then $X(\alpha)$ is irreducible for all $\alpha$. If
$\mathcal{A} \neq \emptyset$, then $X(\alpha)$ can only be
reducible for $\alpha = \alpha_{i}$ for some $i, 1\leq i\leq m$.
So, suppose $P_{i_{1}}, \ldots, P_{i_{m}}$ is the set of all
embedded, connected planar, normal surfaces with projective
classes at a vertex of $\mathcal P$($X$,$\mathcal T$) having slope
$\alpha_{i}$; $X(\alpha_{i})$ will be reducible if and only if
some $P_{i_{j}}$ caps off to an essential 2-sphere in
$X(\alpha_{i})$. This can be checked by the algorithm of Theorem
\ref{T-recognize-1}, stated above.

It follows that there are at most a finite number of slopes
$\alpha$ such that $X(\alpha)$ is reducible; and these slopes are
among the boundary slopes of embedded, planar, normal surfaces
with projective classes at a vertex of $\mathcal P$($X$,$\mathcal
T$).
\end{proof}

In the hypothesis of Theorem \ref{T-reducibleslopes}, it is
assumed that it is known that the knot-manifold $X$ is
irreducible. Of course, Theorem \ref{T-recognize-2} tells us that
it can be decided if a 3--manifold is irreducible; so, the issue
is, in the case $X$ is reducible, can we decide those slopes
$\alpha$ for which $X(\alpha)$ is reducible (or, more accurately,
those slopes $\alpha$ for which $X(\alpha)$ is irreducible; since
the generic case when $X$ is reducible, is for $X(\alpha)$ to be
reducible). We can do this; however, we need the results of
Section \ref{s-heegaard} for the complete proof; in particular, we
need Theorem \ref{s3algo} which provides an algorithm to decide
precisely those slopes for which a Dehn filling gives the
3-sphere.

\begin{thm}
\label{T-irreducibleslopes} Given a reducible knot-manifold $X$,
there is an algorithm to determine precisely those slopes $\alpha$
for which the Dehn filling $X(\alpha)$ is irreducible.
\end{thm}

\begin{proof}
We are assuming the knot-manifold $X$ is reducible. By Theorem
\ref{T-reducibleslopes} there is an algorithm to construct an
irreducible decomposition of $X$. If $X$ contains a
non-separating, embedded 2--sphere the algorithm will find one and
it follows that $X(\alpha)$ will be reducible for all $\alpha$.
So, we may assume every 2--sphere embedded in $X$ separates $X$.
However, if $X$ contains two independent, separating 2--spheres
(i.e., $X$ contains disjoint, essential 2--spheres $S_{1}$ and
$S_{2}$ where $S_{1} \cup S_{2}$ does not bound a product $S^{2}
\times [0,1]$), then, again, the algorithm will construct such a
pair and it follows that $X(\alpha)$ is reducible for all slopes
$\alpha$. Thus, the only possibility for $X(\alpha)$ to be
irreducible when $X$ is reducible is that $X$ has a separating,
essential 2--sphere $S$, each component of $\widehat{X_S}$ (the
manifold obtained from $X$ by splitting at $S$ and capping off
each 2--sphere boundary component with a 3--cell) is irreducible,
and, for notation chosen so that $M$ is the component of
$\widehat{X_S}$ containing $\partial{X}$, $M$ is a knot-manifold
in $S^{3}$, i.e., $M$ embeds in $S^{3}$.

Again by Theorem \ref{T-reducibleslopes}, the algorithm will find
a separating, essential 2--sphere $S$ in $X$ and thus determine
the knot-manifold $M$, as above. Now by Theorem \ref{s3algo}, we
can decide if the knot-manifold $M$ embeds in $S^{3}$ and
determine precisely those slopes $\alpha$ for which $M(\alpha)$ is
homeomorphic with $S^{3}$. If $M$ does not embed in $S^{3}$, then
for all slopes $\alpha$, $X(\alpha)$ is reducible. If the
knot-manifold $M$ embeds in $S^{3}$ and it is not a solid torus,
there is only one slope $\alpha$ \cite{g-l} for which the Dehn
filling $M(\alpha)$ is homeomorphic with $S^{3}$, and the
algorithm finds this slope. If the knot-manifold $M$ is a solid
torus and $\mu$ is the meridional slope (the algorithm of Theorem
\ref{T-decideincomp}, for example, finds the meridional slope),
then for every slope $\alpha$ with $\Delta(\alpha,\mu) \leq 1$,
the line $L_{\mu}$, $X(\alpha)$ is irreducible.
\end{proof}

It follows from the proof of the previous theorem that whenever
the knot-manifold $X$ is reducible, one of the following holds:
$X(\alpha)$ is reducible for every slope $\alpha$; or $X$ is a
connected sum of a non-trivial knot-manifold in $S^{3}$ and an
irreducible manifold and there is precisely one computable slope
$\alpha$ for which $X(\alpha)$ is irreducible; or $X$ is a
connected sum of a solid torus and an irreducible manifold and
there is a computable line of slopes for which $X(\alpha)$ is
irreducible.

\vspace{.25 in}

\noindent {\bf Algorithm R.} {\it Given a knot-manifold $X$,
determine precisely those slopes $\alpha$ for which the Dehn
filling $X(\alpha)$ is reducible.}

\vspace{.2 in}

{\bf Step 1.} We assume the knot-manifold $X$ is given via a
triangulation. Endow $X$ with a triangulation $\mathcal{T}$ that
restricts to a one-vertex triangulation on $\bdy X$. (An algorithm
to do this is given in \cite{j-r3}.)

{\bf Step 2.}  Compute the vertices of
$\mathcal{P}$($X$,$\mathcal{T}$); i.e. find all embedded normal
surfaces whose projective class is a vertex of
$\mathcal{P}$($X$,$\mathcal{T})$.

{\bf Step 3.} Determine if $X$ is reducible (if $X$ has an
essential, embedded 2--sphere). (Recall that $X$ has an essential,
embedded 2--sphere if and only if there is an essential, embedded
normal 2--sphere in ($X$,$\mathcal{T}$) whose projective class is
at a vertex of $\mathcal{P}$($X$,$\mathcal{T}$) \cite{j-t, j-re};
furthermore, it can be decided if a given embedded, normal
2--sphere is essential \cite{rubinstein} and if a finite, pairwise
disjoint collection of normal 2--spheres is independent
\cite{rubinstein, j-re}.) Begin the algorithm to construct a
minimal irreducible decomposition of $X$ \cite{j-t, j-re}.

If an embedded, non-separating 2--sphere is found, then for
every slope $\alpha$, the manifold $X(\alpha)$ will be reducible and the
algorithm terminates.

If a pair of independent, embedded, normal 2-spheres
is found, then for every slope $\alpha$, the manifold $X(\alpha)$ will be
reducible and the algorithm terminates.

So, the only possibility left, if $X$ is reducible (the
irreducible decomposition is not empty), is that there is one
essential (separating) 2--sphere in the irreducible decomposition
of $X$. If this is the situation and we let $S$ denote such a
normal 2--sphere and let $M$ denote the component of
$\widehat{X_S}$ that contains $\partial X$, then $M$ is a knot
manifold with $\partial M = \partial X$. We wish to determine
precisely those slopes $\alpha$ for which $M(\alpha)$ is the
3--sphere. The algorithm in Section 6 (Theorem \ref{s3algo}),
which includes the possibility that $M$ is a solid torus, can be
used to determine such slopes $\alpha$, either precisely one slope
or a line of slopes and the algorithm terminates.

If the irreducible decomposition is empty, then go to the next step.

{\bf Step 4.} List the vertices of
$\mathcal{P}$($X$,$\mathcal{T}$) that correspond to the projective
classes of planar, normal surfaces in ($X$,$\mathcal{T}$). (Recall
that if the knot-manifold $X$ is irreducible and $X(\alpha)$ is
reducible, then a vertex-solution $S$ of $\mathcal
P$($X$,$\mathcal T$) must be planar with $S(\alpha)$ being an
embedded, essential 2--sphere in $X(\alpha)$.) If there are none,
then for every slope $\alpha$, $X(\alpha)$ is irreducible.
Otherwise, let $\{S_{1}, \ldots, S_{k}\}$ be all the planar normal
surfaces whose projective class is a vertex of $\mathcal
P$($X$,$\mathcal T$). Calculate the boundary slope of each
$S_{i}$, $1 \leq i \leq k$; let $\{\alpha_{1}, \ldots,
\alpha_{k}\}$ be the set of slopes where $\alpha_{i}$ is the
boundary slope of $S_{i}$, $1 \leq i \leq k$.

{\bf Step 5.} Determine if $S_{i}(\alpha_{i})$ is essential in
$X(\alpha_{i})$, $1 \leq i \leq k$, using the algorithm of Theorem
\ref{T-recognize-1}. If $S_{i}(\alpha_{i})$ is essential in
$X(\alpha_{i})$, then $X(\alpha_{i})$ is reducible. If $X$ is
irreducible, the finite list of slopes $(\alpha_{i})$ for which
$X(\alpha_{i})$ is reducible is precisely the set of slopes
$\alpha$ for which $X(\alpha)$ is reducible.

This completes Algorithm R.

\vspace{.25 in}

\subsection{Haken-Manifolds in Dehn Surgery Space.}

In this section we provide an algorithm to determine precisely
those manifolds in the space of Dehn fillings that are
Haken-manifolds. The main problem, after the previous section, is
given a knot-manifold $X$ to determine precisely those slopes
$\alpha$ for which the Dehn filling $X(\alpha)$ contains an
embedded, closed, incompressible, two-sided surface. The problem
splits in a manner similar to that in the last section. If the
knot-manifold $X$ does not contain an embedded, closed, essential,
two-sided surface, then the generic Dehn filling is not expected
to contain an embedded, closed, incompressible, two-sided surface.
We give an algorithm to determine precisely those slopes $\alpha$
for which $X(\alpha)$ contains an embedded, incompressible,
two-sided surface. We obtain a new proof that there are only
finitely many slopes $\alpha$ for which $X(\alpha)$ contains such
a surface.

On the other hand, if $X$ contains an embedded, closed, essential,
two-sided surface, then the generic Dehn filling is expected to
contain an embedded, closed, incompressible, two-sided surface. We
give, in this case, an algorithm to determine precisely those
slopes $\alpha$ for which $X(\alpha)$ does {\it not} contain an
embedded, closed, incompressible, two-sided surface. Of
independent interest in this proof is an algorithm to show that
given an embedded, closed, two-sided surface $S$, we can find
precisely those slopes $\alpha$ for which $S$ compresses in
$X(\alpha)$. From this and the results in the previous section it
is easy to determine precisely those slopes $\alpha$ for which the
Dehn filling $X(\alpha)$ is a Haken-manifold.

First, we recall a result due to Haken \cite{haken}. A proof
appears in \cite{j-o} for handle-decompositions; the proof in the
case of triangulations requires some modification to the proof for
handle decompositions. Proofs for triangulations and refinements
in the algorithm appear in \cite{j-t, j-re}.

\begin{prop}\cite{haken}
\label{T-decideincomp} Let $M$ be a 3--manifold with triangulation
$\mathcal{T}$. Given a two-sided, normal surface $F$ in
($M$,$\mathcal{T}$), there is an algorithm to decide if $F$ is
incompressible in $M$.
\end{prop}

\begin{thm}
\label{T-surfacesincompress} Given a knot-manifold $X$ which does
not contain an embedded, closed, essential, two-sided surface
there is an algorithm to determine precisely those slopes $\alpha$
for which $X(\alpha)$ contains an embedded, closed,
incompressible, two-sided surface; in particular, it follows that
there are at most finitely many slopes $\alpha$ for which
$X(\alpha)$ contains such a surface.
\end{thm}

\begin{proof} Suppose for the slope $\alpha$, the Dehn filling $X(\alpha)$
contains an embedded, incompressible, two-sided surface. It
follows from Theorem \ref{T-essinx} that there is an embedded,
essential, two-sided surface $S$ of ($X$,$\mathcal{T}$) whose
projective class in $\mathcal P$($X$, $\mathcal T$) is a vertex,
the boundary slope of $S$ is $\alpha$ ($\partial{S} \neq
\emptyset$, since by assumption $X$ does not contain an embedded,
closed, essential, two-sided surface), and $S(\alpha)$ is an
embedded, incompressible, two-sided surface in $X(\alpha)$. So,
the only slopes $\alpha$ for which it is possible that $X(\alpha)$
contain an embedded, closed, incompressible, two-sided surface are
among a subset of boundary slopes coming from embedded, essential,
non-planar, normal surfaces in ($X$,$\mathcal{T}$) whose
projective classes are at vertices of $\mathcal P$($X$,$\mathcal
T$). This is a finite set.

Let $\{S_{1},\ldots,S_{n}\}$ be the set of embedded, non-planar,
normal surfaces in ($X$,$\mathcal{T}$) whose boundary consists of
only non-trivial curves in $\partial X$ and whose projective class
in $\mathcal P$($X$,$\mathcal T$) is a vertex. Let $\alpha_{i}$ be
the slope of $X_{i}, 1\leq i\leq n.$ We check if
$S_{i}(\alpha_{i})$ is incompressible in $X(\alpha_{i})$. It is
precisely those slopes $\alpha_{i_{j}}$ for which
$S_{i_{j}}(\alpha_{i_{j}})$ is incompressible in
$X(\alpha_{i_{j}})$ that satisfy the conclusion of the theorem.
\end{proof}

We now consider the situation when the knot-manifold $X$ contains
an embedded, closed, essential, two-sided surface. First, we make
some notational conventions regarding a planar surface. If $D$ is
a disk and $p_{1}, \ldots, p_{n}$ are points in the interior of
$D$, we call the planar surface obtained from $D$ by removing the
interior of a small regular neighborhood about each point $p_{i}$,
$i \leq i \leq n$, a {\it punctured-disk}. In this situation, if
$P = D \setminus \bigcup_{1}^{n} Int N(p_{i})$, where $N(p_{i})$
is a small regular neighborhood of $p_{i}$ in $D$, we call the
boundary component $\partial{D}$ of $P$ the {\it boundary of} $P$,
written $bdry(P)$, and the boundary components
$\partial{N(p_{i})}$ the {\it punctures of} $P$. Similarly, if $A$
is an annulus and the $p_i$'s and $N(P_i)$'s are defined the same,
then we call $Q = A \setminus \bigcup_{i=1}^n Int N(p_i)$ a
punctured annulus and call the boundary components of $A$ the
boundary  of $Q$, denoted $bdry(Q)$, and the boundary components
of $\bdy N(p_i)$ the punctures of $Q$.  In this way we distinguish
boundary components of such planar surfaces.

We find that the generic case, when the knot-manifold $X$ contains
an embedded, essential, closed, two-sided surface, is for
$X(\alpha)$ to contain an embedded, closed, incompressible,
two-sided surface. In particular, Wu  has shown \cite{wu} that if
$X$ contains an embedded, essential, closed, two-sided surface $S$
and there is no annulus from $S$ to $\partial X$, then $S$ will
compress in the Dehn filling $X(\alpha)$ for at most $3$ slopes
$\alpha$. We are able to show that given an embedded, closed,
two-sided, normal surface $S$, there is an algorithm to determine
precisely those slopes $\alpha$ for which $S$ compresses in
$X(\alpha)$. Again, our techniques give finiteness in the case
considered by Wu but do not give similar global bounds; and we
obtain complete answers when there is an annulus embedded in $X$
having one boundary a non-trivial curve in $S$ and the other in
$\bdy X$ (see \cite{cgls}).

First, we have the following lemma which has independent interest.

\begin{lem}
\label{closedess} Let $\T$ be a triangulation of the knot-manifold
$X$ that restricts to a one-vertex triangulation of $\bdy X$.
Given a closed, two-sided, normal surface $S$ in $(X,\T)$, there
is an algorithm to decide precisely those slopes $\alpha$ for
which $S$ is incompressible in $X(\alpha)$.
\end{lem}
\begin{proof}

 Given $S$ normal in $(X,\T)$, there are algorithms to
decide if $S$ is incompressible in $X$ (\cite{haken}, see
Proposition \ref{T-decideincomp}) and if $S$ is equivalent to a
boundary parallel surface \cite{j-t}. If $S$ compresses in $X$ or
$S$ is equivalent to a boundary-parallel surface, then $S$ will
compress in $X(\alpha)$ for every $\alpha$. Hence, we may assume
$S$ is essential in $X$ (and $S$ is not $S^{2}$). The surface $S$
will be normal in $(X(\alpha),\T(\alpha))$ for all $\alpha$.

We consider two possibilities:

\renewcommand{\theenumi}{\roman{enumi}}

\begin{enumerate}
\item there is no annulus in $X$ having one boundary component a
non-trivial curve in $S$ and the other a curve in $\bdy X$, or
\item there is an annulus in $X$ having one boundary component a
non-trivial curve in
$S$ and the other a curve in $\bdy X$.
\end{enumerate}

Split $X$ at $S$ and let $X_S$ denote the component of the
3-manifold which contains $\bdy X$.  Then $\bdy X_S$ consists of
one ($S$ separates $X$) or two ($S$ does not separate $X$) copies
of $S$ along with the torus, $\bdy X$. Note that if $S$ does not
separate $X$, then for every $\alpha$ we have that either
$X(\alpha)$ is reducible or $X(\alpha)$ contains an embedded,
closed, incompressible, two-sided surface; however, we do not need
to make a distinction as to $S$ separating or not separating $X$.

Observe that if $S$ compresses in $X(\alpha)$ for some $\alpha$,
there is a punctured disk $P$ embedded in $X_S$ with $bdry(P)$ in
a copy of $S$ in $\bdy X_S$ and punctures in the torus $\bdy X$.
In this case, there is no loss in generality to assume that $P$ is
essential. Hence, we have that the punctures form a non-empty,
pairwise disjoint collection of simple closed curves in $\bdy X$,
each having slope $\alpha$. In particular, in situation ii above,
the existence of such an annulus gives that $S$ compresses (the
annulus must also meet $\partial X$ in a non-trivial curve) in
$X(\alpha)$ where $\alpha$ is the slope of the boundary curve of
the annulus in $\bdy X$. Also, we observe in situation ii that
there is a unique slope on $\bdy X$ for an annulus which joins $S$
to $\bdy X$; for otherwise, the characteristic Seifert-Pair
Theorem \cite{j-sh,johann} gives a contradiction to our assumption
that $S$ is essential (not equivalent to a surface parallel to
$\bdy X$). Finally, if $S$ compresses in $X(\alpha)$ and $\alpha$
is not a boundary slope in $X$, then $S$ completely compresses in
$X(\alpha)$.

Now, in situation i, where there is no annulus in $X_S$ having one
boundary component a non-trivial curve in $S$ and the other in
$\bdy X$, we shall show that there is a finite and computable set
of slopes $\alpha$ for which $S$ compresses in $X(\alpha)$. (In
\cite{wu} this set is shown to have no more than 3 slopes.)

Let $\T_S$ be a triangulation of $X_S$ having precisely one vertex
in the component $\partial X$ of $\bdy X_S$ (\cite{j-r3}).

If $P$ is an essential punctured disk as above, then we may
assume that $P$ is normal in $(X_S,\T_S)$ and $P$ is least weight
in its equivalence class. We have,

$$P = \sum_i k_i F_i + \sum_{i'} l_{i'} K_{i'} + \sum_j m_j A_j +
\sum_{j'} n_{j'} A'_{j'}$$

\noindent where all of the summands are essential, normal,
fundamental surfaces in $(X_S,\T_S)$ (\cite{j-o}), and notation
has been chosen so that $\chi(F_i) < 0$, each $K_{i'}$ is either a
torus or Klein bottle, each $A_j$ is an annulus or M\"{o}bius band
with its boundary in $S$, and each $A'_{j'}$ is an annulus or
M\"{o}bius band with its boundary in $\bdy X$. Of course, it is
possible that there are no factors $K_{i'}, A_{j}$ and $A_{j'}'$.
We have written the most general sum in this situation and in
fact, each $l_{i'}$,$m_{j}$ and $n_{j'}$ might be zero.

Suppose some $n_{j'} \neq 0$. Then we may write $P= F + A'_{j'}$
for some normal surface $F$ in $(X_S,\T_S)$. Hence, either
$\partial A_{j'}' \cap \partial F = \emptyset$ and $\partial
A_{j'}'$ has slope $\alpha$ or all intersections between $A'_{j'}$
and $F$ run from $\bdy X$ to $\bdy X$ (we can assume there are no
trivial curves of intersection). So, by Proposition
\ref{same-or-comp} and in this latter case, $\bdy A'_{j'}$ has
slope $\alpha$ and $\alpha$ is a boundary slope for $X$. In fact,
there is an essential normal annulus or M\"{o}bius band in
$(X,\T)$ (possibly $A'_{j'}$ itself) with boundary slope $\alpha$
and whose projective class is a vertex of $\P(X_S,\T_S)$. So, if
$n_{j'} \neq 0$ we arrive at the conclusion that $\alpha$ is a
computable boundary slope of $X$ and we can check if $S$
compresses in $X(\alpha)$. Hence, we check if $S$ is
incompressible in $X(\alpha)$ for all boundary slopes $\alpha$
corresponding to embedded, normal annuli having projective class
at a vertex of $\mathcal{P}$($X_S$,$\T_S$), a finite, computable
set.   Note that in this situation there can be at most one slope
bounding an essential annulus with boundary in $\bdy X$.
Otherwise, $X$ would have to be a twisted I-bundle over a Klein
bottle \cite{j-sh} which contains no two-sided essential surface.
If we find more than one slope realized by vertex annuli, we have
the option of first checking which of these annuli are essential.

So, we may suppose that $n_{j'} = 0, \forall j'$; for otherwise,
$P$ would have boundary slope the same as $A_{j'}'$. Let $L(\bdy
G)$ denote the length of the boundary of the normal surface, $G$,
in $(X_S,\T_S)$, we have
:

$$L(\bdy P) = \sum_i k_i L(\bdy F_i) + \sum_j m_j L(\bdy A_j).$$
Also,
$$-\chi(P) = \sum_i k_i (-\chi(F_i)).$$
Let
$$C = \max \left\{\frac{L(\bdy F_i)}{-\chi(F_i)}\right\}.$$

\medskip

Notice that $C$ is computable for $F_i$ ranging over the embedded,
normal, fundamental surfaces in $(X_S,\T_S)$ with $\chi(F_i) < 0$
and that $L(\bdy F_i) < -\chi(F_i) C$ for all such $F_i$.  Let
$\gamma$ be the length of $bdry(P)$ and let $\gamma_\alpha$ denote
the length of the slope $\alpha$.  If $P$ has $p$ punctures, then
$-\chi(P) = p-1$ and $L(\bdy P) = \gamma + p \gamma_\alpha$. Thus,
\begin{gather*}
\gamma + p \gamma_\alpha = \sum_i k_i L(\bdy F_i) + \sum_j m_j
L(\bdy A_j) \\
\leq \sum_i (-\chi(F_i)) C + \sum_j m_j L(\bdy A_j) \\
= -\chi(P) C + \sum_j m_j L(\bdy A_j).
\end{gather*}

\noindent From this and the fact that $\gamma - \sum_j m_j L(\bdy
A_j) \geq 0$, we have $$\gamma_\alpha \leq C.$$

So, in situation i and if $S$ compresses in $X(\alpha)$, it either
compresses at a boundary slope $\alpha$ corresponding to the boundary
slope of an essential normal annulus or M\"{o}bius band in $(X,\T)$
whose projective class is a vertex of $\P(X,\T)$ or it compresses in
$X(\alpha)$ where
$\gamma_\alpha$, the length of the slope $\alpha$, satisfies
$\gamma_\alpha \leq C$, where $C$ is computable from certain
fundamental solutions in $(X_S,\T_S)$. In either case, there are
at most finitely many slopes $\alpha$ for which $S$ compresses and
we can determine precisely those slopes $\alpha$
where $S$ compresses in $X(\alpha)$.

Now, we consider situation ii, where there is an annulus in $X$
having one boundary component a non-trivial
curve in $S$ and the other a
curve in $\bdy X $.  Note in this case,
since $S$ is assumed to be essential, it follows that the boundary curve
of the annulus in $\bdy X$ is non-trivial.

First we observe that if there is such an annulus $A$, then there
is a normal one in $(X_{S},\T_{S})$. We claim that if $A$ is least
weight (in its equivalence class) among all such annuli, then $A$
is fundamental in $(X_S,\T_S)$.  To see this suppose $A$ is not
fundamental; then $A = A' + A''$ where we may write such a sum
with both $A'$ and $A''$ incompressible,
$\partial$-incompressible, and not parallel into $\partial X_S$
and $A' \cap A''$ has the smallest number of components under
these conditions \cite{j-o} . But $\chi(A') = \chi(A'') = \chi(A)
= 0$.  It follows that (with choice of notation) the possibilities
are:  both $A'$ and $A''$ are annuli each having one boundary
component in $S$ and one in $\bdy X$ (but this contradicts the
choice of $A$ being least weight);  $A'$ is an annulus having one
boundary component a non-trivial curve in $S$ and the other
boundary a curve in $\bdy X$ and $A''$ is a M\"obius band, a
torus, or a Klein bottle (but again this contradicts the choice of
$A$ being least weight); or both $A'$ and $A''$ are M\"obius
bands, one having its boundary in $S$ and the other having its
boundary in $\bdy X$.  In this last possibility there is no loss
in generality to assume that $A' \cap A''$ has exactly one
component and it is the non-separating (orientation-reversing)
simple closed curve in each.  Thus, a regular exchange along the
intersection gives the normal annulus $A$; but then an irregular
exchange gives an annulus $B$ having the same boundary as $A$ but
containing a fold (see Figure \ref{f-fold}).  So $wt(B) < wt(A)$.
But this also contradicts our choice of $A$. So, as claimed, a
least weight normal annulus in $(X_S,\T_S)$ running from $S$ to
$\bdy X$ is fundamental.

\begin{figure}[h]
 {\epsfxsize = 2 in \centerline{\epsfbox{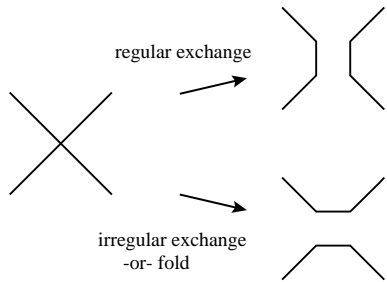}}
} \caption{Regular exchange vs. fold.}  \label{f-fold}
\end{figure}

Notice from this analysis, it is possible to have an annulus $A$ with
one boundary component a non-trivial curve in $S$ and the other a
curve in $\bdy X$ and have $A = M_1 + M_2$ where $M_i$ is a
M\"obius band; but $A$ cannot be least weight, the irregular
switch at $M_1 \cap M_2$ gives a similar annulus with lower
weight. Furthermore, there is a unique (up to isotopy) slope in
$\partial X$ for such an annulus $A$. It follows that we can find
such an $A$ and the slope $\alpha$, where $\alpha$ is the slope of
$\partial A$ on $\bdy X$. This completes our claim.

Now, as we noted above, $S$ compresses in $X(\alpha)$.  If $S$
compresses in $X(\beta)$, where $\beta \neq \alpha$, then, as in
situation i, there is a planar surface $P$ embedded in $X$ having
$bdry(P)$ in $S$ and punctures in $\bdy X$. There is no loss to
assume that $P$ is essential in $X_S$.  It follows from
\cite{cgls} that $\Delta(\alpha,\beta) = 1$ and that $S$
compresses in $X(\beta)$ for all $\beta$ where
$\Delta(\alpha,\beta) = 1$, i.e. the ``line" $L_\alpha$.

So, by considering the fundamental surfaces in $(X_{S},\T_{S})$ we
can determine if there is an annulus embedded in $X_{S}$ having
one boundary a non-trivial curve in $S$ and the other in $\partial
X$. If there is, we can find its boundary slope, say $\alpha$ in
$\partial X$. The surface $S$ compresses in $X(\alpha)$ and either
$S$ does not compress for any slope distinct from $\alpha$ or $S$
compresses in $X(\beta)$ precisely for all those slopes $\beta \in
\{\alpha\} \cup L_\alpha$. We can use the algorithm given in
\cite{haken}, see Proposition \ref{T-decideincomp}, to determine
which is the case; namely, check if $S$ compresses in $X(\beta)$
for some $\beta \in L_{\alpha}$ $(\beta \neq \alpha)$. This
completes the proof of the lemma and provides an algorithm to
decide precisely those slopes $\alpha$ for which $S$ compresses in
$X(\alpha)$.
\end{proof}

\begin{thm}
\label{T-surfacescompress} Given a knot-manifold $X$ that
contains an embedded, essential, closed, two-sided surface distinct
from $S^{2}$, then there is an algorithm to determine precisely those
slopes $\alpha$ for which $X(\alpha)$ does not contain an embedded,
incompressible, closed, two-sided surface; in particular, the set of
slopes $\alpha$ for which $X(\alpha)$ does not contain an embedded,
incompressible, closed, two-sided surface is either a finite set of
slopes or all but possibly finitely many slopes in the set
$\{\alpha_{0}\} \cup L_{\alpha_{0}}$ for some slope $\alpha_{0}$.
\end{thm}

\begin{proof}Suppose the knot-manifold $X$ is given by a
triangulation $\T$ that restricts to a one-vertex
triangulation on $\bdy X$ (recall there is an algorithm to
modify any triangulation of $X$ to such a triangulation).

If the Dehn filling $X(\alpha)$ contains an embedded,
incompressible, two-sided surface, then by Theorem \ref{T-essinx}
there is an embedded, essential, normal surface $S$ in
($X$,$\mathcal{T}$) such that the projective class of $S$ is a
vertex-solution of $\mathcal P$($X$,$\mathcal T$), the boundary
slope of $S$ is $\alpha$ (if $\partial{S} \neq \emptyset$), and
$S(\alpha)$ is an embedded, essential, normal surface in
($X(\alpha)$,$\mathcal{T}(\alpha)$). These surfaces are
constructable. The closed vertex solutions provide candidate
surfaces to which we can apply the algorithms of Lemma
\ref{closedess}. Of course, it is also possible that Dehn fillings
along boundary slopes $\alpha$ of $X$ may create embedded,
incompressible, two-sided surfaces in $X(\alpha)$; hence, those
vertex solutions that are bounded will also need to be taken into
consideration.

Let $\{S_{1},\ldots,S_{J}\}$ denote, the embedded, essential,
two-sided, connected, closed normal surfaces in $(X,\T)$ that are
not 2--spheres and whose projective class is a vertex of
$\P(X,\T)$; and let $\{F_{1},\ldots,F_{K}\}$ denote the embedded,
essential, two-sided, connected, bounded, normal surfaces in
$(X,\T)$ that are not planar and whose projective class is a
vertex of $\P(X,\T)$. By hypothesis and \cite{j-o}, the set
$\{S_{1},\ldots,S_{J}\} \neq \emptyset$.

For each surface $S_{j}, 1\leq j\leq J$, use the algorithm of
Lemma \ref{closedess} to determine those slopes $\alpha$ for which
$S_{j}$ compresses in $X(\alpha)$. For each $S_{j}$ we have that
$S_{j}$ compresses in at most a finite number of computable slopes
or for the set of slopes $\{\alpha_{j}\} \cup L_{\alpha_{j}}$
where $\alpha_{j}$ (and hence, $L_{\alpha_{j}}$) is computable.
Hence, we conclude that there is an embedded, incompressible,
two-sided surface in $X(\alpha)$ for all but a finite number of
computable slopes $\alpha$ (this includes the possibility
$(\{\alpha_{j}\} \cup L_{\alpha_{j}}) \cap (\{\alpha_{j'}\} \cup
L_{\alpha_{j'}})$,  where  $j\neq j'$) or for all but those slopes
in the set  $\{\alpha_{j}\} \cup L_{\alpha_{j}}$. If any of the
slopes where one of the closed surfaces in
$\{S_{1},\ldots,S_{J}\}$ does not remain incompressible in
$X(\alpha)$ is a boundary slope of some $F_{k}, 1\leq k\leq K$,
say the slope $\beta_{k}$ which is the boundary slope of $F_{k}$,
then we check, \cite{haken} (Theorem \ref{T-decideincomp} above)
to determine if  $F_{k}(\beta_{k})$ is incompressible in
$X(\beta_{k})$. This can only add slopes where $X(\alpha)$
contains an embedded, incompressible, two-sided surface.
\end{proof}

We now have the main theorem of this section.

\begin{thm}
 Given a knot-manifold $X$  there is
an algorithm to determine precisely those slopes $\alpha$ for
which $X(\alpha)$ is a Haken-manifold; in particular, if $X$ is
irreducible and does not contain an embedded, essential,
two-sided, closed surface, the set of slopes $\alpha$ for which
$X(\alpha)$ is a Haken-manifold is finite; if $X$ is irreducible
and does contain an embedded, essential, two-sided, closed
surface, the set of slopes $\alpha$ for which $X(\alpha)$ is not a
Haken-manifold is either a finite set of slopes or all but
possibly a finite number of slopes on the line $\{\alpha_{0}\}
\cup L_{\alpha_{0}}$ for some slope $\alpha_{0}$.
\end{thm}

\begin{proof}
From the preceding subsection, we have an algorithm to determine
precisely those slopes $\beta$ for which the Dehn filling
$X(\beta)$ is reducible and from the above Theorems, we have
algorithms to determine precisely those slopes $\gamma$ for which
the Dehn filling $X(\gamma)$ contains an embedded, incompressible,
two-sided surface. The combination of these algorithms will give
us precisely those slopes $\alpha$ for which $X(\alpha)$ is a
Haken-manifold.
\end{proof}

We next outline the steps for an algorithm to decide for a given
knot-manifold $X$ precisely those slopes $\alpha$ for which the
Dehn filling $X(\alpha)$ is a Haken-manifold, Algorithm H. We have
organized the algorithms so that one can determine precisely those
slopes $\alpha$ for which a Dehn filling $X(\alpha)$ contains an
embedded, incompressible, two-sided surface, Algorithm S; and then
we can apply our earlier algorithm to eliminate those slopes where
the manifold is reducible, Algorithm R. A fundamental step in
Algorithm S is to decide for any given embedded, two-sided, closed
surface $S$ in $X$, precisely those slopes $\alpha$ for which the
surface $S$ is incompressible in the Dehn filling $X(\alpha)$; we
give this as an independent algorithm, Algorithm I.

\vspace{.25 in}

\noindent {\bf Algorithm I.} {\it Suppose $X$ is a knot-manifold
with a triangulation $\T$ which restricts to a one-vertex
triangulation on $\partial X$. Given an embedded, two-sided,
closed, normal surface in $(X$,$\T)$, determine precisely those
slopes $\alpha$ for which the surface compresses in the Dehn
filling $X(\alpha)$.}

\vspace{.2 in}

{\bf Step 1.} Let $S$ be the given embedded, two-sided, closed,
normal surface in $(X,\T)$. Split $X$ at $S$ and let $X_{S}$
denote the component containing $\partial X$ and let $X_{S}'$
denote the other component, in the case $S$ separates $X$. The
manifold $X_{S}$ has either one or two copies of $S$ in $\partial
X_{S}$; and, if $S$ separates $X$, the manifold $X_{S}'$ has a
single copy of $S$ in $\partial X_{S}'$. Endow $X_{S}$ and
$X_{S}'$ with triangulations $\mathcal{T}_{S}$ and
$\mathcal{T}_{S}'$, respectively, so that $\T_{S}$ restricts to
the triangulation $\T$ on the boundary component $\partial X$ of
$X_{S}$ \cite{j-r3}.

{\bf Step 2.} Compute the fundamental solutions of $(X_{S},\T_{S})$
and (in the case $S$ separates $X$) of $(X_{S}',\T_{S}')$. We look
for the existence of disks and annuli among these
fundamental solutions.

{\bf Step 3.} If a fundamental solution is an embedded disk with
boundary a non-trivial curve in a copy of $S$, then the surface $S$
compresses in $X$ and therefore will compress in $X(\alpha)$ for every
slope $\alpha$. If this is not the case but a fundamental solution is
an embedded disk with boundary a non-trivial curve in $\partial X$,
then the knot-manifold $X$ is reducible and $S$ is incompressible in
$X(\alpha)$ for every slope $\alpha$. (Notice that if $\partial X$
compresses and $X$ is irreducible then $X$ is a solid torus, $S$ would
necessarily compress and we would have a fundamental solution that is
an embedded disk with boundary a non-trivial curve in a copy of $S$,
i.~e.~, we would have found such a disk in the first part of this
step.)  If either type of disk is found, then the algorithm is
complete and we have either $S$ compresses for every slope or $S$
compresses for no slope.

{\bf Step 4.} We have that no fundamental solution found in Step 2
is a disk with non-trivial boundary in either a copy of $S$ or in
$\partial X$.  Now, look for fundamental solutions that are
embedded annuli having one boundary a non-trivial curve in a copy
of $S$ and the other in $\partial S$. If there are two such annuli
having distinct slopes in $\partial X$, then $S$ is equivalent to
a peripheral torus and compresses in $X(\alpha)$ for every slope
$\alpha$. If there is only one slope for all such annuli, then go
to Step 6. If there are no such annuli, then go to Step 5.

{\bf Step 5.} We have that no fundamental solution found in Step 2
is an embedded disk or an embedded annulus having one boundary a
non-trivial curve in $S$ and the other in $\partial X$. However,
there may be fundamental solutions found in Step 2 that are
embedded annuli or M\"{o}bius bands having their boundary
non-trivial curves in $\partial X$. Let $\{A_{1},\ldots,A_{m}\}$
denote such fundamental solutions and let $\{F_{1},\ldots,
F_{n}\}$ be the set of all embedded fundamental solutions of
$(X_{S},\T_{S})$ with $\chi(F_{i}) < 0, \forall i.$ Set $$C = \max
\left\{\frac{L(\partial F_{i})}{-\chi(F_{i})}\right\}.$$ \noindent
Let $\{\alpha_{1},\ldots ,\alpha_{K}\}$ denote the slopes in
$\partial X$ that either have length $\lambda_{\alpha_{i}} \leq C$
or are a boundary slope for some $A_{j}, 1\leq j\leq m$. Recall,
in this situation, the surface $S$ will compress in $X(\alpha)$ if
and only if there is an $i, 1\leq i\leq K$, $\alpha = \alpha_{i}$
and $S$ compresses in $X(\alpha_{i})$.  (Also recall that there is
at most one slope bounding an essential annulus $A_j$; it may be
advantageous to check whether each vertex annulus is inessential
before listing the slope.)

\noindent For $\alpha_{i} \in \{\alpha_{1},\ldots ,\alpha_{K}\}$,
build $X(\alpha_{i})$ via a layered triangulation and check if $S$
compresses in $X(\alpha_{i})$. Let $\{\alpha_{i_{1}}, \ldots,
\alpha_{i_{J}}\}$ be the set of slopes in $\{\alpha_{1},\ldots
,\alpha_{K}\}$ for which $S$ compresses in $X(\alpha_{i})$. The
algorithm terminates having found this finite set of slopes as
precisely the set of slopes $\alpha$ for which the surface $S$
compresses in $X(\alpha)$.

{\bf Step 6.} Let $A$ denote an embedded annulus in
$(X_{S},\T_{S})$ having one boundary a non-trivial curve in a copy
of $S$ and the other in $\partial S$ (If there is such an $A$,
then one may be constructed.) The component of $\partial A$ in
$\partial S$ is a non-trivial curve, say, with slope $\alpha_{0}$.
Choose any slope $\beta \in L_{\alpha_{0}}$. Determine if $S$
compresses in $X(\beta)$. If $S$ does not compress in $X(\beta)$,
then $S$ is incompressible in all Dehn fillings $X(\alpha), \alpha
\neq \alpha_{0}$, and so, $S$ compresses in $X(\alpha)$ for
precisely one slope, the slope $\alpha_{0}$. If $S$ compresses in
$X(\beta)$, then $S$ compresses in all Dehn fillings $X(\alpha)$
where $\alpha \in \{\alpha_{0}\} \cup L_{\alpha_{0}}$.  This
completes Algorithm I.

\vspace{.25 in}

We now consider an algorithm to decide for a given knot-manifold
$X$ the set of slopes $\alpha$ for which the Dehn filling
$X(\alpha)$ contains an embedded, incompressible, two-sided
surface.

\vspace{.25 in}

\noindent {\bf Algorithm S.} {\it Given a knot-manifold $X$,
determine precisely those slopes $\alpha$ for which the Dehn
filling $X(\alpha)$ contains an embedded, incompressible,
two-sided surface.}

\vspace{.2 in}

{\bf Step 1.} We have the knot-manifold $X$ given via a
triangulation. Endow $X$ with a triangulation $\mathcal{T}$ that
has precisely one vertex in $\partial X$. (An algorithm is given
in \cite{j-r3} to modify a given triangulation of $X$ to such a
triangulation.)

{\bf Step 2.} Compute the vertices of
$\mathcal{P}(X,\mathcal{T})$.

{\bf Step 3.} Make two lists: $\mathcal{G} =
\{G_{1},\ldots,G_{J}\}$, those vertices of
$\mathcal{P}(X,\mathcal{T})$ that are projective classes of
embedded, closed, normal surfaces that are not 2--spheres in
$(X,\T)$; and $\mathcal{B} = \{B_{1},\ldots,B_{K}\}$, those
vertices $\mathcal{P}(X,\mathcal{T})$ which are projective classes
of embedded, non-planar, normal surfaces in $(X,\mathcal{T})$ and
have nonempty boundary consisting entirely of non-trivial simple
closed curves in $\partial X$. For each surface in $\mathcal{B}$
compute its boundary slope. Denote the boundary slope of $B_{k},
1\leq k\leq K$, by $\beta_{k}$; these are a subset of the boundary
slopes of $X$.

{\bf Step 4.} If $\mathcal{G} = \emptyset$, it follows that $\partial
X$ compresses and $X$ is a solid torus or a non-trivial connected sum
of a solid torus and a 3--manifold $M'$. Furthermore, $M'$ does not
contain any embedded, incompressible, two-sided, closed surfaces. So,
in either case, $X(\alpha)$ does not contain an embedded,
incompressible, two-sided closed surface for any slope $\alpha$ and
the algorithm terminates.

{\bf Step 5.} We have $\mathcal{G} = \{G_{1},\ldots,G_{J}\} \neq
\emptyset$. For each $G_{i} \in \mathcal{G}, 1\leq i\leq J$, apply
Algorithm S to decide precisely those slopes $\alpha$ for which
the surface $G_{i}$ compresses in the Dehn filling $X(\alpha)$.
For the surface $G_{i}$, let this set of slopes be denoted
$\mathcal{A}_{i}$. Recall the possibilities are: a finite set of
slopes, or all the slopes on a ``line" $\{\alpha_{i_{0}}\} \cup
L_{\alpha_{i_{0}}}$ for some slope $\alpha_{i_{0}}$, or every
slope (i.~e.~, $G_{i}$ either compresses in $X$ or is peripheral).
Let $$\mathcal{A} = \bigcap_{i=1}^{J}\mathcal{A}_{i}.$$

{\bf Step 6.} For each slope $\beta_{k}$ found in Step 3,
construct $X(\beta_{k})$ via a layered triangulation. It can be
determined if the surface $B_{k}(\beta_{k})$, $B_{k}$ also found
in Step 3, compresses in $X(\beta_{k})$. Let
$\{\beta_{k_{1}},\ldots,\beta_{k_{n}}\}$ be the set of slopes for
which $B_{k_{j}}, 1\leq j\leq n$, does NOT compress in
$X(\beta_{k_{j}})$

{\bf Step 7.} If $\mathcal{A}$ is finite (i.~e.~, there are only
finitely many slopes $\alpha$ for which \underline{all} the
surfaces in $\mathcal{G}$ compress in $X(\alpha)$), then the set
of slopes $\mathcal{A} \setminus
\{\beta_{k_{1}},\ldots,\beta_{k_{n}}\}$ is precisely the set of
slopes $\alpha$ for which $X(\alpha)$ does NOT contain an
embedded, incompressible, two-sided surface.  If $\mathcal{A}$ is
infinite, there are two possibilities: $\mathcal{A} =
\{\alpha_{i_{0}}\} \cup L_{\alpha_{i_{0}}}$ for some slope
$\alpha_{i_{0}}$ or $\mathcal{A}$ is the set of all slopes. In the
first case, the set of slopes $\mathcal{A} \setminus
\{\beta_{k_{1}},\ldots,\beta_{k_{n}}\}$ is precisely the set of
slopes $\alpha$ for which $X(\alpha)$ does NOT contain an
embedded, incompressible, two-sided surface. In the second case,
$\{\beta_{k_{1}},\ldots,\beta_{k_{n}}\}$ is precisely the set of
slopes $\alpha$ for which $X(\alpha)$ does contain an embedded,
incompressible, two-sided surface. This terminates Algorithm S.

\vspace{.25 in}

Finally, we are prepared to give an algorithm to determine the
manifolds in the space of Dehn fillings which are Haken-manifolds.

\vspace{.25 in}

\noindent {\bf Algorithm H.} {\it Given a knot-manifold $X$,
determine precisely those slopes $\alpha$ for which the Dehn
filling $X(\alpha)$ is a Haken-manifold.}

\vspace{.2 in}

{\bf Step 1.} We have the knot-manifold $X$ given via a
triangulation. Endow $X$ with a triangulation $\mathcal{T}$ that has
precisely one vertex in $\partial X$.

{\bf Step 2.} Employ Algorithm S to determine precisely those
slopes $\alpha$ for which the Dehn filling $X(\alpha)$ contains an
embedded, incompressible, two-sided surface.

{\bf Step 3.} Employ Algorithm R to determine precisely those
slopes $\alpha$ for which the Dehn filling $X(\alpha)$ is
irreducible (those for which it is not reducible).

The slopes common to those found in Steps 2 and 3 are precisely
the set of slopes $\alpha$ for which $X(\alpha)$ is a
Haken-manifold.

\subsection{Fibered manifolds in Dehn Surgery Space.} In this
section we provide an algorithm to determine for a given
knot-manifold $X$ precisely those slopes $\alpha$ for which the
Dehn filling $X(\alpha)$ fibers as a surface bundle over a circle.
We wish to thank Robert Myers who suggested that our methods
should solve this problem. Our proof is based on material from
lectures of the first author given at University of Melbourne a
decade ago. Revision of this work appears in \cite{jac}; we state
the results we need below without proof.

\begin{thm}\cite{jac,j-t}
\label{surfacefiber} Suppose $\T$ is a triangulation of the 3--manifold
$M$. Given a properly embedded normal surface $F$ in $(M,\T)$, there is an
algorithm to determine if $F$ is a fiber in a fibration of $M$ over $S^1$.
\end{thm}

\begin{thm}\cite{jac}
\label{fiber}Given a 3--manifold $M$, there is an algorithm to determine
if $M$ is a
fibration over $S^1$.
\end{thm}

\begin{thm} Given a knot-manifold $X$ there is an algorithm to determine
precisely those slopes $\alpha$ for which the Dehn filling
$X(\alpha)$ is
 a fibration over $S^1$.
\end{thm}

\begin{proof} We assume $X$ is given via a triangulation $\T$ that
restricts to a one-vertex triangulation on $\bdy X$. We separate
the argument into two cases depending  on $X$ reducible or $X$
irreducible.

If $X$ is reducible, then by Theorem \ref{T-irreducibleslopes}, we
have:

\renewcommand{\labelitemi}{--}

\begin{itemize}
\item $X$ contains an embedded, non-separating (hence, essential) 2--sphere
and $X(\alpha)$ is reducible for all $\alpha$,
\item $X$ contains two separating, embedded, disjoint, inequivalent,
essential 2--spheres and $X(\alpha)$ is reducible for all $\alpha$,
\item $X$ is a connected sum of a nontrivial knot-manifold in $S^3$ and an
irreducible 3--manifold and there is precisely one computable
slope for which a Dehn filling is irreducible, or
\item $X$ is a connected sum of a solid torus and an irreducible
3--manifold and there is a computable line of slopes for which $X(\alpha)$
is irreducible.
\end{itemize}

If $X$ contains an embedded, non-separating 2--sphere $S$, it may
be possible that $X(\alpha)$ fibers over $S^1$ with fiber the
surface $S$. There is an algorithm to determine this; again, we
call upon Theorem 6.4 of the next section. Hence, if $X$ contains
an embedded, non-separating 2--sphere, by Theorem
\ref{T-recognize-2}, there is an algorithm to find one, say $S$ is
such a 2--sphere. Split the knot-manifold $X$ at $S$ to form the
3--manifold $X_{S}$. The manifold $X_{S}$ has two copies of $S$ in
its boundary, along with $\bdy X$. We can fill the two copies of
$S$ with 3--cells to get a new knot-manifold $\widehat{X_S}$. The
slopes $\alpha$ for which $X'(\alpha)$ is $S^3$ are precisely the
slopes for which $X(\alpha)$ fibers over $S^1$ with fiber the
2--sphere $S$. Hence, we have $X(\alpha)$ does not fiber for all
$\alpha$, or $X(\alpha)$ fibers for a unique, and computable,
slope $\alpha$, or $X(\alpha)$ fibers for a computable line of
slopes $L_{\alpha}$.

If $X$ contains two separating, embedded, disjoint, independent,
essential 2--spheres, it is not possible for $X(\alpha)$ to fiber for any
$\alpha$.

If $X$ is a connected sum of a nontrivial knot-manifold in $S^3$
and an irreducible 3--manifold, we have that there is precisely
one computable slope for which a Dehn filling is irreducible; say
for $\alpha_{0}$, we have $X(\alpha_{0})$ irreducible.  If we
denote the irreducible 3--manifold by $N$, then $X(\alpha)$ will
fiber over $S^1$ with fiber a surface only for the slope
$\alpha_{0}$ and then if and only if we have $N$ fibers over $S^1$
with fiber a surface. By \cite{jac}, Theorem \ref{fiber} above,
there is an algorithm to determine if $N$ fibers over $S^1$ with
fiber a surface.

If $X$ is a connected sum of a solid torus, say $M$, and an
irreducible 3--manifold, say $N$, there is a computable line of
slopes $L_{\mu}$, where $\mu$ is the (computable) slope of the
meridian of the solid torus $M$, for which Dehn fillings on $X$
are irreducible. Hence, $X(\alpha)$ can fiber over $S^1$ with
fiber a surface only for those slopes $\alpha \in  L_{\mu}$ and
then if and only if we have $N$ fibers over $S^1$ with fiber a
surface. Again, this can be determined by Theorem \ref{fiber}
above.

Hence, if $X$ is reducible, we can determine precisely those
slopes $\alpha$ for which the Dehn filling $X(\alpha)$ fibers over
$S^1$.

So, suppose $X$ is irreducible. The argument in this case is very
similar to the combination of arguments used in Theorem
\ref{T-surfacesincompress}, Lemma \ref{closedess}, and Theorem
\ref{T-surfacescompress}. By Theorem \ref{T-essinx}, if
$X(\alpha)$ fibers over $S^1$, there is a vertex solution $F$ of
$\mathcal{P}(X,\mathcal T)$ that is an embedded, essential,
two-sided surface and $F(\alpha)$ is a fiber in a fibration over
$S^1$. We may also assume that $F$ does not separate $X$ for
otherwise $F(\alpha)$ could not be a fiber in a fibration of
$X(\alpha)$ over $S^1$.

Suppose $\bdy F\neq\emptyset$. By Theorem \ref{surfacefiber} there
is an algorithm to determine if $F(\alpha)$ is a fiber in a
fibration of $X(\alpha)$ over $S^1$. There are only finitely many
such surfaces we need to check.

Suppose $\bdy F =\emptyset$. We have the embedded, essential,
closed normal surface $F$ in $X$ and we wish to determine
precisely those slopes $\alpha$ for which $F=F(\alpha)$ is a fiber
in a fibration of $X(\alpha)$ over $S^1$.

As in the proof of Lemma \ref{closedess}, we consider two
possibilities:
\renewcommand{\theenumi}{\roman{enumi}}
\begin{enumerate}
\item there is no annulus in $X$ having one boundary component a
non-trivial curve in $F$ and the other a curve in $\bdy X$, or
\item there is an annulus in $X$ having one boundary component a
non-trivial curve in
$F$ and the other a curve in $\bdy X$.
\end{enumerate}

Suppose we are in situation i where there is no annulus in $X$
having one boundary component a non-trivial curve in $F$ and the
other a curve in $\bdy X$. Split $X$ at $F$ to get the 3--manifold
$X_{F}$. The manifold $X_{F}$ has two copies ($F$ does not
separate $X$) of $F$, say $F_0$ and $F_1$,  along with $\bdy X$ as
its boundary.

We extend our notion of a Dehn filling to this situation where the
manifold $X_{F}$ has components of the boundary other than the
torus $\bdy X$. A {\it slope} $\alpha$ will be an isotopy class of
a simple closed curve in $\bdy X$ and a Dehn filling of $X_F$
along $\alpha$ is the 3--manifold obtained by attaching a solid
torus $V_{\alpha}$ to $X_F$ via a homeomorphism of $\bdy X$ to
$\bdy V_{\alpha}$ taking the slope $\alpha$ to a meridian of
$V_{\alpha}$. We denote the Dehn filling of $X_F$ along $\alpha$
by $X_{F}(\alpha)$. With this notation, we have $X(\alpha)$ will
fiber over $S^1$ with fiber $F=F(\alpha)$ if and only if the Dehn
filling $X_{F}(\alpha)$ is homeomorphic to the product $F\times
[0,1]$. Hence, we wish to determine precisely those slopes
$\alpha$ for which the Dehn filling $X_{F}(\alpha)$ is a product.

Give $X_{F}$ a triangulation $\mathcal{T}_F$ that restricts to a
one-vertex triangulation on the component of $\bdy X_F$
corresponding to $\bdy X$. If $X_{F}(\alpha)$ is a product, then
either there is an embedded, essential, punctured annulus $Q$ in
$X_F$ with one component of $bdry(Q)$ a nontrivial curve in
$F_{0}$ and the other a nontrivial curve in $F_{1}$ and punctures
in $\bdy X$ with slope $\alpha$ or we have $V_{\alpha}$, the
attached solid torus, is contained in a 3--cell in $X(\alpha)$.
However, the latter situation could only happen if $X$ were
reducible. We conclude that for any Dehn filling with $X(\alpha)$
a product, there is such a punctured annulus $Q$ having punctures
in $\bdy X$ with slope $\alpha$.

We now use an average length estimate similar to that in the proof
of Lemma \ref{closedess} to give an algorithm to find such a
punctured annulus.

If $Q$ is an essential punctured annulus as above, then we may
assume that $Q$ is normal in $(X_F,\T_F)$ and $Q$ is least weight
in its equivalence class. We have,

$$Q = \sum_i k_i G_i + \sum_{i'} l_{i'} K_{i'}  + \sum_j p_j A^{0}_j +
\sum_{j'} q_{j'} A^{1}_{j'} +
\sum_{j''} r_{j''} A^{0,1}_{j''} + \sum_{k} s_{k} A^{\bdy}_{k}$$

\noindent where all of the summands are essential, normal,
fundamental surfaces in $(X_S,\T_S)$ (\cite{j-o}), and notation
has been chosen so that $\chi(G_i) < 0$, each $K_{i'}$ is either a
torus or Klein bottle, $A^{0}_j$ and $A^{1}_{j'}$ are annuli or
M\"{o}bius bands with their boundaries in $F_{0}$ or $F_{1}$,
respectively, each  $A^{0,1}_{j''}$ is an annulus with one
boundary component in $F_0$ and the other in $F_{1}$, and each
$A^{\bdy}_{k}$ is an annulus or M\"obius band with its boundary in
the copy of $\bdy X$. Of course, it is possible that there are no
factors $K_{i'}, A^{0}_j, A^{1}_{j'}, A^{0,1}_{j''}$ and
$A^{\bdy}_{k}$.  We have written the most general sum in this
situation. Also, by assumption for this case, there are no annuli
in $X$ (and hence in $X_F$) having one boundary a nontrivial curve
in $F$ and the other boundary in $\bdy X$.

As in the proof of Lemma \ref{closedess}, if there are any annuli
or M\"obius bands of type $A^{\bdy}_{k}$, then the slope $\alpha$
is the same as the boundary slope of $A^{\bdy}_{k}$, which is a
computable slope of a fundamental surface of $(X_{F},\T_{F})$.
(Actually, if there is an embedded, essential annulus
$A^{\bdy}_{k}$, then there is one whose projective class is also a
vertex solution of $\mathcal{P}(X,\T)$.) So, we may assume that
each $s_{k} = 0$.

Again we let $L(\bdy G)$
denote the length of the boundary of a normal surface, $G$, in
$(X_F,\T_F)$, we have :

$$L(\bdy Q) = \sum_i k_i L(\bdy G_i) + \sum_j p_j L(\bdy A^{0}_j) +
\sum_{j'} q_{j'} L(\bdy A^{1}_{j'}) +
\sum_{j''} r_{j''} L(\bdy A^{0,1}_{j''}).$$
Also,
$$-\chi(Q) = \sum_i k_i (-\chi(G_i)).$$
Let
$$C = \max \left\{\frac{L(\bdy G_i)}{-\chi(G_i)}\right\}.$$

Notice that $C$ is computable for $G_i$ ranging over the embedded,
normal, fundamental surfaces in $(X_S,\T_S)$ with $\chi(G_i) < 0$
and $L(\bdy G_i) < -\chi(G_i) C$ for all such $G_i$.  Let $\gamma_{0}$ and
$\gamma_{1}$ be the length of the components of $bdry(Q)$ in $F_{0}$ and
$F_{1}$, respectively; and let
$\gamma_\alpha$ denote the length of the slope $\alpha$.  If $Q$ has $q$
punctures, then
$-\chi(Q) = q$ and $L(\bdy Q) = \gamma_{0} + \gamma_{1} + q
\gamma_\alpha$.
Thus, if we set $L' =  \sum_j p_j L(\bdy A^{0}_j) +
\sum_{j'} q_{j'} L(\bdy A^{1}_{j'}) +
\sum_{j''} r_{j''} L(\bdy A^{0,1}_{j''})$, we have
\begin{gather*}
\gamma_{0} + \gamma_{1} + q \gamma_\alpha =\sum_i k_i L(\bdy G_i) + L' \\
\leq \sum_i (-\chi(G_i)) C + L' \\
= -\chi(Q) C + L'.
\end{gather*}

From this and the fact that $\gamma_{0} + \gamma_{1} - L' \geq 0$,
we have $$\gamma_\alpha \leq C.$$

So, in situation i and if $X_{F}(\alpha)$ is a product, $\alpha$
it either the boundary slope of an essential normal annulus or
M\"{o}bius band in $(X,\T)$ whose projective class is a vertex of
$\P(X,\T)$ or $\gamma_\alpha$, the length of the slope $\alpha$,
satisfies $\gamma_\alpha \leq C$, where $C$ is computable from
certain fundamental solutions in $(X_F,\T_F)$. In either case,
there are at most finitely many computable slopes $\alpha$ for
which a Dehn filling of $X_F$ can be homeomorphic to a product $F
\times [0,1]$ and the Dehn filling $X(\alpha)$ can be a fibration
over $S^1$ with fiber $F$.

Now, we consider situation ii where there is an annulus having one
boundary component a nontrivial curve in $F_{i}$ and the other in
$\bdy X$, for either $i=0, i=1$ or both. Note that if there is any
combination of such annuli, then by \cite{j-sh, johann} and the
fact that $F$ is not peripheral, there is a unique slope for the
components of all such annuli in $\bdy X$. Furthermore, by the
same argument as that in  Lemma \ref{closedess}, if there is such
an annulus, then there is one that is a fundamental solution of
$(X_F,\T_F)$; hence, the boundary slope on $\bdy X$ of such an
annulus, say $\alpha_{0}$, can be computed.

Notice that our analogy with Lemma \ref{closedess} diverges at
this point, as the existence of such an annulus gives that $F$
compresses in $X(\alpha_{0})$ and so could not be a fiber in a
fibration of $X(\alpha_{0})$ over $S^1$. However, if for some Dehn
filling along a slope $\alpha$, we do have that $X(\alpha)$ fibers
over $S^1$ with fiber $F$, then, as above, there is an embedded,
essential, punctured annulus in $X_{F}$ having one boundary
component in $F_0$ and the other in $F_1$ and punctures in $\bdy
X$ having slope $\alpha$. It follows from \cite{cgls} that $\Delta
(\alpha,\alpha_{0}) \leq 1$; hence, it is only possible for
$X(\alpha)$ to fiber over $S^1$ with fiber $F$ for $\alpha \in
L_{\alpha_0}$. Furthermore, one such Dehn filling will fiber over
$S^1$ with fiber $F$ if and only if all do. By Theorem
\ref{fiber}, there is an algorithm to check if any one does fiber
over $S^1$ with fiber $F$. This completes the proof. \end{proof}

\vspace{.25 in}

\noindent {\bf Algorithm F.} {\it Given a knot-manifold $X$,
determine precisely those slopes $\alpha$ for which the Dehn
filling $X(\alpha)$ is a fibration over $S^1$.}

\vspace{.2 in}

\noindent {\bf Step 1.} $X$ is given via a triangulation. Endow
$X$ with a triangulation $\T$ that restricts to a one-vertex
triangulation on $\bdy X$.

\noindent {\bf Step 2.} Compute $\mathcal{P}(X,\T)$.

\noindent {\bf Step 3.} Construct the irreducible decomposition of $X$.

If $X$ has a non-separating 2--sphere, say $S$, then split $X$ at
$S$ to form $X_S$ and then fill the resulting 2--sphere boundary
components with 3--cells to get the knot-manifold $\widehat{X_S}$.
Use the algorithm of Theorem 6.4 to determine those slopes
$\alpha$ for which $\widehat{X_S}(\alpha)$ is $S^3$. It is
precisely these slopes $\alpha$ for which the Dehn filling
$X(\alpha)$ fibers over $S^1$ with fiber the 2--sphere $S$; and
the algorithm terminates.

If $X$ has two separating, independent, essential 2--spheres, then
$X(\alpha)$ does not fiber over $S^1$ for any Dehn filling
$\alpha$; and the algorithm terminates.

If $X$ has precisely one, separating, essential 2--sphere, say
$S$, then split $X$ at $S$ and fill the resulting 2--sphere
boundary components with 3--cells to get the two 3--manifolds $M$
and $N$, where we choose notation so that $M$ contains the copy of
$\bdy X$. Determine if $N$ is a  fibration over $S^1$. If $N$ does
not fiber over $S^1$, then $X(\alpha)$ will not fiber over $S^1$
for any $\alpha$; and the algorithm terminates. If $N$ does fiber
over $S^1$, use Theorem 6.4 to determine those slopes $\alpha$ for
which Dehn filling on $M$ along $\alpha$ gives $S^3$; it is
precisely these slopes for which $X(\alpha)$ is a fibration over
$S^1$ and the algorithm terminates.

If $X$ is irreducible, go to the next step.

\noindent {\bf Step 4.} Let $\mathcal{F} = \{F_1,\ldots,F_J\}$
denote the collection of all embedded, closed, non-separating,
two-sided, normal surfaces whose projective class is a
vertex-solution of $\mathcal{P}(X,\T)$. Let $\mathcal{B} =
\{B_1,\ldots,B_{K}\}$ denote the collection of all embedded,
non-separating, two-sided, normal surfaces with boundary
consisting of nontrivial curves in $\bdy X$ whose projective class
is a vertex-solution of $\mathcal{P}(X,\T)$. Compute the boundary
slopes of the surfaces in $\mathcal{B}$, let $\beta_k$ denote the
boundary slope of $B_k$.

\noindent {\bf Step 5.} Check if $B_k(\beta_{k})$ is a fiber in a
fibration of $X(\beta_k)$ over $S^1$. In this way we get a
possible finite number of slopes $\beta_{k_{i}}$ for which
$X(\beta_{k_{i}})$ fibers over $S^1$.

\noindent {\bf Step 6.} For each $F_j \in\mathcal F$, split $X$ at
$F_j$ and form $X_{F_{j}}$. Triangulate $X_{F_{j}}$ with a triangulation
that restricts to a one-vertex triangulation on $\bdy X$, say
$\mathcal{T}_{F_{j}}$. Compute the fundamental solutions of
$(X_{F_{j}},\mathcal{T}_{F_{j}})$.

\noindent {\bf Step 7.} Consider those $j, 1\leq j\leq J$, for
which a fundamental solution is an embedded annulus with one
boundary component in $F_{j}$ and the other in $\bdy X$, compute
the slope of the component of the boundary in $\bdy X$, say
$\alpha_j$. Compute the line $L_{\alpha_{j}} = \{\beta :
\Delta(\alpha_{j},\beta) \leq 1\}$. For some $\beta_{0} \in
L_{\alpha_{j}}$ check if $F_j$ is a fiber in a fibration of
$X(\beta_{0})$ over $S^1$. If yes, then $F_j$ is a fiber in a
fibration over $S^1$ for all $\beta \in L_{\alpha_{j}}$. If no,
then $X(\alpha)$ does not fiber over $S^1$ with $F_{j}$ a fiber
for any $\alpha$.

\noindent {\bf Step 8.} Consider those $j, 1\leq j\leq J$, for
which no fundamental solution is an embedded annulus with one
boundary component in $F_{j}$ and the other in $\bdy X$.

If some fundamental solution is an embedded annulus with both its
boundary components in $\bdy X$, then compute the boundary slope
of such an annulus in $\bdy X$, say, $\alpha_{j_{0}}$. Note there
is only one such slope for such embedded essential annuli and it
may be necessary to determine if the annulus is essential.

Let $\{G_{j_{1}},\ldots,G_{j_{N_{j}}}\}$ denote the fundamental
solutions of $(X_{F_{j}},\mathcal{T}_{F_{j}})$ for which $\chi
(G_{j_{i}}) < 0$. Compute $$C_{j} = \left\{\frac{L(\partial
G_{j_{i}})}{-\chi(G_{j_{i}})}\right\}.$$

\vspace{.1 in} \noindent For each $j, 1\leq j\leq J$, compute all
slopes in $\bdy X$ having length less than $C_{j}$. Let $\{
\alpha_{j_0}, \alpha_{j_{1}},\ldots,\alpha_{j{K_{j}}}\}$ be this
set of slopes along with the slope $\alpha_{j_{0}}$, if  found
above. Check if $F_j$ is a fiber in a fibration of
$X(\alpha_{j_{i}})$ for each of these slopes. It is precisely
these slopes $\alpha$ for which the surface $F_{j}$ is a fiber in
a fibration over $S^1$.

\noindent {\bf Step 9.} The union of the slopes found in Step 5,
in Step 7, and in Step 8 determine all slopes $\alpha$ for which
$X(\alpha)$ fibers over $S^1$; and the algorithm terminates.


\section{Decision Problems in the Space of Dehn Fillings: Heegaard Surfaces}
\label{s-heegaard}

In the last section we used normal surface theory to determine
precisely those slopes for which Dehn fillings had ``interesting"
essential surfaces. In this section we use {\it almost normal
surface theory} in order to add Heegaard surfaces to our list of
interesting surfaces. We employ the ideas of J.~H.~Rubinstein
(almost normal surfaces and sweep outs) and of D.~Gabai (thin
position), along with the work of A.~Thompson \cite{thompson} and
M.~Stocking \cite{stocking}. We are able to give algorithms to
determine for a given knot-manifold $X$ precisely those slopes
$\alpha$ for which the Dehn filling $X(\alpha)$ is either the
3--sphere or a lens space.

We will use two important solutions to the homeomorphism problem
for 3--manifolds. The first is a restatement of Theorem
\ref{T-recognize-1}, which is directly applicable to this section.

\begin{thm}\label{r-3-sphere}
\cite{rubinstein,thompson} Given a compact 3--manifold $M$, it can
be decided if $M$ is homeomorphic to $S^{3}$.
\end{thm}

\begin{thm}\label{r-lensspace}
\cite{rubinstein}
Given a compact 3--manifold $M$ it can be decided if $M$ is
homeomorphic to a lens space.
\end{thm}


Both of these algorithms are based on the fact that given a
triangulation of $S^3$ or of a lens space, a strongly irreducible
Heegaard surface (a 2--sphere in $S^3$ or a torus in a lens space)
is isotopic to an almost normal surface.  (See \cite{stocking} for
the general case.)   From these algorithms, if we are given a
knot-manifold $X$ and a slope $\alpha$, we are able to determine
if the Dehn filling $X(\alpha)$ is $S^3$ or a lens space; however,
these algorithms are not sufficient (do not provide finite
algorithms) to answer the general questions as to precisely which
slopes $\alpha$ the Dehn filling $X(\alpha)$ is either $S^3$ or a
lens space or whether there is a Dehn filling of $X$ that is
either $S^3$ or a lens space.

The generic model \cite{cgls} is that there are only a finite
number of slopes along which a knot-manifold $X$ can be filled to
produce $S^3$ or a lens space. We obtain a finiteness result and
more by showing that the slopes giving Dehn fillings that are
either $S^{3}$ or a lens space arise as the slopes of embedded
normal or almost normal surfaces or as the slope of an edge in
$\bdy X$ of the triangulation, a so-called {\it boundary edge}.
This is a finite computable set of slopes. We identify and analyze
a few exceptional cases that arise when the core of the solid
torus that is attached to $\bdy X$ is isotopic into the minimal
genus Heegaard splitting of the Dehn filling. In this event, thin
position does not provide the desired conclusion. For example, for
fillings giving $S^3$, this occurs when the core of the attached
solid torus is isotopic into a 2--sphere; so, the core is an
unknot and its exterior is a solid torus.

Fortunately, we are able to identify and analyze these exceptions.
A knot-manifold is a solid torus if and only if it has
compressible boundary and is irreducible.  We can determine when a
manifold has compressible boundary \cite{haken, j-o} and when it
is irreducible \cite{rubinstein,thompson} (see Theorem
\ref{T-recognize-2}).  Note that Haken's original algorithm to
recognize the unknot \cite{haken} consisted of finding a
compressing disk for the boundary  of the knot-manifold combined
with the advance knowledge that the knot-manifold was contained in
$S^3$ (was irreducible).

We make the following convention. If $X$ is embedded in $M$ as a
knot-manifold, the exterior of a knot in $M$, then there is a
unique slope in $\bdy X$ that we call the meridional slope or a
meridian. In $S^3$ there is a unique slope that has distance 1
from the meridian and bounds a properly embedded, orientable
surface in $X$; its slope is called a longitude. However, in
general, there is no such unique curve in $\bdy X$; so, we shall
refer to any slope in the line of slopes having distance $1$ from
the meridian as a {\it generalized longitude}.

We have the following lemma.

\begin{lem} \label{knot-normal}
Let $X$ be the exterior of a non-trivial knot in $S^3$ and $\T$ a
one-vertex triangulation of $X$. Then $(X,\T)$ contains a normal
or almost normal planar surface with an essential boundary curve
that has slope a meridian.
\end{lem}

\begin{rem}
For the purposes of this lemma an almost normal surface possesses
a single octagon (no tubes).
\end{rem}

\begin{proof}
The proof of this lemma is adapted from Thompson's proof of the
existence of an almost normal sphere in a triangulation of $S^3$
\cite{thompson}. It differs in that we guarantee that there is a
level surface which intersects the boundary torus $\bdy X$ in a
collection of curves that includes an essential curve with
meridional slope. Both are applications of Gabai's notion of {\it
thin position} \cite{gabai} to an embedding of the 1-skeleton of a
triangulation, and we assume that the reader has a familiarity
with the basic concepts. For more detailed information on thin
position for graphs the reader is directed to \cite{s-t}.

By assumption, $X$ is the exterior of a non-trivial knot in $S^3$
and it is endowed with a one-vertex triangulation, $\T$. Note
however, that $\T$ is {\it not} a triangulation of $S^3$; the
exterior of the 2-skeleton of $\T$ is a collection of tetrahedra
and a single solid torus (the neighborhood of the knot).

Consider the singular foliation of $S^3$ induced by its genus 0
Heegaard splitting. Each leaf of the foliation, $S_t, 0 < t
< 1$, is a $2$-sphere except for $S_0$ and $S_1$ which are single
points. We think of this foliation in terms of the height function
that it induces: $h: S^3 \rightarrow [0,1]$. Arrange $\T$ to be in
general position with respect to this foliation and so that the
boundary vertex is held fixed at $S_1$. We define the width of the
one-skeleton $\T^{(1)}$ to be
\begin{gather}
w(\T^{(1)}) = \sum | \T^{(1)} \cap S_t|,
\end{gather}
where the sum is taken over level surfaces, $S_t$, where one level
surface is chosen between each pair of successive critical values
of $h: \T^{(1)} \rightarrow [0,1]$. Among such generic embeddings
of $\T$ choose one which minimizes the width of the 1-skeleton,
$w(\T^{(1)})$. This is called a {\it thin position} for
$\T^{(1)}$.

\begin{claim}
Suppose that a sphere $S$ intersects $\bdy X$ in a non-empty
collection of curves, at least one of which is essential in $\bdy
X$. Then the essential curves of intersection have slope a
meridian on $\bdy X$.
\end{claim}

We have a 2--sphere $S$ which intersects $\bdy X$ in a non-empty
collection of curves, at least one of which is essential in $\bdy
X$. There is no loss in generality to assume that among all
2--spheres meeting $\bdy X$ in the same slope as $S$,
 $S$ has the minimal number of curves that are inessential
in $\bdy X$.  Let $c$ be a curve of intersection which is
innermost on the sphere $S$. If $c$ is inessential on $\bdy X$,
then we may perform an isotopy of $S$ that removes $c$ (and
perhaps some other inessential curves of intersection), a
contradiction.  We conclude that $c$ is an essential curve in
$\bdy X$ bounding an embedded disk whose interior is disjoint from
$\bdy X$.  But, $X$ is the exterior of a non-trivial knot in
$S^3$, so $c$ must be the meridional slope on $\bdy X$.  Any other
essential curve in the intersection is parallel to $c$, and is
therefore also a meridian curve.  This completes the proof of the
claim.

Each of the boundary edges of the triangulation is a loop ($\T$ is
a one-vertex triangulation on $\bdy X$) and, therefore, defines a
knot in $S^3$. The {\it bridge number} of a knot $K$ relative to a
height function $h$, is its minimum number of maxima, taken over
all generic embeddings of knots $K'$ that are ambient isotopic to
$K$.

\begin{claim}
If a boundary edge $e$ of $\T$ has bridge number 1, then $e$ is a
meridian or a generalized longitude.
\end{claim}

We assume that the boundary edge $e$ of $\T$ has bridge number
$1$. Then there is an ambient isotopy of $S^{3}$ so that with
respect to the given genus zero Heegaard decomposition, $e$ has
only one maximum and one minimum. We may take this as the original
embedding of $X$ and $\bdy X$. Choose a level surface that
intersects $\bdy X$ in a collection of curves that contains at
least two essential components in $\bdy X$ (a Heegaard surface is
separating).  By the previous claim, each of these essential
curves is a meridian.  Moreover, there are at most 2 intersections
between $e$ and these curves, hence at most 1 intersection between
$e$ and each of these curves. For otherwise, $e$ would necessarily
contain more than one maximum and one minimum. As $e$ intersects
one of these  meridians at most once, it is itself either a
meridian (does not intersect) or a curve that meets a meridian
exactly once, a generalized longitude.  This completes the proof
of the claim.

The three edges of $\T$ in $\bdy X$ meet pairwise in exactly one
point. It is possible that one of these edges is a meridian and
has bridge number 1.  However, there must be at least two edges
which are not bridge number 1.  If two edges are bridge number 1,
then at most one is a meridian, and so one is necessarily a
generalized longitude.  But, the fact that a generalized longitude
has bridge number 1 implies that $X$ is the exterior of an unknot,
a contradiction. It follows that at least two boundary edges have
bridge number at least 2, and in particular possess a maximum that
is not the vertex of the triangulation.

Consider the height function as restricted to the 1-skeleton of
the triangulation, $h:\T^{(1)} \rightarrow [0,1]$.  A {\it thick
region for a set of edges $E$} is a sub-interval $(a,b) \subset
[0,1]$ which consists only of regular values of $h:\T^{(1)}
\rightarrow [0,1]$ and so that $a$ is a critical value
corresponding to a minimum of some edge  $e \in E$ and $b$ is a
critical value corresponding to a maximum of some edge $ e' \in E$
and this maximum is not the vertex.

We may choose $e$, a boundary edge that has bridge number at least
2. There is necessarily a thick region for the edge $e$.  Identify
all of the thick regions for $e$ and within each of these choose a
thick region that is a thick region for all edges of the
triangulation. This yields a collection of thick regions
$\{(a_1,b_1),(a_2,b_2),\dots,(a_n,b_n)\}$.  Within each of the
thick regions $(a_i,b_i)$ we apply the four claims of Thompson,
each of which follows from thin position.

\begin{claim}
For some $t_i \in (a_i,b_i)$ there is a level 2-sphere
$S_i=S_{t_i}$ which intersects the boundary of each tetrahedron in
normal curves or curves disjoint from the 1-skeleton.
\end{claim}

We can assume (see \cite{s-t}) that at the top of the thick region
$(a_i,b_i)$, just below $b_i$, there is a {\it high disk} for the
1-skeleton which is contained in the 2-skeleton.  A high disk is a
boundary compression for $S_{t_i}$ in the exterior of the
1-skeleton that starts above $S_{t_i}$. We may also assume that
there is a {\it low disk} contained in the 2-skeleton at the
bottom of this thick level, just above $a_i$. Thin position
guarantees that for some value of $t_i$ in this thick region there
is a level surface $S_i = S_{t_i}$ for which there is no high or
low disk contained in the 2-skeleton. For otherwise, at some level
between there would be a pair of cancelling high and low disks. In
particular, the intersection of $S_i$ with the boundary of each
tetrahedron does not contain any curves which intersect the
1-skeleton but are not normal.  Such a curve implies an innermost
arc joining an edge to itself which defines a bigon in the
2-skeleton that is either a high or low disk.  This completes the
proof of the claim.

\begin{claim}
$S_i$ does not intersect any tetrahedron $\Delta$ in a normal
curve of length greater than 8.
\end{claim}

If there is a normal curve $c \subset \bdy \Delta$ of length
greater than 8, then this curve must intersect some edge $e$ at
least three times \cite{thompson}.  Following $e$ through three
consecutive intersections with $c$ we note that they cobound two
bigons on $\bdy \Delta$, one above $c$ and one below $c$. Moreover
these bigons may be chosen to be disjoint except for a single
point of intersection on $c$.  They may contain portions of other
edges (including $e$), but, by pushing them slightly into
$\Delta$, see Figure \ref{f-long8}, they become a cancelling pair
of high and low disks for the 1-skeleton.  In particular, they can
be used to guide an isotopy of $e$ that reduces the width of the
1-skeleton. (It is possible that there are other curves of
intersection $c'$ that also intersect the portion of $e$ that
bounds the bigons.  In this case the isotopy is even more
beneficial in reducing width.)  This completes the proof of the
claim.

\begin{figure}[h]
{\epsfxsize = 2.75 in \centerline{\epsfbox{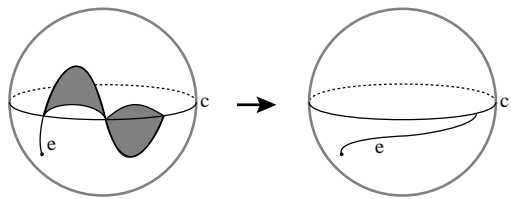}} }
\caption{A normal curve of length greater than 8.} \label{f-long8}
\end{figure}

\begin{claim}
The sphere $S_i$ does not intersect any tetrahedron $\Delta$ in
parallel curves of length 8.
\end{claim}

A normal curve of length 8 on $\bdy \Delta$ intersects two
distinct edges twice \cite{thompson}.  If $c$ and $c'$ are an
outermost pair of parallel curves of length 8 then some edge $e$
hits each twice and there are bigons bounded by both $c$ and $e$
and $c'$ and $e$.  As $c$ and $c'$ are parallel and outermost, one
bigon is a subdisk of the other and when the larger one is pushed
slightly into $\Delta$, it acts simultaneously as a high disk and
low disk that can be used to reduce the width.  This completes the
proof of the claim.


\begin{claim}
The sphere $S_i$ does not intersect distinct tetrahedra, $\Delta$
and $\Delta'$, in curves of length 8.
\end{claim}

If $c$ is a curve of length 8 in $\bdy \Delta$, it intersects two
edges $e$ and $e'$ exactly twice.  This defines two bigons, when
pushed into $\Delta$ one is a high disk for $e$ and the other a
low disk for $e'$.  These disks are not disjoint when pushed into
$\Delta$ and do not by themselves contradict thin position.
However, we have the same situation in $\Delta'$ so we may choose
a high disk in $\Delta$ and a low disk in $\Delta'$ which reduce
width and contradict thin position.  This completes the proof of
the claim.

\begin{claim}
For some $i$ the  intersection $S_i \cap \bdy X$ contains a
meridional curve in $\bdy X$.
\end{claim}

Consider the collection of level surfaces $\{S_1,S_2,\dots,S_n \}$
one chosen for each of the thick regions
$\{(a_1,b_1),(a_2,b_2),\dots,(a_n,b_n)\}$.  The surface $S_i$ was
chosen within a thick region for the boundary edge $e$ and
necessarily intersects $e$. By the first claim, if for some $i$,
the intersection $S_i \cap \bdy X$ contains an essential curve
then that curve is meridional and we are done. The alternative is
that each of these intersections $S_i \cap \bdy X$ consists
entirely of trivial curves, the normal ones are vertex linking and
the others disjoint from the boundary edges.  Choose the outermost
vertex linking curve $c$. The curve $c$ bounds a disk $D$ in $\bdy
X$. The boundary edge $e$ intersects $D$ in two arcs that are
joined to the vertex, call the union of these arcs $e'$. The
remainder of $e$ is a single arc in $\bdy X - D$ which is
connected to the endpoints of $e'$, call this arc $e''$.

Now $e''$ can only possess a single maximum or minimum. For
otherwise there would be a thick region for $e$ between some
maximum and minimum of $e''$, and we have chosen  thick regions
$(a_i,b_i)$ and a level surface $S_i$ within each such thick
region. The level surface $S_i$ would intersect the interior of
the edge $e''$.

Then $e'$ is parallel in $D$ to a subarc of $c = \bdy D$. Perform
this isotopy, see Figure \ref{f-bridge}.  Then the boundary edge
$e$ can be isotoped so that it has only a single minimum and
maximum.  This contradicts our choice of an edge $e$ with bridge
number at least 2, and we conclude that there must be a curve of
intersection that is essential in $\bdy X$ and by the above it
must be meridional.  This completes the proof of the claim.

\begin{figure}[h]
{\epsfxsize = 5.5 in \centerline{\epsfbox{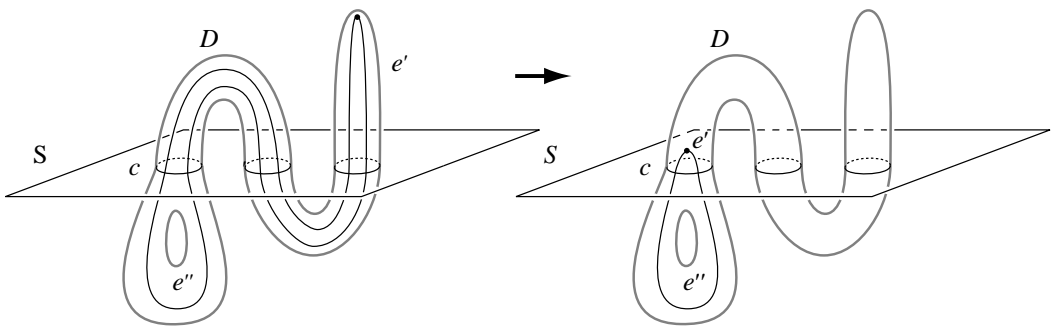}} }
\caption{When $S \cap \bdy X$ consists of trivial curves.}
\label{f-bridge}
\end{figure}

Let $S$ be one of the level spheres $S_i$ which possesses a
meridional curve of intersection with $\bdy X$.   The arguments
above guarantee that $S$ intersects the 2-skeleton of $\T$ in
normal curves and curves disjoint from the 1-skeleton.  However,
the intersection of $S$ with each tetrahedron of $\T$ may not
consist entirely of disks, there may be planar surfaces with more
than one boundary component (tubes).  So, within each tetrahedron
$\Delta$ compress $S$ to a collection of disks. Throw away any
component (a disk or a sphere) which does not intersect the
1-skeleton.  Any component of the resulting surface is a normal or
almost normal sphere or planar surface in $(X,\T)$. At least one
of these components $S'$ is planar and has non-empty boundary
containing at least two meridional curves.
\end{proof}

If $X$ is not the solid torus, then  Lemma \ref{knot-normal} and
Corollary \ref{slopes-are-finite} imply that there is a  finite
computable set of slopes $\alpha$ so that $X(\alpha)$ is $S^3$. In
fact, by the work of Gordon and Luecke \cite{g-l} there is at most
one filling on $X$ which can produce $S^3$. We complement this
result by giving an algorithm that either computes this slope or
demonstrates that it does not exist.

\begin{thm} \label{s3algo}
Given a knot-manifold $X$, there is an algorithm to determine
precisely those slopes $\alpha$ for which $X(\alpha)$ is $S^3$. In
particular, this gives an algorithm to determine whether $X$
embeds in $S^3$.
\end{thm}

\begin{proof}
We are given $X$ via a triangulation, we may assume that this is a
one-vertex triangulation $\T$.

First, we determine whether $X$ has compressible boundary
\cite{haken,j-o} and/or whether $X$ is irreducible (Theorem
\ref{T-recognize-2}). If $X$ has compressible boundary and is
irreducible, then $X$ is the exterior of the unknot in $S^3$. Dehn
filling along any slope $\alpha$ which intersects the slope of the
compressing disk once will produce $S^3$. This set is computable
as the slope of the compressing disk is a by-product of these
computations.  If $X$ is reducible, then $X$ is not the exterior
of a knot in $S^3$, no filling produces $S^3$.

We may therefore assume that $X$ is irreducible and not the
exterior of an unknot ($\bdy X$ is incompressible). If for any
$\alpha$ the manifold $X(\alpha)$ is $S^3$, then by Lemma
\ref{knot-normal}, $X$ possesses a normal or almost normal surface
with slope $\alpha$.  By Corollary \ref{slopes-at-vertices} the
slope $\alpha$ is the slope of a normal or almost normal surface
(using a single octagon) whose projective class is at a vertex of
$\P(X,\T)$.  There are only finitely many such slopes.

For each slope $\alpha$ bounding an embedded normal or almost
normal surface whose projective class is at a vertex of
$\P(X,\T)$, use the filling described in Section \ref{s-layered}
to construct the manifold $X(\alpha)$.  Use the 3-sphere
recognition, Theorem \ref{r-3-sphere} (see
\cite{rubinstein,thompson}) to determine whether $X(\alpha)$ is
$S^3$.  If any $X(\alpha)$ is $S^3$ then this is the sole filling
producing $S^3$ \cite{g-l}. If after checking all of this finite
number of fillings, none is $S^3$, then no Dehn filling gives the
3--sphere and $X$ does not embed in $S^3$.
\end{proof}

\noindent {\bf Algorithm $\S$.} {\it Given a knot-manifold $X$,
determine precisely those slopes $\alpha$ for which $X(\alpha)$ is
the 3--sphere.}

\vspace{.2 in}

{\bf Step 1.} Endow $X$ with a one-vertex triangulation and
compute the vertices of $\P(X,\T)$.  (Using both normal surfaces
{\it and} almost normal surfaces with only octagons.)

{\bf Step 2.}  Determine whether $X$ is reducible.  If so, no
filling can produce the 3--sphere and the algorithm terminates.

{\bf Step 3.} Determine whether $X$ has compressible boundary,
i.e., whether there is a normal disk $D$ with its boundary an
essential curve $\mu$ in $\bdy X$.  If so, then $X$ is a solid
torus and Dehn filling along any slope on the line $L_\mu$
produces the 3--sphere.

{\bf Step 4.}  List the slopes $\{\alpha_1, \dots, \alpha_n\}$
that correspond to embedded vertex surfaces of $\P(X,\T)$.  For
each $\alpha_i$ construct $X(\alpha_i)$ via a layered
triangulation and determine whether it is the 3--sphere using the
algorithm given by Theorem \ref{r-3-sphere}. If any such filling
is found, terminate the algorithm, it is the only filling
producing the 3--sphere \cite{g-l}.

The slopes from Step 3 or slope from Step 4 are the only Dehn
fillings yielding the 3--sphere.

\vspace{.25 in}

Suppose $X$ is the exterior of a knot $K$ in a lens space. If $K$
is isotopic into a Heegaard torus and $X$ has incompressible
boundary then we say that $K$ is a {\it generalized torus knot}
(in a lens space). When the knot-manifold $X$ has compressible
boundary or the exterior of a generalized torus knot we get
special cases for lens space fillings. Note that the exterior of a
generalized torus knot is the union of two solid tori glued along
an incompressible and $\partial$-incompressible annulus $A$, i.e.
a Seifert fibered space over the disk with 2 exceptional fibers.
If $\alpha$ is the slope of the annulus $A$ on $\bdy X$ and $\beta
\in L_\alpha$ then $X(\beta)$ will possess a genus one Heegaard
splitting, i.e. is either $S^3$ or a lens space. Thus, $X$
possesses an infinite number of slopes yielding lens spaces (at
most one is $S^3$).  We must be able to recognize this situation.

\begin{lem}
\label{T-torusknot} Let $X$ be a knot-manifold. There is an
algorithm to determine whether $X$ is a generalized torus knot
exterior.
\end{lem}

\begin{proof}
We may assume $X$ is given via a triangulation $\T$ that restricts
to a one-vertex triangulation on $\bdy X$.

Recall that the annulus $A$ characterizing a generalized torus
knot exterior (an embedded, essential annulus separating the
manifold into two solid tori) is vertical (composed entirely of
regular fibers) with respect to the Seifert fibering of the
manifold.  Moreover, it is the unique essential annulus with
boundary a regular fiber.

Hence, if $X$ is a generalized torus knot exterior it contains an
essential annulus of the above type and if $\T$ is a triangulation
of $X$ that restricts to a one-vertex triangulation on $\bdy X$,
then there is a normal annulus having these same properties and
whose projective class is a vertex in $\P(X,\T)$. To see this let
$A$ be such an annulus and suppose that some multiple of $A$ can
be written as a sum $$k A = \sum_i k_i V_i,$$ where each $V_i$ is
an incompressible and $\bdy$-incompressible vertex surface
(Theorem \ref{T-esscarrier}) and $\chi(V_i) \leq 0$ ($X$ has
incompressible boundary and no summand can be a 2--sphere or
$\mathbb R P^2$). Then $\chi(V_i) =0, \forall i$ and some $V_i$,
say $V_1$, has non-empty boundary with the same slope as $A$
(Proposition \ref{same-or-comp}). So $V_1$ is either an annulus or
a M\"obius band and either $V_1$ or $2 V_1$, respectively, is an
essential annulus with the same boundary slope. This implies that
$V_1$ or $2 V_1$ is isotopic to $A$. Hence, there is such an
annulus that has its projective class a vertex of $\P(X,\T)$.

Now, to determine whether $X$ is a generalized torus knot
exterior, first enumerate the vertices of $\P(X,\T)$ that
correspond to separating annuli.  For each of these annuli $A$,
split $X$ at $A$, retriangulate the components, and determine
whether each is a solid torus (has compressible boundary and is
irreducible).  The knot-manifold $X$ is a generalized torus knot
exterior if and only if we find such a decomposition.
\end{proof}

\begin{lem} \label{lens-space-slope}
Suppose $X$ is a knot-manifold with incompressible boundary which
is not a generalized torus knot exterior, and that for some slope
$\alpha$ the Dehn filling $X(\alpha)$ is a lens space. For any
one-vertex triangulation $\T$ of $X$, either
\begin{enumerate}
\item $(X,\T)$
contains a normal or almost normal surface (a punctured sphere or
torus) with slope $\alpha$, or
\item $\alpha$ is the slope of an edge of
the triangulation $\T$ in $\bdy X$.
\end{enumerate}
\end{lem}

\begin{rem}
For the purposes of this lemma an almost normal surface possesses
a single octagon (no tubes).
\end{rem}

\begin{proof}
This lemma is an adaptation of Lemma \ref{knot-normal}, above,
which was the case for non-trivial knots in $S^3$. We have that
$X(\alpha)$ is a lens space and we proceed as before, putting the
1-skeleton of $\T$ in thin position. This time using a foliation
of the lens space by level Heegaard tori $H_t$.

The first adjustment is a  variation of the first claim of Lemma
\ref{knot-normal} for Heegaard tori in lens spaces.

\begin{claim}
Suppose that a Heegaard torus $H$ intersects $\bdy X$ in a
non-empty collection of curves, at least one of which is essential
in $\bdy X$. Then that curve has meridional slope  on $\bdy X$.
\end{claim}

We have a Heegaard torus $H$ which intersects $\bdy X$ in a
non-empty collection of curves, at least one of which is essential
in $\bdy X$. There is no loss in generality to assume that among
all Heegaard tori meeting $\bdy X$ in the same slope as $H$, that
$H$, itself, has the minimal number of curves. Suppose $c$ is a
curve of intersection between  $H$ and $\bdy X$ which is
inessential and innermost on the torus $H$. If $c$ is inessential
on $\bdy X$, then we may perform an isotopy of $H$ that removes
$c$ (and perhaps some other inessential curves of intersection), a
contradiction.  So $c$ is an essential curve in $\bdy X$ bounding
an embedded disk whose interior is disjoint from $\bdy X$.  But
$X$ has incompressible boundary, so $c$ must be the meridional
slope on $\bdy X$.

The alternative is that there is no curve $c$ which is inessential
in the Heegaard torus $H$; hence, $H$ is cut into a collection of
annuli by its intersection with $\bdy X$.  If any intersection
curve is inessential in $\bdy X$ then at least one of these
annuli, call it $A$, joins an essential curve in $\bdy X$ to an
inessential curve in $\bdy X$.  This also shows that the
intersection is meridional; perform a surgery on the annulus at
the inessential end to produce a disk bounding the essential
curve.

We are left in the case that every curve of intersection is
essential in both $H$ and $\bdy X$; thus cutting each into a
collection of annuli.  If any annulus is compressible in one of
the solid tori bounded by $H$ then $\bdy X$ is compressible or the
slope is meridional.    We are left assuming that each annulus is
boundary parallel in the solid tori bounded by $H$.  We can reduce
the number of intersections (a contradiction) by pushing an
outermost annulus out of one solid torus and into the other unless
the surfaces intersect in exactly 2 curves, cutting each surface
into 2 annuli. Each of the annuli from $\bdy X$ are then isotopic
to one of the annuli in $H$.  If the two annuli are isotopic to
distinct annuli, then $\bdy X$ is isotopic to $H$, a
contradiction, $X$ is not a solid torus. So they are both isotopic
to the same annulus. This implies that the core of the attached
solid torus is isotopic into the Heegaard torus and is either a
generalized torus knot or $X$ has compressible boundary, a
contradiction.  This completes the proof of the claim.

The second claim follows exactly as before, an edge with bridge
number 1 is either a meridian or a generalized longitude.

We now need to show that there is a boundary edge that has bridge
number at least 2. If all three edges have bridge number 1, then
two are generalized longitudes and the other a meridian.  At this
point, there is a notable difference with the $S^3$ case.  In a
lens space it is distinctly possible for the longitude of a knot,
hence the knot itself, to have bridge number 1, yet not be trivial
(the boundary of its exterior is not compressible). In this case,
we have the second conclusion of the theorem:  one of the boundary
edges is the slope of the meridian.

With this exception noted, we continue as before.  Choose a
boundary edge with bridge number at least 2, and identify its
thick regions. Within each of these regions we choose a thick
region for all edges of the triangulation.  This produces a list
of thick regions $\{(a_1,b_1),(a_2,b_2),$ $\dots,(a_n,b_n)\}$. In
each thick region $(a_i,b_i)$ a level Heegaard torus $H_i$ is
found which intersects the boundary of each tetrahedron in normal
curves and curves disjoint from the 1-skeleton. Furthermore, no
normal curve is longer than 8, and there is at most one of length
8.  One of these level surfaces $H=H_i$ must intersect $\bdy X$ in
a curve that is essential in $\bdy X$. These curves are meridional
by the first claim. We compress $H \cap X$ inside each tetrahedron
and choose a normal or almost normal component $H' \subset (X,\T)$
with meridional slope.  Compressing may have lowered genus so $H'$
is either a punctured torus or a punctured sphere.
\end{proof}

If the knot-manifold $X$ has incompressible boundary and is not a
generalized torus knot then Lemma \ref{lens-space-slope} and
Corollary \ref{slopes-are-finite} imply that there is a finite
computable set of slopes $\alpha$ so that $X(\alpha)$ is a lens
space.  In fact, it is well known \cite{cgls} that any pair of
such slopes have distance at most 1, for a total of at most 3
slopes. Again, we complement this result by supplying an algorithm
which either determines precisely these slopes or demonstrates
that they do not exist.

\begin{thm}
Given a knot-manifold $X$ there is an algorithm to determine
precisely those slopes $\alpha$ for which the Dehn filling
$X(\alpha)$ is a lens space.
\end{thm}

\begin{proof}
We may assume $X$ is given via a one-vertex triangulation $\T$.

First, one determines whether $X$ has compressible boundary, in
which case a Dehn filling on $X$ will produce a lens space only if
$X$ is a trivial knot in $S^3$ or in a lens space.  If $X$ does
have compressible boundary then we can find a normal disk $D$ with
essential boundary in $\bdy X$, call its slope $\mu$.  Cut $X$
along $D$ and cap off the resulting 2--sphere boundary with a ball
to obtain a closed manifold $\widehat{X_D}$.  Next, using the
algorithm from \cite{rubinstein}, Theorem \ref{r-3-sphere} above,
determine whether $\widehat{X_D}$ is the 3--sphere, if so, $X$ is
a solid torus and every filling on $L_\alpha$ produces the
3--sphere and every other filling produces a lens space. If
$\widehat{X_D}$ is not the 3--sphere, using the algorithm from
\cite{rubinstein}, Theorem \ref{r-lensspace} above, determine
whether it is a lens space. If so, then $X$ is the exterior of a
trivial knot in a lens space and $X(\beta)$ is a lens space
precisely for $\beta \in L_\mu$. If not, then no filling can
produce a lens space.

If $X$ is the exterior of a generalized torus knot, then we may
determine so by Lemma \ref{T-torusknot}. Moreover, that algorithm
will produce the slope $\alpha$ of the essential annulus. In this
case $X(\beta)$ is a lens space or $S^3$ for precisely the slopes
$\beta \in L_\alpha$. By Theorem \ref{s3algo} we may identify
which slope, if any, to fill along to obtain $S^3$.

The remaining case is that $X$ has incompressible boundary and is
not the exterior of a generalized torus knot.  If for any $\alpha$
the manifold $X(\alpha)$ is a lens space, then by Lemma
\ref{lens-space-slope}, $X$ possesses a normal or almost normal
surface (using a single octagon) with slope $\alpha$ or $\alpha$
is the slope of a boundary edge. In the former case, Corollary
\ref{slopes-at-vertices} implies that $\alpha$ is the slope of a
normal or almost normal surface whose projective class is at a
vertex of $\P(X,\T)$.  The lens space fillings can then be
identified by performing the following steps.  For each slope
$\alpha$ which is either the slope of an embedded vertex normal or
almost normal surface or the slope of one of the three boundary
edges, use the filling described in Section \ref{s-layered} to
construct the manifold $X(\alpha)$ and use the lens space
recognition Theorem \cite{rubinstein} to determine whether
$X(\alpha)$ is a lens space.
\end{proof}

\noindent {\bf Algorithm L.} {\it Given a knot-manifold $X$,
determine precisely those slopes $\alpha$ for which $X(\alpha)$ is
a lens space.}

\vspace{.2 in}

{\bf Step 1.} Endow $X$ with a one-vertex triangulation and
compute the vertex solutions of $\P(X,\T)$.  (Again, we are
considering both normal and almost normal surfaces with octagons.)

{\bf Step 2.} If a vertex solution of $\P(X,\T)$ is a disk $D$
whose boundary is an essential curve in $\bdy X$, then $\bdy X$ is
compressible.  Determine the slope of $\bdy D$, say $\mu$.

In this case, cut $X$ along $D$ and cap off the remaining
2--sphere boundary component with a ball. This yields a closed
manifold $\widehat{X_D}$. Determine whether $\widehat{X_D}$ is the
3--sphere or a lens space. If $\widehat{X_D}$ is the 3--sphere,
then filling along every slope $\beta \in L_\mu$ yields $S^3$ and
every other filling yields a lens space.  If $\widehat{X_D}$ is a
lens space, then filling along every slope $\beta \in L_\mu$
yields that same lens space.  In no other case is a lens space
filling obtained; and the algorithm terminates.

{\bf Step 3.}  For each vertex solution of $\P(X,\T)$ that is a
separating annulus, split $X$ along the annulus and determine if
each component is a solid torus; i.e. if $X$ is a generalized
torus knot exterior.  If an annulus, say $A$, is found so that $X$
split at $A$ yields two solid tori, then compute $\alpha$, the
boundary slope of $A$ (there is a unique such boundary slope).
Then $X(\beta)$ is a lens space or $S^3$ for precisely those
$\beta \in L_\alpha$.   Algorithm $\S$ can identify which slope,
if any, produces $S^3$; and, the algorithm terminates.

If $X$ is not a generalized torus knot exterior (and $\bdy X$ is
incompressible), go to the next step.

{\bf Step 4.}  Enumerate the slopes $\{\alpha_1, \dots,
\alpha_n\}$ of vertex normal and almost normal surfaces.  For each
slope $\alpha_i$ construct $X(\alpha_i)$ via a layered
triangulation of a solid torus and determine whether it is a lens
space. If at any time three such slopes are found, terminate the
algorithm, this is the maximum number of slopes yielding lens
space fillings \cite{cgls}. This completes the algorithm.

\section{Summary comments}

The preceding considerations are well adapted to normal (and
almost normal) surface theory.  In each, our algorithms were based
on finding interesting surfaces and are rather comprehensive in
their application to exceptional and Haken Dehn fillings. However,
there are some notable exclusions.  We have not considered Dehn
fillings that are Seifert fibered or have finite fundamental group
(except for $S^3$ and lens spaces).

Our methods can be used to determine for a given knot-manifold $X$
those Dehn fillings that are Haken-manifolds and are Seifert
fibered.  The proof uses the methods of Section \ref{s-essential}
and, while quite tedious, does not require new ideas.  However,
there is a major gap for applying our methods to determine small
Seifert fibered manifolds, Seifert fibered manifolds that are not
Haken-manifolds.  The major problems here are probably associated
to the lack of understanding of immersed (not embedded) normal
surfaces.  We give the following remaining open problems.

\begin{prob}
Given a 3--manifold $M$ that is known to be irreducible and not a
Haken-manifold, is there an algorithm to determine if $M$ is a
small Seifert fibered space ?
\end{prob}

\begin{prob}
Given a knot-manifold $X$ is there an algorithm to determine
precisely those slopes $\alpha$ for which the Dehn filling
$X(\alpha)$ is Seifert fibered ?
\end{prob}

While similar, the next problem probably calls for even a wider
range of new ideas.

\begin{prob}
Given a knot-manifold $X$, is there an algorithm to determine
precisely those slopes $\alpha$ for which the Dehn filling
$X(\alpha)$ has finite fundamental group ?
\end{prob}

Finally, our objective has been to determine interesting phenomena
in the space of Dehn fillings on a given knot-manifold $X$. One
very interesting open problem is the homeomorphism problem for
manifolds in the space of Dehn fillings on $X$.

\begin{prob}
Given the knot-manifold $X$ and slopes $\alpha$ and $\beta$ is
there an algorithm to determine if $X(\alpha)$ and $X(\beta)$ are
homeomorphic ?
\end{prob}

\bibliographystyle{plain}
\bibliography{dd}

\end{document}